\documentclass[11pt]{article}%
\usepackage{amssymb}
\usepackage{amsfonts}
\usepackage{amsmath}
\usepackage{graphicx}%
\providecommand{\U}[1]{\protect\rule{.1in}{.1in}}
\begin{document}

\title{Voltage and current density in a rectangle with point in- and ejection via
Green's functions and $q$-analysis }
\author{P. Van Mieghem\thanks{ Faculty of Electrical Engineering, Mathematics and
Computer Science, P.O Box 5031, 2600 GA Delft, The Netherlands; \emph{email}:
P.F.A.VanMieghem@tudelft.nl }}
\date{Delft University of Technology\\
25 June 2026}
\maketitle

\begin{abstract}
The current density in a rectangle in the $D=2$ dimensional Euclidean space,
where opposite point charges are placed in a source $s$ and a destination $d$,
is computed as the gradient of a potential that obeys the Poisson equation
(first law of Maxwell). Although the computation of the current density in a
2D rectangle is a classical problem, several barriers (e.g. very slow
convergence of double Fourier series) were eventually alleviated by a
$q$-analysis, related to Gaussian polynomials and Jacobi's theta-functions,
that form the basic building blocks for any elliptic function. We present an
exact and computationally very efficient analytic formula of the potential and
the magnitude of the current density in a rectangle.

\end{abstract}

\section{Introduction}

We compute the potential and current density in each point of a rectangle in
the two-dimensional Euclidean space, where at the source point $s$ a unit
current is injected, which leaves the rectangle at the destination point $d$.
The electric unit current thus spreads from the injection point $s$ to the
ejection point $d$ over the entire rectangle, steered by the first law of
Maxwell, briefly reviewed in Section \ref{sec_flow_electrical_current}.
Although the potential and current computation is a classical problem, that is
usually solved by Fourier series, our main result is the Green's function
(\ref{quasi_Green_in_Tq_elliptic}) in terms of a $q$-function
(\ref{def_Tq_function}), which belongs to the class of elliptic functions. A
famous earlier victory in physics is the \emph{phenomenal} paper
\cite{Onsager1944} of Noble-Prize winner Ralph Onsager, in which he solved the
2D Ising model \emph{exactly} using elliptic functions (Jacobi theta
functions, Jacobian elliptic functions and elliptic integrals \cite[from p.
144 on]{Onsager1944}).

Our motivation arose from Kitsak \emph{et al}. \cite{Kitsak2023_sp}, who
conjectured that the shortest path $\mathcal{P}_{sd}^{\ast}$ in many graphs,
where only partial information of the graph and its link weights is available
as in most real-world networks, is close to the geodesic from $s$ to $d$. The
knowledge of the set of network nodes that most likely lie on the shortest
path between a source and destination is important, both for the network
operator to secure and for the adversary to eliminate or disrupt those nodes.
A simple example of a partially known graph is the hard random geometric
graphs (HRGG), in which $N$ nodes are uniformly distributed in a 2D rectangle
and connected to other nodes that lie within a given distance $r_{0}$.
Randomness here reflects our uncertainty about the precise nodal positions in
a partially observed network. Random geometric graphs (RGGs) have been studied
extensively (see e.g. references in \cite{PVM_geometric_random_graphs_Frechet}%
,\cite{Kitsak2023_sp}). RGGs can model transportation networks such as
wireless and airline networks as well as infrastructural networks like power
grids. RGGs can be applied to analyze the structure of large data sets and in
modelling ad hoc networks (such as vehicular, disaster relief, sensor and
flying swarm robotics networks). The probability that a point $\left(
x,y\right)  $ in the rectangle belongs to the shortest path $\mathcal{P}%
_{sd}^{\ast}$ between a given pair $\left(  s,d\right)  $ of nodes in the HRGG
is difficult to compute analytically, while the current density and the
electric potential in a continuous space can be computed as shown below. In
the high-density limit ($N\rightarrow\infty$) of an HRGG with link
connectivity distance $r_{0}=O\left(  \frac{1}{\sqrt{N}}\right)  $, we
conjecture that the magnitude of current density $\left\vert j\left(
x,y\right)  \right\vert $ in each point $\left(  x,y\right)  $ of the
rectangle, due to current in- and ejection at $s$ and $d$, is proportional to
the probability that the point $\left(  x,y\right)  $ lies on the shortest
path $\mathcal{P}_{sd}^{\ast}$ in an HRGG between source $s$ and drain $d$.

The paper presents an exploratory journey from Fourier series towards
$q$-functions, belonging to the class of elliptic functions. The paper starts
in Section \ref{sec_Laplacian_in_3D} by computing the Laplacian equation
(\ref{Laplacian_differential_eq}) in 3D for the current in- and ejection in
the beam in Fig. \ref{fig_balk_afmetingen}. Although the computation is rather
involved, but a standard exercise in solving the Laplacian equation in
rectangular or Cartesian coordinates, the importance is that an exact solution
of the potential\footnote{The potential $V\left(  x,y,z\right)  $ is only
unique if a reference potential $V_{\text{ref}}$ is chosen.} $V\left(
x,y,z\right)  +V_{\text{ref}}$ exists, which is written as a double Fourier
series in (\ref{Voltage_beam}). Section \ref{sec_reduction_2D} considers the
limit where the size $c$ of the parallellepipidum in Fig.
\ref{fig_balk_afmetingen} tends to zero, thus forming a rectangle in the
$xy$-plane with sizes $a$ and $b$. Here, we encounter physical limitations,
because the injection current is orthogonal to the $xy$-plane and thus, there
are no vector components in the $x$ and $y$-directions, which complicates the
limit $c\rightarrow0$ in our setting! Nevertheless, we show that we can
continue due to the mathematical beauty of a renormalization argument and
Green's function theory. First, we present a double Fourier series for the
potential or voltage $V\left(  x,y,0\right)  $ in (\ref{Voltage_2D}) and show
agreement with the corresponding Green's function. However, the double Fourier
series of the current density converges very slowly, which motivated the
computation of a single Fourier series, presented in Section
\ref{sec_Green_2D_single_sum}. Although the single Fourier is computationally
far superior than the double Fourier series, inaccuracies at $x$- and
$y$-lines passing through injection point $s$ and ejection point $d$,
explained in Section \ref{sec_discontinuities}, have motivated the analytic
evaluation of the single Fourier series (\ref{Green_K_single_sum}) and
(\ref{Green_K_single_sum_dual}). Surprisingly, perhaps, our analysis leads to
the product $\left(  a;q\right)  _{n}=\prod_{k=0}^{n-1}\left(  1-aq^{k}%
\right)  $, which is a basic product in $q$-analysis, briefly reviewed in
Appendix \ref{sec_Gaussianpoly} on Gaussian polynomials. The product $\left(
a;q\right)  _{n}$ also appears in Jacobi's theta functions. Apart from this
connection to elliptic functions, the resulting series for the current density
converges extremely fast, at least in a square $a=b$, to the extent that only
1 term seems\footnote{There is no visual difference between the current
density (computed in Fig. \ref{Fig_currentdensity_fourier_qserie}) with only 1
term compared to the series with more than one term, due to an expansion in
powers of $e^{-2\pi}\approx\frac{1}{535}$. The same amazing convergence
appears in Jacobi theta series.} accurate enough in plots.

In the many Appendices, we have positioned important knowledge, that wanders
around the direct path from the Fourier analysis towards the $q$-analysis. A
brief review of Green's function theory (Appendix \ref{sec_Green_functions}),
of $q$-analysis (Appendix \ref{sec_function_w(z,A)}) and of Gaussian
polynomials (Appendix \ref{sec_Gaussianpoly}) as well as tedious computations
are placed in the Appendices. The effective resistance matrix of a graph,
containing the effective resistance $\omega_{ij}$ between all pairs $\left(
i,j\right)  $ of nodes, is an important graph matrix, as illustrated in
\cite[Chapter 5]{PVM_graphspectra_second_edition}, whose theory is built upon
Kirchoff's circuit laws and the law of Ohm and which is intimately related to
the simplex geometry \cite{PVM_SimplexGeometry} of any graph. The effective
resistance $\omega_{ij}$ lower bounds the weight \cite[Sec. 5.7]%
{PVM_graphspectra_second_edition} (i.e. hopcount, the number of links in a
path, in an unweighted graph) of a shortest path $\mathcal{P}_{ij}^{\ast}$.
The rather straightforward applications of the theory, with explicit formulae
for the voltage in 2D and 3D, is shown to lead to surprises in Appendix
\ref{sec_effective_resistance_2D_3D}.

\section{The flow of electrical current}

\label{sec_flow_electrical_current}The first law of Maxwell
\cite{Feynman_lec_phys2} is Gauss's law,%
\[
\oint_{A}\overrightarrow{E}.\overrightarrow{dS}=\frac{1}{\epsilon_{0}}%
\iiint_{v}\rho dr
\]
which states that the electric flux, which is the scalar product of the
electric field $\overrightarrow{E}$ and the orthogonal vector $\overrightarrow
{dS}$ on the surface, through a closed surface $A$ that surrounds a region
with volume $v$ equals the total amount of charges, specified by the charge
density $\rho\left(  x,y,z\right)  $ at each point $\left(  x,y,z\right)  $,
that are contained in that region with volume $v$, divided by the permittivity
$\epsilon_{0}$ in vacuum. The permittivity in vacuum is a constant and equal
to $\epsilon_{0}=8.85\;10^{-12}C^{2}/Nm^{2}$. The differential form of
Maxwell's first law is%
\[
\operatorname{div}\overrightarrow{E}=\frac{1}{\epsilon_{0}}\rho
\]
The electric field is defined as $\overrightarrow{E}=-\operatorname{grad}V$,
where the potential $V$ is specified by the Poisson equation,%
\[
\Delta V=-\frac{1}{\epsilon_{0}}\rho
\]
that follows from Maxwell's first law $\operatorname{div}\overrightarrow
{E}=\frac{1}{\epsilon_{0}}\rho$ as $\operatorname{div}\operatorname{grad}%
V=-\frac{1}{\epsilon_{0}}\rho$, where%
\[
\operatorname{div}\operatorname{grad}V=\Delta V=\frac{\partial^{2}V}{\partial
x^{2}}+\frac{\partial^{2}V}{\partial y^{2}}+\frac{\partial^{2}V}{\partial
z^{2}}%
\]
and the potential $V\left(  x,y,z\right)  $ is expressed in Cartesian
coordinates. The differential operator $\Delta=\frac{\partial^{2}}{\partial
x^{2}}+\frac{\partial^{2}}{\partial y^{2}}+\frac{\partial^{2}}{\partial z^{2}%
}$ is called the Laplacian, which is a continuous limit of the discrete
Laplacian matrix $Q$ of a graph \cite{PVM_graphspectra_second_edition}.

\begin{figure}
[ptb]
\begin{center}
\includegraphics[
trim=1.650760in 4.400791in 1.651585in 4.401960in,
height=5.3422cm,
width=9.1028cm
]%
{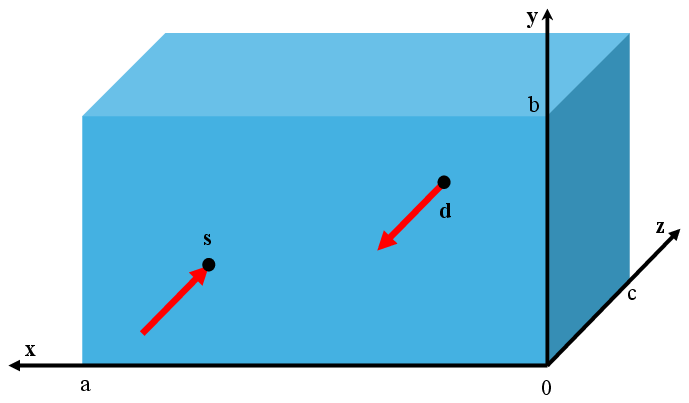}%
\caption{The coordinate axes in a parallellepipidum, in which a current is
injected at the source point $s$ and ejected at the drain point $d$.}%
\label{fig_balk_afmetingen}%
\end{center}
\end{figure}

A beam in Fig. \ref{fig_balk_afmetingen} is a representation in three spacial
dimensions, defined by the Cartesian coordinate frame with $x$-axis, $y$-axis
and $z$-axis. The current injection is orthogonal to the $xy$-plane, i.e.
$z=0$. Fig. \ref{fig_balk_afmetingen}%
%TCIMACRO{\FRAME{ftbpFU}{9.1028cm}{5.3422cm}{0pt}{\Qcb{The coordinate axes in a
%parallellepipidum, in which a current is injected at the source point $s$ and
%ejected at the drain point $d$.}}{\Qlb{fig_balk_afmetingen}}%
%{balk_afmetingen.ps}{\special{ language "Scientific Word";  type "GRAPHIC";
%width 9.1028cm;  height 5.3422cm;  depth 0pt;  original-width 8.2538in;
%original-height 11.6949in;  cropleft "0.2000";  croptop "0.6236";
%cropright "0.7999";  cropbottom "0.3763";
%filename 'Balk_afmetingen.ps';file-properties "XNPEU";}} }%
%BeginExpansion
%EndExpansion
illustrates a beam or parallellepipidum with size $a,b,c$ of a certain
material with constant resistivity $\varrho$. A unit current is injected in a
source point $s$ in the $z=0$-plane at coordinates $\left(  x\,_{in}%
,y_{in},0\right)  $. The current leaves the beam at the drain point $d$ at
$\left(  x_{out},y_{out},0\right)  $. Both injection and ejection point lie at
the surface of the beam. Here, the injection point $s$ and ejection point $d$
lie in the same $xy$-plane. However, the drain point $d$ can also be placed on
any other surface plane of the beam, i.e. at $y=0$ or $y=b$, at $x=0$ or $x=a$
and at $z=0$ or $z=c$.

We compute the potential $V\left(  x,y,z\right)  $ at a point $\left(
x,y,z\right)  $ inside and on the surface of the beam. Since no charges are
contained inside the beam, the Poisson equation reduces to the Laplacian
partial differential equation%
\begin{equation}
\frac{\partial^{2}V}{\partial x^{2}}+\frac{\partial^{2}V}{\partial y^{2}%
}+\frac{\partial^{2}V}{\partial z^{2}}=0 \label{Laplacian_differential_eq}%
\end{equation}
with boundary conditions $\frac{\partial V}{\partial x}=0$ at $x=0$ and $x=a$,
$\frac{\partial V}{\partial y}=0$ at $y=0$ and $y=b$ and $\frac{\partial
V}{\partial z}=0$ at $z=c$, implying that no current\footnote{The current
density is a vector $\overrightarrow{j}\left(  x,y,z\right)  =j_{x}\left(
x,y,z\right)  \overrightarrow{e}_{x}+j_{y}\left(  x,y,z\right)
\overrightarrow{e}_{y}+j_{z}\left(  x,y,z\right)  \overrightarrow{e}_{z}$,
written by the arrow above the symbol and $\overrightarrow{e}_{i}$ is the
basic, unit norm vector in direction $i$, orthogonal to any other direction.
The magnitude $\left\vert \overrightarrow{j}\left(  x,y,z\right)  \right\vert
$ of the current density is a vector norm, for which we employ here the
simplest possible form of a scalar product $\left\vert \overrightarrow
{j}\left(  x,y,z\right)  \right\vert ^{2}=\overrightarrow{j}\left(
x,y,z\right)  .\overrightarrow{j}\left(  x,y,z\right)  $. Sometimes, we omit
the arrow when talking about the magnitude $\left\vert j\left(  x,y,z\right)
\right\vert $ of the current density. The current density through a surface
$A$ in the $z$-direction of the beam is $I=\int_{A}j_{z}\left(  x,y,0\right)
dxdy$. Hence, if $j_{z}\left(  x,y,0\right)  =j$ in a point can be regarded as
constant and the same for all points of the surface, then the current through
that surface is $I=jA$. The resistivity $\varrho=\sigma^{-1}$ has units
[$\Omega m$], i.e. Ohm meter, whereas the conductivity $\sigma$ is expressed
in [$S/m$], i.e. Siemens per meter.} $\overrightarrow{j}\left(  x,y,z\right)
=-\sigma\overrightarrow{E}=-\sigma\operatorname{grad}V$ can flow through the
surface of the beam, except at the plane $z=0$. The current, orthogonal to the
$xy$-plane and equal to $\frac{\partial V}{\partial z}=-\frac{1}{\sigma
}I=-\varrho I$, is injected at the source point $s=\left(  x\,_{in}%
,y_{in},0\right)  $ and that same current $\frac{\partial V}{\partial
z}=\varrho I$ is leaving at the drain point $d=\left(  x_{out},y_{out}%
,0\right)  $ in the other direction. The law of Ohm then states that the
corresponding resistance equals $R_{sd}=\frac{V\left(  x\,_{in},y_{in}%
,0\right)  -V\left(  x_{out},y_{out},0\right)  }{I}$.

\section{Solution of the Laplacian in 3D}

\label{sec_Laplacian_in_3D}

\subsection{Separation of dimensions}

We employ the standard method of separation of dimensions and assume that a
solution of the form $V\left(  x,y,z\right)  =X\left(  x\right)  Y\left(
y\right)  Z\left(  z\right)  $ exists. The Laplacian
(\ref{Laplacian_differential_eq}) then reduces to%
\[
X^{\prime\prime}\left(  x\right)  Y\left(  y\right)  Z\left(  z\right)
+X\left(  x\right)  Y^{\prime\prime}\left(  y\right)  Z\left(  z\right)
+X\left(  x\right)  Y\left(  y\right)  Z^{\prime\prime}\left(  z\right)  =0
\]
where $f^{\prime\prime}\left(  x\right)  =\frac{d^{2}f}{dx^{2}}$.
Dividing\footnote{We can always choose a reference potential in such a way
that $V\left(  x,y,z\right)  \neq0$ at any point $p=\left(  x,y,z\right)  $
inside and on the beam.} by $V\left(  x,y,z\right)  $,%
\begin{equation}
\frac{X^{\prime\prime}\left(  x\right)  }{X\left(  x\right)  }+\frac
{Y^{\prime\prime}\left(  y\right)  }{Y\left(  y\right)  }+\frac{Z^{\prime
\prime}\left(  z\right)  }{Z\left(  z\right)  }=0 \label{Laplacian_separated}%
\end{equation}
and differentiating with respect to $x$, $y$ and $z$ shows that%
\[
\left\{
\begin{array}
[c]{l}%
\frac{d}{dx}\left(  \frac{X^{\prime\prime}\left(  x\right)  }{X\left(
x\right)  }\right)  =0\\
\frac{d}{dy}\left(  \frac{Y^{\prime\prime}\left(  y\right)  }{Y\left(
y\right)  }\right)  =0\\
\frac{d}{dz}\left(  \frac{Z^{\prime\prime}\left(  z\right)  }{Z\left(
z\right)  }\right)  =0
\end{array}
\right.
\]
Integration then gives%
\begin{equation}
\left\{
\begin{array}
[c]{l}%
\frac{X^{\prime\prime}\left(  x\right)  }{X\left(  x\right)  }=\lambda_{x}\\
\frac{Y^{\prime\prime}\left(  y\right)  }{Y\left(  y\right)  }=\lambda_{y}\\
\frac{Z^{\prime\prime}\left(  z\right)  }{Z\left(  z\right)  }=\lambda_{z}%
\end{array}
\right.  \label{diff_eq_3D}%
\end{equation}
while substitution into (\ref{Laplacian_separated}) indicates that the
constants of integration must satisfy%
\begin{equation}
\lambda_{x}+\lambda_{y}+\lambda_{z}=0 \label{eigenvalue_condition}%
\end{equation}
Since the potential $V\left(  x,y,z\right)  $ is a real function in the
spacial coordinates $x$, $y$ and $z$, the integration or separation constants
$\lambda_{x}$, $\lambda_{y}$ and $\lambda_{z}$ in (\ref{diff_eq_3D}) are real.

In summary, the separation of dimensions has reduced the partial differential
equation into three linear second-order differential equations in one variable
in (\ref{diff_eq_3D}) of the form%
\[
W^{\prime\prime}\left(  w\right)  -\lambda_{w}W\left(  w\right)  =0
\]
whose general solution for $\lambda_{w}\neq0$ and constants $\widetilde
{A},\widetilde{B}$ is%
\[
W\left(  w\right)  =\widetilde{A}e^{-\sqrt{\lambda_{w}}w}+\widetilde
{B}e^{\sqrt{\lambda_{w}}w}%
\]
but if $\lambda_{w}=-\left\vert \lambda_{w}\right\vert =i^{2}\left\vert
\lambda_{w}\right\vert $,%
\[
W\left(  w\right)  =A\sin\sqrt{\left\vert \lambda_{w}\right\vert }w+B\cos
\sqrt{\left\vert \lambda_{w}\right\vert }w
\]
If $\lambda_{w}=0$, then $W\left(  w\right)  =\alpha w+\beta$, for real
constants $\alpha$ and $\beta$.

\subsection{Introduction of boundary conditions}

The next step is to incorporate the boundary conditions, $\frac{\partial
V}{\partial x}=0$ at $x=0$ and $x=a$, $\frac{\partial V}{\partial y}=0$ at
$y=0$ and $y=b$ and $\frac{\partial V}{\partial z}=0$ at $z=c$, that translate to%

\begin{equation}
\left\{
\begin{array}
[c]{l}%
X^{\prime}\left(  0\right)  =X^{\prime}\left(  a\right)  =0\\
Y^{\prime}\left(  0\right)  =Y^{\prime}\left(  b\right)  =0\\
Z^{\prime}\left(  c\right)  =0
\end{array}
\right.  \label{boundary_conditions_X_Y_Z}%
\end{equation}
For the $x$-dimension with general solution $X\left(  x\right)  =A_{X}%
\sin\sqrt{\left\vert \lambda_{x}\right\vert }x+B_{X}\cos\sqrt{\left\vert
\lambda_{x}\right\vert }x$ for $\lambda_{x}\neq0$, it holds that%
\[
X^{\prime}\left(  x\right)  =A_{X}\sqrt{\left\vert \lambda_{x}\right\vert
}\cos\sqrt{\left\vert \lambda_{x}\right\vert }x-B_{X}\sqrt{\left\vert
\lambda_{x}\right\vert }\sin\sqrt{\left\vert \lambda_{x}\right\vert }x
\]
Requiring that $X^{\prime}\left(  0\right)  =0$ leads to $X^{\prime}\left(
0\right)  =A_{X}\sqrt{\left\vert \lambda_{x}\right\vert }=0$, while the other
boundary condition $X^{\prime}\left(  a\right)  =0$ implies that $-B_{X}%
\sqrt{\left\vert \lambda_{x}\right\vert }\sin\sqrt{\left\vert \lambda
_{x}\right\vert }a=0$. The condition $\sin\sqrt{\left\vert \lambda
_{x}\right\vert }a=0$ is met, when $\sqrt{\left\vert \lambda_{x}\right\vert
}a=k\pi$, i.e. $\sqrt{\left\vert \lambda_{x}\right\vert }=\frac{k\pi}{a}$ for
$k\in\mathbb{N}_{0}$. If $\lambda_{x}=0$, then the solution is $X\left(
x\right)  =\alpha x+\beta$ and $X^{\prime}\left(  x\right)  =\alpha$. The
boundary condition $X^{\prime}\left(  0\right)  =0$ requires that $\alpha=0$,
in which case $X^{\prime}\left(  x\right)  =0$ for all real $x$ so that the
boundary condition $X^{\prime}\left(  a\right)  =0$ is also satisfied.

Combining all shows that%
\[
X\left(  x\right)  =B_{X}\cos\left(  \frac{k\pi}{a}x\right)
\]
for all non-negative integers $k$, i.e. $k=0,1,2,\ldots$. Similarly, for the
$y$-dimension, we find that $Y\left(  y\right)  =B_{Y}\cos\frac{l\pi}{b}y$.
Also, the integration or separation constants $\left\vert \lambda
_{x}\right\vert =\left(  \frac{k\pi}{a}\right)  ^{2}$ and $\left\vert
\lambda_{y}\right\vert =\left(  \frac{l\pi}{b}\right)  ^{2}$ are non-negative
and both $k,l\in\mathbb{N}$.

The remaining $z$-dimension is a little different. Since $\lambda_{z}%
=-\lambda_{x}-\lambda_{y}$ by (\ref{eigenvalue_condition}) and the boundary
problem requires that $\lambda_{x}=-\left\vert \lambda_{x}\right\vert
=-\left(  \frac{k\pi}{a}\right)  ^{2}$ and $\lambda_{y}=-\left\vert
\lambda_{y}\right\vert =-\left(  \frac{l\pi}{b}\right)  ^{2}$, the
corresponding differential equation $Z^{\prime\prime}\left(  z\right)
-\lambda_{z}Z\left(  z\right)  =0$ in (\ref{diff_eq_3D}) with $\lambda
_{z}=\left(  \frac{k\pi}{a}\right)  ^{2}+\left(  \frac{l\pi}{b}\right)
^{2}\geq0$ has the solution%
\[
Z\left(  z\right)  =\widetilde{A}_{Z}e^{-\sqrt{\lambda_{z}}z}+\widetilde
{B}_{Z}e^{\sqrt{\lambda_{z}}z}%
\]
and%
\[
Z^{\prime}\left(  z\right)  =\sqrt{\lambda_{z}}\left(  -\widetilde{A}%
_{Z}e^{-\sqrt{\lambda_{z}}z}+\widetilde{B}_{Z}e^{\sqrt{\lambda_{z}}z}\right)
\]
The boundary condition $Z^{\prime}\left(  c\right)  =0$ translates to
$\widetilde{B}_{Z}=\widetilde{A}_{Z}e^{-2\sqrt{\lambda_{z}}c}$. Hence,%
\begin{align*}
Z\left(  z\right)   &  =\widetilde{A}_{Z}\left(  e^{-\sqrt{\lambda_{z}}%
z}+e^{\sqrt{\lambda_{z}}\left(  z-2c\right)  }\right)  =\widetilde{A}%
_{Z}e^{-\sqrt{\lambda_{z}}c}\left(  e^{-\sqrt{\lambda_{z}}\left(  z-c\right)
}+e^{\sqrt{\lambda_{z}}\left(  z-c\right)  }\right) \\
&  =2\widetilde{A}_{Z}e^{-\sqrt{\lambda_{z}}c}\cosh\left(  \sqrt{\lambda_{z}%
}\left(  z-c\right)  \right)
\end{align*}
Writing for the unknown coefficient $B_{Z}=2\widetilde{A}_{Z}e^{-\sqrt
{\lambda_{z}}c}$, the solution in $z$-dimension is%
\[
Z\left(  z\right)  =B_{Z}\cosh\left(  \sqrt{\left(  \frac{k\pi}{a}\right)
^{2}+\left(  \frac{l\pi}{b}\right)  ^{2}}\left(  z-c\right)  \right)
\]

In summary, a particular solution for the potential of the Laplace equation
(\ref{Laplacian_differential_eq}), subject to the boundary conditions
(\ref{boundary_conditions_X_Y_Z}), is%
\[
V\left(  x,y,z\right)  =X\left(  x\right)  Y\left(  y\right)  Z\left(
z\right)  =B_{X}B_{Y}B_{Z}\cos\left(  \frac{k\pi}{a}x\right)  \cos\left(
\frac{l\pi}{b}y\right)  \cosh\left(  \sqrt{\left(  \frac{k\pi}{a}\right)
^{2}+\left(  \frac{l\pi}{b}\right)  ^{2}}\left(  z-c\right)  \right)
\]
where the integers $k,l\in\mathbb{N}$ and the coefficients $B_{X}$, $B_{Y}$
and $B_{Z}$ are constants, independent of the spacial coordinates $x,y$ and
$z$.

\subsection{General solution of the Laplacian differential equation}

Since the Laplacian differential equation is a linear equation, any particular
solution for a pair of integers $\left(  k,l\right)  $ is also a solution, as
well as any linear combination of them. The most general solution is then%
\begin{equation}
V\left(  x,y,z\right)  =\sum_{k=0}^{\infty}\sum_{l=0}^{\infty}\varphi_{kl}%
\cos\frac{k\pi}{a}x\cos\frac{l\pi}{b}y\cosh\left(  \sqrt{\left(  \frac{k\pi
}{a}\right)  ^{2}+\left(  \frac{l\pi}{b}\right)  ^{2}}\left(  z-c\right)
\right)  \label{V_general_linear_combination}%
\end{equation}
This general solution of the voltage in $D=3$ dimensions can be extended to a
$D$-dimensional Euclidean space, where the $D$-rectangle has sizes
$a_{1},a_{2},\ldots,a_{D}$, as%
\begin{equation}
V\left(  x_{1},x_{2},\ldots,x_{D}\right)  =\sum_{k_{1}=0}^{\infty}\sum
_{k_{2}=0}^{\infty}\cdots\sum_{k_{D-1}=0}^{\infty}\varphi\left(  k_{1}%
,\ldots,k_{D-1}\right)  \prod_{j=1}^{D-1}\cos\left(  \frac{k_{j}\pi}{a_{j}%
}x_{j}\right)  \cosh\left(  \sqrt{\sum_{j=1}^{D-1}\left(  \frac{k_{j}\pi
}{a_{j}}\right)  ^{2}}\left(  x_{D}-a_{D}\right)  \right)
\label{V_D_general_Euclidean}%
\end{equation}
The coefficient function $\varphi\left(  k_{1},\ldots,k_{D-1}\right)  $ must
be determined from the current injection and -ejection. The boundary
conditions that no current can flow through any face of the $D$-rectangle,
except for one face containing the injection and ejection points, are
satisfied. We proceed with the $D=3$ Euclidean case further.

It remains to determine the coefficients $\varphi_{kl}$ in
(\ref{V_general_linear_combination}) so that the two remaining
current-injection conditions, $\frac{\partial V}{\partial z}=\varrho I$ at the
source point $s=\left(  x\,_{in},y_{in},0\right)  $ and $\frac{\partial
V}{\partial z}=-\varrho I$ at the drain point $d=\left(  x_{out}%
,y_{out},0\right)  $, are satisfied,%
\[
\frac{\partial V\left(  x,y,0\right)  }{\partial z}=\varrho I\left(
\delta\left(  x-x_{in},y-y_{in},0\right)  -\delta\left(  x-x_{out}%
,y-y_{out},0\right)  \right)
\]
where the Dirac delta function $\delta\left(  x\right)  $ underlines
\textquotedblleft point\textquotedblright\ in- and ejection of the current.

The derivative%
\begin{equation}
\frac{\partial V}{\partial z}=\sum_{k=0}^{\infty}\sum_{l=0}^{\infty}%
\sqrt{\left(  \frac{k\pi}{a}\right)  ^{2}+\left(  \frac{l\pi}{b}\right)  ^{2}%
}\varphi_{kl}\cos\frac{k\pi}{a}x\cos\frac{l\pi}{b}y\sinh\left(  \sqrt{\left(
\frac{k\pi}{a}\right)  ^{2}+\left(  \frac{l\pi}{b}\right)  ^{2}}\left(
z-c\right)  \right)  \label{current-z_direction}%
\end{equation}
at the source $s$ and drain $d$ is%
\begin{align*}
\left.  \frac{\partial V}{\partial z}\right\vert _{s=\left(  x\,_{in}%
,y_{in},0\right)  }  &  =\sum_{k=0}^{\infty}\sum_{l=0}^{\infty}\sqrt{\left(
\frac{k\pi}{a}\right)  ^{2}+\left(  \frac{l\pi}{b}\right)  ^{2}}\varphi
_{kl}\cos\frac{k\pi}{a}x_{in}\cos\frac{l\pi}{b}y_{in}\sinh\left(
\sqrt{\left(  \frac{k\pi}{a}\right)  ^{2}+\left(  \frac{l\pi}{b}\right)  ^{2}%
}c\right) \\
&  =\varrho I\delta\left(  x-x_{in},y-y_{in},0\right) \\
\left.  \frac{\partial V}{\partial z}\right\vert _{d=\left(  x\,_{out}%
,y_{out},0\right)  }  &  =\sum_{k=0}^{\infty}\sum_{l=0}^{\infty}\sqrt{\left(
\frac{k\pi}{a}\right)  ^{2}+\left(  \frac{l\pi}{b}\right)  ^{2}}\varphi
_{kl}\cos\frac{k\pi}{a}x_{out}\cos\frac{l\pi}{b}y_{out}\sinh\left(
\sqrt{\left(  \frac{k\pi}{a}\right)  ^{2}+\left(  \frac{l\pi}{b}\right)  ^{2}%
}c\right) \\
&  =-\varrho I\delta\left(  x-x_{out},y-y_{out},0\right)
\end{align*}
while $\left.  \frac{\partial V}{\partial z}\right\vert _{u}=0$ at any other
point $u=\left(  x,y,0\right)  $.

\subsection{Current injection in the plane $z=0$}

\label{sec_current_point_injection_3D}After letting%
\[
\psi_{kl}=\varphi_{kl}\sqrt{\left(  \frac{k\pi}{a}\right)  ^{2}+\left(
\frac{l\pi}{b}\right)  ^{2}}\sinh\left(  c\sqrt{\left(  \frac{k\pi}{a}\right)
^{2}+\left(  \frac{l\pi}{b}\right)  ^{2}}\right)
\]
we will reformulate the current injection $j_{z}=\sigma\frac{\partial
V}{\partial z}$ in the plane $z=0$ in (\ref{current-z_direction}) more
generally as%
\[
\left.  \frac{\partial V}{\partial z}\right\vert _{z=0}=\sum_{k=0}^{\infty
}\sum_{l=0}^{\infty}\psi_{kl}\cos\frac{k\pi}{a}x\cos\frac{l\pi}{b}y=g\left(
x,y\right)
\]
The goal is to express the coefficients $\psi_{kl}$ in terms of the boundary
current in- and ejection function $g\left(  x,y\right)  $.

Multiplying both sides with $\cos\frac{k^{\prime}\pi}{a}x\cos\frac{l^{\prime
}\pi}{b}y$,
\[
g\left(  x,y\right)  \cos\frac{k^{\prime}\pi}{a}x\cos\frac{l^{\prime}\pi}%
{b}y=\sum_{k=0}^{\infty}\sum_{l=0}^{\infty}\psi_{kl}\cos\frac{k^{\prime}\pi
}{a}x\cos\frac{k\pi}{a}x\cos\frac{l\pi}{b}y\cos\frac{l^{\prime}\pi}{b}y
\]
Integrating $x$ from $0$ to $a$ and $y$ from 0 to $b$,%
\[
\int_{0}^{a}dx\;\int_{0}^{b}dy\;g\left(  x,y\right)  \cos\frac{k^{\prime}\pi
}{a}\cos\frac{l^{\prime}\pi}{b}y=\sum_{k=0}^{\infty}\sum_{l=0}^{\infty}%
\psi_{kl}\left(  \int_{0}^{a}\cos\frac{k^{\prime}\pi}{a}x\cos\frac{k\pi}%
{a}xdx\right)  \left(  \int_{0}^{b}\cos\frac{l\pi}{b}y\cos\frac{l^{\prime}\pi
}{b}y\right)
\]
invoking the orthogonality relations $\int_{0}^{a}\cos\frac{k^{\prime}\pi}%
{a}x\cos\frac{k\pi}{a}xdx=\frac{a}{\pi}\int_{0}^{\pi}\cos k^{\prime}\theta\cos
k\theta d\theta=a\left(  1_{\left\{  k=0\right\}  }+\frac{1}{2}1_{\left\{
k\neq0\right\}  }\right)  \delta_{k^{\prime}k}$, leads to%
\begin{equation}
\psi_{k^{\prime}l^{\prime}}=\left\{
\begin{array}
[c]{lc}%
\frac{4}{ab}\int_{0}^{a}dx\;\int_{0}^{b}dy\;g\left(  x,y\right)  \cos
\frac{k^{\prime}\pi}{a}\cos\frac{l^{\prime}\pi}{b}y & \text{for }k^{\prime
}>0\text{ and }l^{\prime}>0\\
\frac{2}{ab}\int_{0}^{a}dx\;\int_{0}^{b}dy\;g\left(  x,y\right)  \cos
\frac{k^{\prime}\pi}{a}\cos\frac{l^{\prime}\pi}{b}y & \left\{
\begin{array}
[c]{c}%
\text{for }k^{\prime}=0\text{ and }l^{\prime}>0\\
\text{or }k^{\prime}>0\text{ and }l^{\prime}=0
\end{array}
\right. \\
\frac{1}{ab}\int_{0}^{a}dx\;\int_{0}^{b}dy\;g\left(  x,y\right)  & \text{for
}k^{\prime}=0\text{ and }l^{\prime}=0
\end{array}
\right.  \label{Fourier_coeff_psi}%
\end{equation}

In our setting, we must choose the function $g\left(  x,y\right)  $ to satisfy
the current injection at the source point $s=\left(  x_{in},y_{in}\right)  $
and the current ejection at the drain point $d=\left(  x_{out},y_{out}\right)
$,%
\[
g\left(  x,y\right)  =\varrho I\left\{  \delta\left(  x-x_{in}\right)
\delta\left(  y-y_{in}\right)  -\delta\left(  x-x_{out}\right)  \delta\left(
y-y_{out}\right)  \right\}
\]
Hence, the double integral in $\psi_{k^{\prime}l^{\prime}}$ becomes, with
$c_{k^{\prime}l^{\prime}}$ representing the three different conditions for
$k^{\prime}$ and $l^{\prime}$,
\begin{align*}
\psi_{k^{\prime}l^{\prime}}  &  =\frac{\rho Ic_{k^{\prime}l^{\prime}}}{ab}%
\int_{0}^{a}dx\;\int_{0}^{b}dy\left(  \delta\left(  x-x_{in}\right)
\delta\left(  y-y_{in}\right)  -\delta\left(  x-x_{out}\right)  \delta\left(
y-y_{out}\right)  \right)  \cos\frac{k^{\prime}\pi}{a}x\;\cos\frac{l^{\prime
}\pi}{b}y\\
&  =\frac{\rho Ic_{k^{\prime}l^{\prime}}}{ab}\int_{0}^{a}\delta\left(
x-x_{in}\right)  \cos\frac{k^{\prime}\pi}{a}x\;dx\int_{0}^{b}\delta\left(
y-y_{in}\right)  \cos\frac{l^{\prime}\pi}{b}y\;dy\\
&  \hspace{0.5cm}-\frac{\rho Ic_{k^{\prime}l^{\prime}}}{ab}\int_{0}^{a}%
\delta\left(  x-x_{out}\right)  \cos\frac{k^{\prime}\pi}{a}x\;dx\int_{0}%
^{b}\delta\left(  y-y_{out}\right)  \;\cos\frac{l^{\prime}\pi}{b}y\;dy\\
&  =\frac{\rho Ic_{k^{\prime}l^{\prime}}}{ab}\left(  \cos\frac{k^{\prime}\pi
}{a}x_{in}\;\cos\frac{l^{\prime}\pi}{b}y_{in}-\cos\frac{k^{\prime}\pi}%
{a}x_{out}\;\cos\frac{l^{\prime}\pi}{b}y_{out}\right)
\end{align*}
which simplifies the above with%
\[
q_{kl}=\frac{\rho I}{ab}\left(  \cos\frac{k\pi}{a}x_{in}\;\cos\frac{l\pi}%
{b}y_{in}-\cos\frac{k\pi}{a}x_{out}\;\cos\frac{l\pi}{b}y_{out}\right)
\]
to%
\[
\psi_{kl}=\left\{
\begin{array}
[c]{lc}%
4q_{kl} & \text{for }k>0\text{ and }l>0\\
2q_{kl} & \text{for }\left\{  k=0\text{ and }l>0\right\}  \text{ or }\left\{
k>0\text{ and }l=0\right\} \\
q_{00} & \text{for }k^{\prime}=0\text{ and }l^{\prime}=0
\end{array}
\right.
\]

In conclusion, we find the potential in (\ref{V_general_linear_combination})%
\[
V\left(  x,y,z\right)  =\sum_{k=0}^{\infty}\sum_{l=0}^{\infty}\varphi_{kl}%
\cos\frac{k\pi}{a}x\cos\frac{l\pi}{b}y\cosh\left(  \sqrt{\left(  \frac{k\pi
}{a}\right)  ^{2}+\left(  \frac{l\pi}{b}\right)  ^{2}}\left(  z-c\right)
\right)
\]
with the coefficients%
\[
\varphi_{kl}=\left\{
\begin{array}
[c]{lc}%
4p_{kl} & \text{for }k>0\text{ and }l>0\\
2p_{kl} & \text{for }\left\{  k=0\text{ and }l>0\right\}  \text{ or }\left\{
k>0\text{ and }l=0\right\} \\
p_{00} & \text{for }k^{\prime}=0\text{ and }l^{\prime}=0
\end{array}
\right.
\]
where, for $\left(  k,l\right)  \neq\left(  0,0\right)  $,%
\[
p_{kl}=\frac{\frac{\rho I}{ab}\left(  \cos\frac{k\pi}{a}x_{in}\;\cos\frac
{l\pi}{b}y_{in}-\cos\frac{k\pi}{a}x_{out}\;\cos\frac{l\pi}{b}y_{out}\right)
}{\sinh\left(  c\sqrt{\left(  \frac{k\pi}{a}\right)  ^{2}+\left(  \frac{l\pi
}{b}\right)  ^{2}}\right)  \sqrt{\left(  \frac{k\pi}{a}\right)  ^{2}+\left(
\frac{l\pi}{b}\right)  ^{2}}}%
\]
Hence,
\[
p_{k0}=\frac{\frac{\rho I}{ab}\left(  \cos\frac{k\pi}{a}x_{in}\;-\cos
\frac{k\pi}{a}x_{out}\;\right)  }{\frac{k\pi}{a}\sinh\left(  c\frac{k\pi}%
{a}\right)  }%
\]
and
\[
p_{0l}=\frac{\frac{\rho I}{ab}\left(  \cos\frac{l\pi}{b}y_{in}-\cos\frac{l\pi
}{b}y_{out}\right)  }{\frac{l\pi}{b}\sinh\left(  c\frac{l\pi}{b}\right)  }%
\]
but $p_{00}$ is constant, which will produce the reference potential
$V_{\text{ref}}$.

\subsection{Summary of the analytic solution in 3D}

\subsubsection{Potential $V(x,y,z)$ in 3D}

The explicit solution of the potential at the point $\left(  x,y,z\right)  $
inside and on the boundary of the beam, i.e. for $0\leq x\leq a$, $0\leq y\leq
b$ and $0\leq z\leq c$, then is%
\begin{align}
V\left(  x,y,z\right)  +V_{\text{ref}}  &  =\frac{2\varrho I}{ab}\sum
_{k=1}^{\infty}\frac{\left(  \cos\frac{k\pi}{a}x_{in}\;-\cos\frac{k\pi}%
{a}x_{out}\;\right)  }{\frac{k\pi}{a}\sinh\left(  c\frac{k\pi}{a}\right)
}\cos\frac{k\pi}{a}x\cosh\left(  \frac{k\pi}{a}\left(  z-c\right)  \right)
\nonumber\\
&  \hspace{0.5cm}+\frac{2\varrho I}{ab}\sum_{l=1}^{\infty}\frac{\left(
\cos\frac{l\pi}{b}y_{in}-\cos\frac{l\pi}{b}y_{out}\right)  }{\frac{l\pi}%
{b}\sinh\left(  c\frac{l\pi}{b}\right)  }\cos\frac{l\pi}{b}y\cosh\left(
\frac{l\pi}{b}\left(  z-c\right)  \right) \nonumber\\
&  \hspace{0.5cm}+\frac{4\varrho I}{ab}\sum_{k=1}^{\infty}\sum_{l=1}^{\infty
}\frac{\left(  \cos\frac{k\pi}{a}x_{in}\;\cos\frac{l\pi}{b}y_{in}-\cos
\frac{k\pi}{a}x_{out}\;\cos\frac{l\pi}{b}y_{out}\right)  }{\sinh\left(
c\sqrt{\left(  \frac{k\pi}{a}\right)  ^{2}+\left(  \frac{l\pi}{b}\right)
^{2}}\right)  \sqrt{\left(  \frac{k\pi}{a}\right)  ^{2}+\left(  \frac{l\pi}%
{b}\right)  ^{2}}}\nonumber\\
&  \hspace{0.5cm}\hspace{1cm}\times\cos\frac{k\pi}{a}x\cos\frac{l\pi}{b}%
y\cosh\left(  \sqrt{\left(  \frac{k\pi}{a}\right)  ^{2}+\left(  \frac{l\pi}%
{b}\right)  ^{2}}\left(  z-c\right)  \right)  \label{Voltage_beam}%
\end{align}
Perhaps needless to recall, the voltage is only defined with respect to some
reference potential $V_{\text{ref}}$. The periodicity of the double Fourier
series (\ref{Voltage_beam}) for $V\left(  x,y,z\right)  $ in the beam implies
that the extension of the real number $x$, $y$ and $z$ to all real numbers,
beyond $0\leq x\leq a$, $0\leq y\leq b$ and $0\leq z\leq c$, leads to a
coverage of the 3D space by a repeated parallellepipidum, which can be
regarded as the 3D extension of elliptic functions in the complex plane (see
Appendix \ref{sec_Theta_functions_in_Green_function}). In general, for
sufficiently large difference $\left\vert z-c\right\vert $, the double Fourier
series (\ref{Voltage_beam}) in 3D converges reasonably well.

\subsubsection{Current density in 3D}

The current density, for $0\leq x\leq a$, $0\leq y\leq b$ and $0\leq z\leq
c$,
\[
\overrightarrow{j}\left(  x,y,z\right)  =\sigma\operatorname{grad}%
V=-\sigma\left(  \frac{\partial V}{\partial x},\frac{\partial V}{\partial
y},\frac{\partial V}{\partial z}\right)
\]
is a vector with magnitude equal to%
\[
\left\vert \overrightarrow{j}\left(  x,y,z\right)  \right\vert =\left\vert
j\left(  x,y,z\right)  \right\vert =\sigma\sqrt{\left(  \frac{\partial
V}{\partial x}\right)  ^{2}+\left(  \frac{\partial V}{\partial y}\right)
^{2}+\left(  \frac{\partial V}{\partial z}\right)  ^{2}}%
\]
The partial derivatives $\frac{\partial V}{\partial x},\frac{\partial
V}{\partial y}$ and $\frac{\partial V}{\partial z}$ are readily computed from
(\ref{Voltage_beam}), e.g.:%
\begin{align*}
\frac{\partial V}{\partial x}  &  =-\frac{2\varrho I}{ab}\sum_{k=1}^{\infty
}\frac{\left(  \cos\frac{k\pi}{a}x_{in}\;-\cos\frac{k\pi}{a}x_{out}\;\right)
}{\sinh\left(  c\frac{k\pi}{a}\right)  }\sin\frac{k\pi}{a}x\cosh\left(
\frac{k\pi}{a}\left(  z-c\right)  \right) \\
&  \hspace{0.5cm}-\frac{4\varrho I}{ab}\sum_{k=1}^{\infty}\sum_{l=1}^{\infty
}\frac{k\pi}{a}\frac{\left(  \cos\frac{k\pi}{a}x_{in}\;\cos\frac{l\pi}%
{b}y_{in}-\cos\frac{k\pi}{a}x_{out}\;\cos\frac{l\pi}{b}y_{out}\right)  }%
{\sinh\left(  c\sqrt{\left(  \frac{k\pi}{a}\right)  ^{2}+\left(  \frac{l\pi
}{b}\right)  ^{2}}\right)  \sqrt{\left(  \frac{k\pi}{a}\right)  ^{2}+\left(
\frac{l\pi}{b}\right)  ^{2}}}\\
&  \hspace{0.5cm}\hspace{1cm}\times\sin\frac{k\pi}{a}x\cos\frac{l\pi}{b}%
y\cosh\left(  \sqrt{\left(  \frac{k\pi}{a}\right)  ^{2}+\left(  \frac{l\pi}%
{b}\right)  ^{2}}\left(  z-c\right)  \right)
\end{align*}
Unfortunately, the series of $\left\vert j\left(  x,y,z\right)  \right\vert $
for $z=0$ converge \emph{too slowly to be of practical use}!

\section{Reduction to two dimensions}

\label{sec_reduction_2D}The Laplace equation (\ref{Laplacian_differential_eq})
in 2D, $\frac{\partial^{2}V}{\partial x^{2}}+\frac{\partial^{2}V}{\partial
y^{2}}=0$, has as general solution a harmonic function $V\left(  x,y\right)
$. A non-constant harmonic function $V\left(  x,y\right)  $ cannot have an
extremum, either minimum or maximum, at an interior point of the region
\cite[p. 167]{Titchmarshfunctions}. Imposing that no current can enter the
rectangle, i.e. $\left.  \frac{\partial V}{\partial x}\right\vert
_{x=0}=\left.  \frac{\partial V}{\partial x}\right\vert _{x=a}=0$ and $\left.
\frac{\partial V}{\partial y}\right\vert _{y=0}=\left.  \frac{\partial
V}{\partial y}\right\vert _{y=b}=0$, implies that the extreme value of the
harmonic function $V\left(  x,y\right)  $ occurs at the boundary. Hence,
current can only flow from one boundary to another, without an internal high
potential at $s=\left(  x\,_{in},y_{in}\right)  $ and internal low potential
at $d=\left(  x_{out},y_{out}\right)  $. Apart from the physical vector
decomposition argument in the introduction, the harmonic function argument
shows why a 3D computation is necessary for the current in- and ejection.

However, the current injection without vector decomposition into the plane can
be regarded as a pile-up of charges at point $s$ when the parallellepipidum is
shrinked to a 2D rectangle. Hence, in 2D, we may consider a plate or rectangle
with two opposite electric charges at $s$ and $d$. The Laplacian equation is
then replaced by a (normalized\footnote{The permittivity $\epsilon_{0}$ in
$\Delta V=-\frac{1}{\epsilon_{0}}\rho$ has disappeared.}) Poisson's equation
$\Delta V=-\varphi\left(  r\right)  $, where the charge density $\varphi
\left(  r\right)  $ in 2D consists of two point charges,
\[
\Delta V=\frac{\partial^{2}V}{\partial x^{2}}+\frac{\partial^{2}V}{\partial
y^{2}}=\delta\left(  x\,_{in},y_{in}\right)  -\delta\left(  x\,_{out}%
,y_{out}\right)
\]
with boundary conditions $\left.  \frac{\partial V}{\partial x}\right\vert
_{x=0}=\left.  \frac{\partial V}{\partial x}\right\vert _{x=a}=0$ and $\left.
\frac{\partial V}{\partial y}\right\vert _{y=0}=\left.  \frac{\partial
V}{\partial y}\right\vert _{y=b}=0$, indicating that no current can flow
through the boundary. Our new formulation is related to the \textquotedblleft
Green's function for the Laplacian $\Delta V$ in the interior of a
rectangle\textquotedblright\ briefly reviewed in Appendix
\ref{sec_Green_functions}, whose solution \cite[p. 384]{Courant_HilbertI}
involves elliptic functions when the boundary conditions are $V\left(
0,y\right)  =V\left(  a,y\right)  =V\left(  x,0\right)  =V\left(  x,b\right)
=0$, i.e. $V=0$ at the boundary of the rectangle. The double Fourier series of
the Green function for our Neumann boundary conditions is computed in \cite[p.
399]{Zauderer}. We first derive that double Fourier series.

\subsection{Potential $V(x,y)$ in 2D}

\label{sec_potential_2D}The voltage in the $xy$-plane (i.e. $z=0$) follows
from the 3D voltage (\ref{Voltage_beam}), for $0\leq x\leq a$ and $0\leq y\leq
b$, as%
\begin{align}
V\left(  x,y,0\right)  +V_{\text{ref}}  &  =\frac{2\varrho I}{ab}\sum
_{k=1}^{\infty}\frac{\left(  \cos\frac{k\pi}{a}x_{in}\;-\cos\frac{k\pi}%
{a}x_{out}\;\right)  }{\frac{k\pi}{a}\tanh\left(  c\frac{k\pi}{a}\right)
}\cos\frac{k\pi}{a}x\nonumber\\
&  \hspace{0.5cm}+\frac{2\varrho I}{ab}\sum_{l=1}^{\infty}\frac{\left(
\cos\frac{l\pi}{b}y_{in}-\cos\frac{l\pi}{b}y_{out}\right)  }{\frac{l\pi}%
{b}\tanh\left(  c\frac{l\pi}{b}\right)  }\cos\frac{l\pi}{b}y\nonumber\\
&  \hspace{0.5cm}+\frac{4\varrho I}{ab}\sum_{k=1}^{\infty}\sum_{l=1}^{\infty
}\frac{\left(  \cos\frac{k\pi}{a}x_{in}\;\cos\frac{l\pi}{b}y_{in}-\cos
\frac{k\pi}{a}x_{out}\;\cos\frac{l\pi}{b}y_{out}\right)  }{\tanh\left(
c\sqrt{\left(  \frac{k\pi}{a}\right)  ^{2}+\left(  \frac{l\pi}{b}\right)
^{2}}\right)  \sqrt{\left(  \frac{k\pi}{a}\right)  ^{2}+\left(  \frac{l\pi}%
{b}\right)  ^{2}}}\cos\frac{k\pi}{a}x\cos\frac{l\pi}{b}y\nonumber
\end{align}
The Taylor series $\pi\cot\left(  \pi x\right)  =\frac{1}{x}-2\sum
_{n=1}^{\infty}\zeta\left(  2n\right)  \,x^{2n-1}$, valid\footnote{The Taylor
series follows directly from the generating function
(\ref{gf_Bernoullinumbers}) below of the Bernoulli numbers \cite[art.
240]{PVM_Mittag-Leffler_Gamma}.} for $\left\vert x\right\vert <1$ and where
$\zeta\left(  s\right)  $ is the Riemann Zeta-function \cite{Titchmarshzeta},
indicates, with $\pi\cot\left(  i\pi x\right)  =\frac{\pi}{i}\coth\left(  \pi
x\right)  $, that%
\[
\coth\left(  z\right)  =\frac{1}{\tanh\left(  z\right)  }=\frac{1}{z}%
-2\sum_{n=1}^{\infty}\frac{\zeta\left(  2n\right)  }{\pi^{2n}}\left(
-1\right)  ^{n}\,z^{2n-1}\text{ for }\left\vert z\right\vert <\pi
\]
Hence, for small $c$, we find that
\begin{align}
V\left(  x,y,0\right)  +V_{\text{ref}}  &  =\frac{2\varrho I}{abc}\sum
_{k=1}^{\infty}\frac{\left(  \cos\frac{k\pi}{a}x_{in}\;-\cos\frac{k\pi}%
{a}x_{out}\;\right)  }{\left(  \frac{k\pi}{a}\right)  ^{2}}\cos\frac{k\pi}%
{a}x+O\left(  c\right) \nonumber\\
&  \hspace{0.5cm}+\frac{2\varrho I}{abc}\sum_{l=1}^{\infty}\frac{\left(
\cos\frac{l\pi}{b}y_{in}-\cos\frac{l\pi}{b}y_{out}\right)  }{\left(
\frac{l\pi}{b}\right)  ^{2}}\cos\frac{l\pi}{b}y\nonumber\\
&  \hspace{0.5cm}+\frac{4\varrho I}{abc}\sum_{k=1}^{\infty}\sum_{l=1}^{\infty
}\frac{\left(  \cos\frac{k\pi}{a}x_{in}\;\cos\frac{l\pi}{b}y_{in}-\cos
\frac{k\pi}{a}x_{out}\;\cos\frac{l\pi}{b}y_{out}\right)  }{\left(  \frac{k\pi
}{a}\right)  ^{2}+\left(  \frac{l\pi}{b}\right)  ^{2}}\cos\frac{k\pi}{a}%
x\cos\frac{l\pi}{b}y \label{Voltage_2D_orderterm}%
\end{align}
where the volume of the parallelepiped is $v_{3D}=abc$. For $c\rightarrow0$,
i.e. omitting the order terms $O\left(  c\right)  $ in
(\ref{Voltage_2D_orderterm}), and after renormalization by $\frac{\varrho
I}{c}=1$, we denote the 2D voltage as $V\left(  x,y\right)  =$ $\lim
_{c\rightarrow0}V\left(  x,y,0\right)  $, subject to $\frac{\varrho I}{c}=1$,
whose result is%
\begin{align}
V\left(  x,y\right)  +V_{\text{ref}}  &  =\frac{2}{ab}\left(  \sum
_{k=1}^{\infty}\frac{\left(  \cos\frac{k\pi}{a}x_{in}\;-\cos\frac{k\pi}%
{a}x_{out}\;\right)  }{\left(  \frac{k\pi}{a}\right)  ^{2}}\cos\frac{k\pi}%
{a}x+\sum_{l=1}^{\infty}\frac{\left(  \cos\frac{l\pi}{b}y_{in}-\cos\frac{l\pi
}{b}y_{out}\right)  }{\left(  \frac{l\pi}{b}\right)  ^{2}}\cos\frac{l\pi}%
{b}y\right) \nonumber\\
&  \hspace{0.5cm}+\frac{4}{ab}\sum_{k=1}^{\infty}\sum_{l=1}^{\infty}%
\frac{\left(  \cos\frac{k\pi}{a}x_{in}\;\cos\frac{l\pi}{b}y_{in}-\cos
\frac{k\pi}{a}x_{out}\;\cos\frac{l\pi}{b}y_{out}\right)  }{\left(  \frac{k\pi
}{a}\right)  ^{2}+\left(  \frac{l\pi}{b}\right)  ^{2}}\cos\frac{k\pi}{a}%
x\cos\frac{l\pi}{b}y \label{Voltage_2D}%
\end{align}
which is precisely
\begin{equation}
V\left(  x,y\right)  =K\left(  x,y;x_{in},y_{in}\right)  -K\left(
x,y;x_{out},y_{out}\right)  \label{Voltage_2D_Green}%
\end{equation}
in the Green's function approach (\ref{Green'sFunction_K_2D}). In summary, we
have shown mathematically that the dimension reduction from 3D to 2D indeed is
equivalent to placing point charges at $s$ and $d$.

We explain the renormalization of $\frac{\varrho I}{c}=1$. In order that
$\lim_{c\rightarrow0}\frac{\varrho I}{c}$ is finite, the current $I$ must tend
to zero, given that the resistivity $\varrho$ of the material is
fixed\footnote{A zero resistivity $\varrho=0$ would imply superconductivity,
which we exclude.}, also in 2D. In other words, there is no injected current
anymore, but only not-moving charges. In any case, a current injected
orthogonal to the plane must bend in a single point over 90 degrees which is
mechanically speaking or in vector terminology not possible: hence, the
current accumulates to a charge peak at that single point $s$.%

%TCIMACRO{\FRAME{fhFU}{15.1281cm}{6.8491cm}{0pt}{\Qcb{Example of the voltage
%$V\left(  x,y\right)  $, computed from (\ref{Voltage_2D}) with $K=100$ terms
%in the Fourier series.}}{\Qlb{Fig_voltage1}}{voltage1.ps}%
%width 15.1281cm;  height 6.8491cm;  depth 0pt;  original-width 8.2538in;
%original-height 11.6949in;  cropleft "0";  croptop "0.6471";  cropright "1";
%cropbottom "0.3294";  filename 'Voltage1.ps';file-properties "XNPEU";}} }%
%BeginExpansion
\begin{figure}
[h]
\begin{center}
\includegraphics[
trim=0.000000in 3.852300in 0.000000in 4.127131in,
height=6.8491cm,
width=15.1281cm
]%
{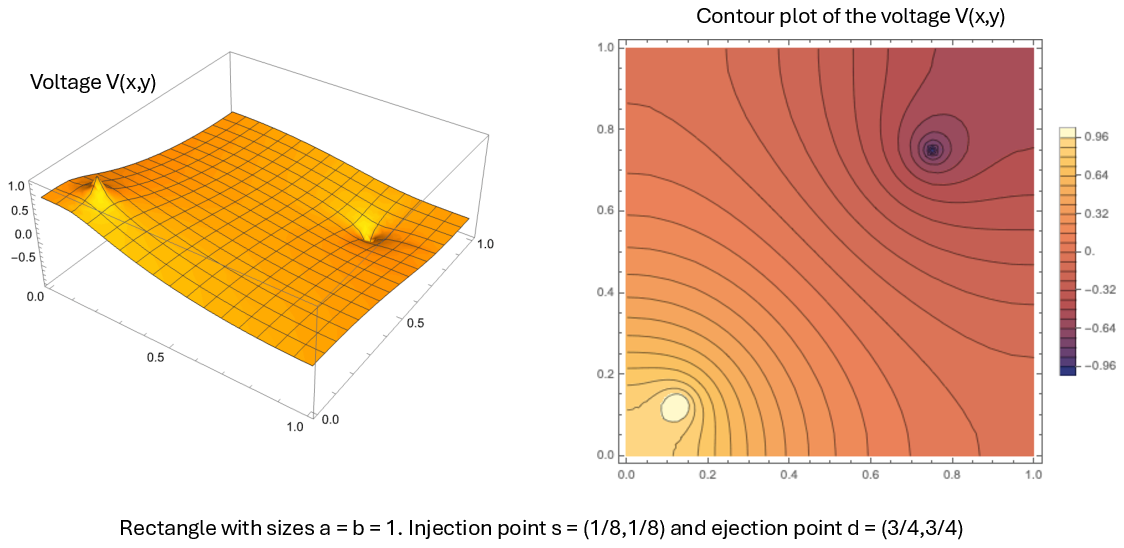}%
\caption{Example of the voltage $V\left(  x,y\right)  $, computed from
(\ref{Voltage_2D}) with $K=100$ terms in the Fourier series.}%
\label{Fig_voltage1}%
\end{center}
\end{figure}
%EndExpansion
%TCIMACRO{\FRAME{fhFU}{15.1391cm}{7.3565cm}{0pt}{\Qcb{The magnitude $\left\vert
%j\left(  x,y\right)  \right\vert $ of the current density corresponding to
%Fig. \ref{Fig_voltage1}, computed from (\ref{Single_sum_current_density}) with
%$K=100$ terms.}}{\Qlb{Fig_cdensity1}}{cdensity1.ps}%
%width 15.1391cm;  height 7.3565cm;  depth 0pt;  original-width 8.2538in;
%original-height 11.6949in;  cropleft "0";  croptop "0.6707";  cropright "1";
%cropbottom "0.3292";  filename 'Cdensity1.ps';file-properties "XNPEU";}} }%
%BeginExpansion
\begin{figure}
[h]
\begin{center}
\includegraphics[
trim=0.000000in 3.849961in 0.000000in 3.851130in,
height=7.3565cm,
width=15.1391cm
]%
{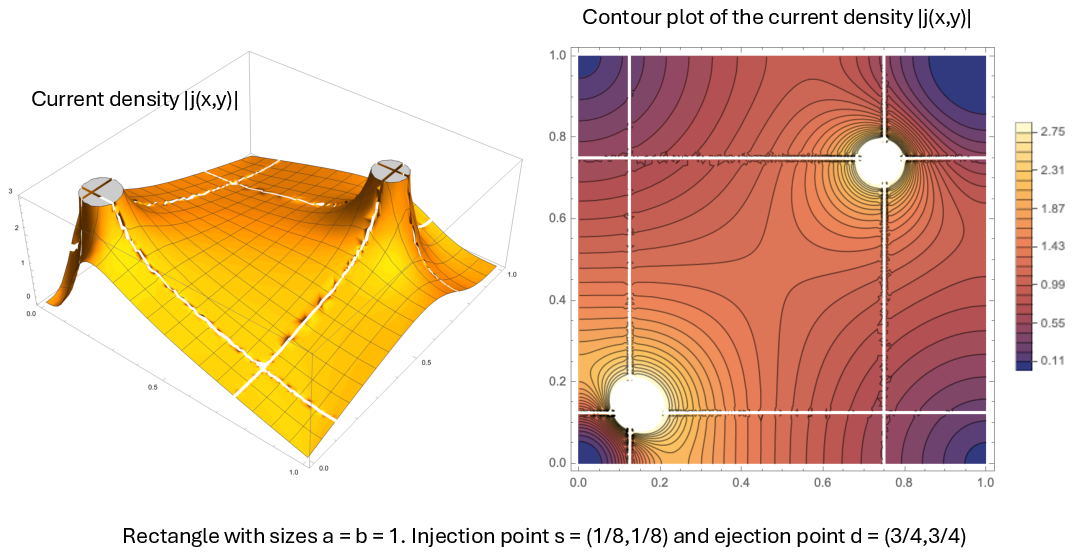}%
\caption{The magnitude $\left\vert j\left(  x,y\right)  \right\vert $ of the
current density corresponding to Fig. \ref{Fig_voltage1}, computed from
(\ref{Single_sum_current_density}) with $K=100$ terms.}%
\label{Fig_cdensity1}%
\end{center}
\end{figure}
%EndExpansion

Although our voltage (\ref{Voltage_beam}) in 3D and the 2D reduction
(\ref{Voltage_2D}) are exact, the double Fourier series converge very slowly
for $V\left(  x,y,0\right)  $, i.e. in the $xy$-plane. Fig. \ref{Fig_voltage1}
shows both a 3D plot and the contour plot of the potential $V\left(
x,y\right)  =V\left(  x,y,0\right)  $ in (\ref{Voltage_2D}) for $c\rightarrow
0$ or in (\ref{Voltage_2D_Green}) in a square with sizes $a=b=1$. The current
injection is in $s=\left(  x_{in},y_{in}\right)  =\left(  1/8,1/8\right)  $
and current ejection occurs at $d=\left(  x_{out},y_{out}\right)  =\left(
3/4,3/4\right)  $. We underline that the computation of the voltage $V\left(
x,y\right)  $ is involved, but possible in a reasonable time, whereas the
current density is problematic, due to the \emph{extremely slow convergence}.
Fortunately, Fig. \ref{Fig_cdensity1} computed with the single Fourier series
(\ref{Single_sum_current_density}) below draws the current density
corresponding to the voltage in Fig. \ref{Fig_voltage1}.

We remark that Mathematica in Fig. \ref{Fig_cdensity1} plots white lines at
the peculiar discontinuities of the Fourier coefficients, the functions
$g_{n}$ in (\ref{Fourier_coeff_gn}) and $h_{n}$ in (\ref{Fourier_coeff_hn}),
whereas those discontinuities disappear in the Fourier series (with infinite
terms), because both the voltage $V\left(  x,y\right)  $ and the magnitude
$\left\vert j\left(  x,y\right)  \right\vert $ of the current density are
continuous\footnote{Indeed, the potential $V\left(  x,y\right)  $ satisfies
the Poisson second-order differential equation and any differentiable function
(such as the partial derivatives $\frac{\partial V}{\partial x}$ and
$\frac{\partial V}{\partial y})$ is continuous.} for any $x$ and $y$ in the
rectangle, apart from the in- and ejection points $s$ and $d$. The
inaccuracies at the white lines in Fig. \ref{Fig_cdensity1} has motivated
further analysis, that eventually led to a much faster converging, so called
$q$-series in Appendix \ref{sec_q_series}. Fig.
\ref{Fig_currentdensity_fourier_qserie} illustrates the smoothness around the
white lines, computed with the $q$-series method.

\subsection{From a double to single Fourier series with quasi-Green's
function}

\label{sec_Green_2D_single_sum}

The very slow convergence of partial derivatives of the double Fourier series
(\ref{Voltage_2D})-(\ref{Voltage_2D_Green}) has motivated the search for a
faster computation, without loosing accuracy. Following Morse and Feshbach
\cite[p. 799-800]{Morse} and recalling the rewriting
(\ref{Green'sFunction_K_2D_short}) of the 2D voltage (\ref{Voltage_2D}) with
coefficients $\gamma_{nm}$,%
\[
K\left(  x,y;\xi,\eta\right)  =\sum_{n=0}^{\infty}\left(  \frac{4}{ab}%
\sum_{m=0}^{\infty}\gamma_{nm}\frac{\cos\left(  \frac{\pi n}{a}\xi\right)
\cos\left(  \frac{\pi m}{b}y\right)  \cos\left(  \frac{\pi m}{b}\eta\right)
}{\left(  \frac{\pi n}{a}\right)  ^{2}+\left(  \frac{\pi m}{b}\right)  ^{2}%
}\right)  \cos\left(  \frac{\pi n}{a}x\right)
\]
we propose the solution%
\[
F\left(  x,y\right)  =\sum_{n=0}^{\infty}f_{n}\left(  y\right)  \cos\left(
\frac{\pi n}{a}x\right)
\]
of the Poisson equation with general charge density $\varphi\left(
x,y\right)  $,
\begin{equation}
\frac{\partial^{2}F}{\partial x^{2}}+\frac{\partial^{2}F}{\partial y^{2}%
}=-\varphi\left(  x,y\right)  \label{Poisson_2D}%
\end{equation}
The Neumann boundary conditions $\operatorname{grad}F=0$ at the boundary of
the rectangle $\left(  0,0\right)  ,\left(  0,b\right)  ,\left(  a,b\right)  $
and $\left(  a,0\right)  $ reflect that no current can flow through the
boundaries of the rectangle. In addition, we propose%
\begin{equation}
\varphi\left(  x,y\right)  =\sum_{n=0}^{\infty}\rho_{n}\left(  y\right)
\cos\left(  \frac{\pi n}{a}x\right)  \label{Proposal_phi(x,y)}%
\end{equation}
where it holds, via orthogonality of Fourier series in Section
\ref{sec_orthogonality_relations},%
\begin{align*}
\int_{0}^{a}\varphi\left(  x,y\right)  \cos\left(  \frac{\pi k}{a}x\right)
dx  &  =\sum_{n=0}^{\infty}\rho_{n}\left(  y\right)  \int_{0}^{a}\cos\left(
\frac{\pi k}{a}x\right)  \cos\left(  \frac{\pi n}{a}x\right)  dx\\
&  =\frac{1}{2}\sum_{n=0}^{\infty}\rho_{n}\left(  y\right)  \int_{0}^{a}%
\cos\left(  \frac{\pi\left(  k-n\right)  }{a}x\right)  dx+\frac{1}{2}%
\sum_{n=0}^{\infty}\rho_{n}\left(  y\right)  \int_{0}^{a}\cos\left(  \frac
{\pi\left(  k+n\right)  }{a}x\right)  dx\\
&  =\frac{1}{2}\sum_{n=0}^{\infty}\rho_{n}\left(  y\right)  \frac{\sin\left(
\pi\left(  k-n\right)  \right)  }{\frac{\pi\left(  k-n\right)  }{a}}+\frac
{1}{2}\sum_{n=0}^{\infty}\rho_{n}\left(  y\right)  \frac{\sin\left(
\pi\left(  k+n\right)  \right)  }{\frac{\pi\left(  k+n\right)  }{a}}=\frac
{a}{2}\rho_{k}\left(  y\right)
\end{align*}
Thus, for $k\geq1$, we obtain that%
\[
\rho_{k}\left(  y\right)  =\frac{2}{a}\int_{0}^{a}\varphi\left(  x,y\right)
\cos\left(  \frac{\pi k}{a}x\right)  dx
\]
and, for $k=0$,%
\[
\rho_{0}\left(  y\right)  =\frac{1}{a}\int_{0}^{a}\varphi\left(  x,y\right)
dx
\]

However, our Neumann boundaries (i.e. $\operatorname{grad}F=0$; no current
through the boundary) are different from the theory in Morse and Feshbach
\cite[p. 799-800]{Morse}, which is based on Dirichlet boundaries (i.e. $F=0$;
the voltage is zero at the boundary), exemplified in Appendix
\ref{sec_Green_functions}. That difference causes a major complication,
because there does not exist a Green's function in our setting as shown below,
but we can define a \textquotedblleft quasi\textquotedblright\ Green's
function with two opposite delta functions instead of a single one in the
classical definition (\ref{Green_diff_equation_Poisson}) of a Green's
function. Since the current density is our main target, which is a partial
derivative of the potential, the unknown appearing constant below and the
non-existence of a Green's function, is alleviated. After removing this
barrier, we may continue with our \textquotedblleft quasi\textquotedblright%
\ Green's functions.

Introducing the proposal (\ref{Proposal_phi(x,y)}) into the Poisson
differential equation (\ref{Poisson_2D})%
\[
-\sum_{n=0}^{\infty}f_{n}\left(  y\right)  \left(  \frac{\pi n}{a}\right)
^{2}\cos\left(  \frac{\pi n}{a}x\right)  +\sum_{n=0}^{\infty}\frac{d^{2}%
f_{n}\left(  y\right)  }{dy^{2}}\cos\left(  \frac{\pi n}{a}x\right)
=-\sum_{n=0}^{\infty}\rho_{n}\left(  y\right)  \cos\left(  \frac{\pi n}%
{a}x\right)
\]
leads to the second-order differential equation%
\begin{equation}
\frac{d^{2}f_{n}\left(  y\right)  }{dy^{2}}-\left(  \frac{\pi n}{a}\right)
^{2}f_{n}\left(  y\right)  =-\rho_{n}\left(  y\right)
\label{second_order_dvg}%
\end{equation}
The boundary conditions $\frac{\partial F}{\partial y}=-\sum_{n=0}^{\infty
}\frac{df_{n}\left(  y\right)  }{dy}\cos\left(  \frac{\pi n}{a}x\right)  =0$
for $y=0$ and $y=b$ and all $x\in\left[  0,a\right]  $ translate to
$\frac{df_{n}\left(  y\right)  }{dy}=0$ for $y=0$ and $y=b$ for all integer
$n\geq0$. The homogeneous solution, satisfying $\frac{d^{2}h_{n}\left(
y\right)  }{dy^{2}}-\left(  \frac{\pi n}{a}\right)  ^{2}h_{n}\left(  y\right)
=0$ is%
\[
h_{n}\left(  y\right)  =A\cosh\frac{\pi n}{a}y+B\cosh\frac{\pi n}{a}\left(
b-y\right)
\]
where $y_{1}=\cosh\frac{\pi n}{a}y$ and $y_{2}=\cosh\frac{\pi n}{a}\left(
b-y\right)  $ are independent solutions with Wronskian $W\left(  y_{1}%
,y_{2}\right)  =\left\vert
\begin{array}
[c]{cc}%
y_{1} & y_{2}\\
y_{1}^{\prime} & y_{2}^{\prime}%
\end{array}
\right\vert =y_{1}y_{2}^{\prime}-y_{2}y_{1}^{\prime}$ equal to%
\begin{align*}
W\left(  y_{1},y_{2}\right)   &  =-\frac{\pi n}{a}\left(  \cosh\frac{\pi n}%
{a}y\sinh\frac{\pi n}{a}\left(  b-y\right)  +\cosh\frac{\pi n}{a}\left(
b-y\right)  \sinh\frac{\pi n}{a}y\right) \\
&  =-\frac{\pi n}{a}\sinh\left(  \frac{\pi n}{a}\left(  b-y\right)  +\frac{\pi
n}{a}y\right)  =-\frac{\pi n}{a}\sinh\left(  \frac{\pi n}{a}b\right)
\end{align*}
and constant (as it must be \cite{Morse}). Only for $n=0$, both $y_{1}=y_{2}$
are dependent! However, for $n=0$, the second-order differential equation
(\ref{second_order_dvg}) simplifies to%
\[
\frac{d^{2}f_{0}\left(  y\right)  }{dy^{2}}=-\rho_{0}\left(  y\right)
=-\frac{1}{a}\int_{0}^{a}\varphi\left(  x,y\right)  dx
\]
and integration yields%
\[
\frac{df_{0}\left(  y\right)  }{dy}-\left.  \frac{df_{0}\left(  y\right)
}{dy}\right\vert _{y=0}=-\frac{1}{a}\int_{0}^{a}\int_{0}^{y}\varphi\left(
x,u\right)  dxdu
\]
The boundary condition $\frac{df_{0}\left(  y\right)  }{dy}=0$ for $y=0$ leads
to%
\begin{equation}
\frac{df_{0}\left(  y\right)  }{dy}=-\frac{1}{a}\int_{0}^{a}\int_{0}%
^{y}\varphi\left(  x,u\right)  dxdu \label{D_f0(y)}%
\end{equation}
However, the boundary condition $\frac{df_{0}\left(  y\right)  }{dy}=0$ at
$y=b$ requires that $\int_{0}^{a}\int_{0}^{b}\varphi\left(  x,y\right)
dxdy=0$, which is physically correct, if the sum of all charges in the beam is
zero. This implies that a single delta function as for the Green's function in
(\ref{Green_diff_equation_Poisson}) is incompatible. Since our solution will
always requires an injection and ejection, implying a difference between two
delta charges, we ignore this incompatibility and the requirement $\int
_{0}^{a}\int_{0}^{b}\varphi\left(  x,y\right)  dxdy=0$ for the moment. After
an additional integration of the both sides of (\ref{D_f0(y)}), we find
\[
f_{0}\left(  y\right)  -f_{0}\left(  0\right)  =-\frac{1}{a}\int_{0}^{a}%
\int_{0}^{y}dt\int_{0}^{t}\varphi\left(  x,u\right)  dxdu
\]
Partial integration yields%
\begin{align*}
f_{0}\left(  y\right)  -f_{0}\left(  0\right)   &  =-\frac{1}{a}\int_{0}%
^{a}\left(  \left.  t\int_{0}^{t}\varphi\left(  x,u\right)  du\right\vert
_{0}^{y}-\int_{0}^{y}t\varphi\left(  x,t\right)  dt\right)  dx\\
&  =-\frac{1}{a}\int_{0}^{a}dx\int_{0}^{y}dt\left(  y-t\right)  \varphi\left(
x,t\right)
\end{align*}
In summary, we arrive for $n=0$, where the constant $f_{0}\left(  0\right)
=c$ arises due to the omission of the boundary condition $\left.  \frac
{df_{0}\left(  y\right)  }{dy}\right\vert _{y=b}=0$, at%
\begin{equation}
f_{0}\left(  y\right)  =c-\frac{1}{a}\int_{0}^{a}dx\int_{0}^{y}dt\left(
y-t\right)  \varphi\left(  x,t\right)  \label{f0(y)}%
\end{equation}

We substitute the functions $y_{1}=\cosh\frac{\pi n}{a}y$ and $y_{2}%
=\cosh\frac{\pi n}{a}\left(  b-y\right)  $ in the general solution \cite[eq.
(5.2.19)]{Morse}%
\[
\psi\left(  t\right)  =y_{1}\left(  c_{1}-\int\frac{ry_{2}}{W\left(
y_{1},y_{2}\right)  }dt\right)  +y_{2}\left(  c_{2}+\int\frac{ry_{1}}{W\left(
y_{1},y_{2}\right)  }dt\right)
\]
of a second order differential equation $\frac{d^{2}\psi}{dt^{2}}+p\left(
t\right)  \frac{d\psi}{dt}+q\left(  t\right)  \psi\left(  t\right)  =r\left(
t\right)  $, where the indefinite integrals and constants $c_{1}$ and $c_{2}$
may be adjusted to fit the boundary conditions. Applying the general theory,
with integration limits that obey the boundary conditions, we obtain%
\begin{align*}
f_{n}\left(  y\right)   &  =\cosh\frac{\pi n}{a}y\left(  c_{1}+\frac{1}%
{\frac{\pi n}{a}\sinh\left(  \frac{\pi n}{a}b\right)  }\int_{y}^{b}\rho
_{n}\left(  t\right)  \cosh\frac{\pi n}{a}\left(  b-t\right)  dt\right) \\
&  \hspace{0.5cm}+\cosh\frac{\pi n}{a}\left(  b-y\right)  \left(  c_{2}%
+\frac{1}{\frac{\pi n}{a}\sinh\left(  \frac{\pi n}{a}b\right)  }\int_{0}%
^{y}\rho_{n}\left(  t\right)  \cosh\frac{\pi n}{a}tdt\right)
\end{align*}
The boundary condition $f_{n}^{\prime}\left(  0\right)  =f_{n}^{\prime}\left(
b\right)  =0$ requires the computation of the derivative\footnote{Indeed,%
\begin{align*}
f_{n}^{\prime}\left(  y\right)   &  =\frac{\pi n}{a}\sinh\frac{\pi n}%
{a}y\left(  c_{1}+\frac{1}{\frac{\pi n}{a}\sinh\left(  \frac{\pi n}%
{a}b\right)  }\int_{y}^{b}\rho_{n}\left(  t\right)  \cosh\frac{\pi n}%
{a}\left(  b-t\right)  dt\right) \\
&  -\cosh\frac{\pi n}{a}y\left(  \frac{1}{\frac{\pi n}{a}\sinh\left(
\frac{\pi n}{a}b\right)  }\rho_{n}\left(  y\right)  \cosh\frac{\pi n}%
{a}\left(  b-y\right)  \right) \\
&  -\frac{\pi n}{a}\sinh\frac{\pi n}{a}\left(  b-y\right)  \left(  c_{2}%
+\frac{1}{\frac{\pi n}{a}\sinh\left(  \frac{\pi n}{a}b\right)  }\int_{0}%
^{y}\rho_{n}\left(  t\right)  \cosh\frac{\pi n}{a}tdt\right) \\
&  +\cosh\frac{\pi n}{a}\left(  b-y\right)  \left(  \frac{1}{\frac{\pi n}%
{a}\sinh\left(  \frac{\pi n}{a}b\right)  }\rho_{n}\left(  y\right)  \cosh
\frac{\pi n}{a}y\right)
\end{align*}
After cancellation of the second and fourth term, we find the derivative.},%
\begin{align*}
f_{n}^{\prime}\left(  y\right)   &  =\frac{\pi n}{a}\sinh\frac{\pi n}%
{a}y\left(  c_{1}+\frac{1}{\frac{\pi n}{a}\sinh\left(  \frac{\pi n}%
{a}b\right)  }\int_{y}^{b}\rho_{n}\left(  t\right)  \cosh\frac{\pi n}%
{a}\left(  b-t\right)  dt\right) \\
&  -\frac{\pi n}{a}\sinh\frac{\pi n}{a}\left(  b-y\right)  \left(  c_{2}%
+\frac{1}{\frac{\pi n}{a}\sinh\left(  \frac{\pi n}{a}b\right)  }\int_{0}%
^{y}\rho_{n}\left(  t\right)  \cosh\frac{\pi n}{a}tdt\right)
\end{align*}
The first boundary condition $f_{n}^{\prime}\left(  0\right)  =-\frac{\pi
n}{a}\sinh\frac{\pi n}{a}b\;\left(  c_{2}\right)  =0$ requires that $c_{2}=0$.
Similarly, the second, boundary condition $f_{n}^{\prime}\left(  b\right)
=\frac{\pi n}{a}\sinh\frac{\pi n}{a}b\;\left(  c_{1}\right)  =0$ requires that
$c_{1}=0$. Thus, with $c_{1}=c_{2}=0$, we obtain the solution, for $n>0$,%
\begin{equation}
f_{n}\left(  y\right)  =\frac{\int_{0}^{y}\rho_{n}\left(  t\right)  \cosh
\frac{\pi n}{a}\left(  b-y\right)  \cosh\frac{\pi n}{a}tdt+\int_{y}^{b}%
\rho_{n}\left(  t\right)  \cosh\frac{\pi n}{a}y\cosh\frac{\pi n}{a}\left(
b-t\right)  dt}{\frac{\pi n}{a}\sinh\left(  \frac{\pi n}{a}b\right)  }
\label{fn(y)_n>0}%
\end{equation}
which we write as%
\[
f_{n}\left(  y\right)  =\int_{0}^{b}g_{n}\left(  y,t\right)  \rho_{n}\left(
t\right)
\]
where%
\[
g_{n}\left(  y,t\right)  =\frac{1}{\frac{\pi n}{a}\sinh\left(  \frac{\pi n}%
{a}b\right)  }\left\{
\begin{array}
[c]{cc}%
\cosh\frac{\pi n}{a}\left(  b-y\right)  \cosh\frac{\pi n}{a}t & \text{for
}0\leq t\leq y\\
\cosh\frac{\pi n}{a}y\cosh\frac{\pi n}{a}\left(  b-t\right)  & \text{for
}y<t\leq b
\end{array}
\right.
\]
One may verify that $f_{n}\left(  y\right)  $ in (\ref{fn(y)_n>0}) satisfies
the differential equation $\frac{d^{2}f_{n}\left(  y\right)  }{dy^{2}}-\left(
\frac{\pi n}{a}\right)  ^{2}f_{n}\left(  y\right)  =-\rho_{n}\left(  y\right)
$.

Introducing $\rho_{n}\left(  y\right)  =\frac{2}{a}\int_{0}^{a}\varphi\left(
x,y\right)  \cos\left(  \frac{\pi n}{a}x\right)  dx$ into $f_{n}\left(
y\right)  =\int_{0}^{b}g_{n}\left(  y,t\right)  \rho_{n}\left(  t\right)  dt$
yields, for $n>0$,%
\[
f_{n}\left(  y\right)  =\frac{2}{a}\int_{0}^{a}d\xi\int_{0}^{b}d\eta
\;\varphi\left(  \xi,\eta\right)  g_{n}\left(  y,\eta\right)  \cos\left(
\frac{\pi n}{a}\xi\right)
\]
The case $f_{0}\left(  y\right)  =c-\frac{1}{a}\int_{0}^{a}d\xi\int_{0}%
^{y}d\eta\varphi\left(  \xi,\eta\right)  \left(  y-\eta\right)  $ in
(\ref{f0(y)}) for $n=0$ shows that%
\[
f_{0}\left(  y\right)  =c-\frac{1}{a}\int_{0}^{a}d\xi\int_{0}^{b}d\eta
\varphi\left(  \xi,\eta\right)  \left(  y-\eta\right)
\]
which defines%
\[
g_{0}\left(  y,\eta\right)  =\left\{
\begin{array}
[c]{cc}%
\eta-y & \text{for }0\leq\eta\leq y\\
0 & \text{for }y<\eta\leq b
\end{array}
\right.
\]
Combining all pieces into the proposal leads to%
\begin{align*}
F\left(  x,y\right)   &  =f_{0}\left(  y\right)  +\sum_{n=1}^{\infty}%
f_{n}\left(  y\right)  \cos\left(  \frac{\pi n}{a}x\right) \\
&  =c+\int_{0}^{a}d\xi\int_{0}^{b}d\eta\varphi\left(  \xi,\eta\right)
\frac{1}{a}g_{0}\left(  y,\eta\right)  +\int_{0}^{a}d\xi\int_{0}^{b}%
d\eta\;\varphi\left(  \xi,\eta\right)  \frac{2}{a}\sum_{n=1}^{\infty}%
g_{n}\left(  y,\eta\right)  \cos\left(  \frac{\pi n}{a}\xi\right)  \cos\left(
\frac{\pi n}{a}x\right)
\end{align*}
Green's function theory tells us that the general solution
(\ref{generalV_solution_Green_function}) can be written as%
\[
F\left(  x,y\right)  =\int_{0}^{a}d\xi\int_{0}^{b}d\eta\;\varphi\left(
\xi,\eta\right)  K\left(  x,y,\xi,\eta\right)
\]
Although a Green's function obeying a Poisson equation
(\ref{Green_diff_equation_Poisson}) with a single delta function does not
exist for our boundary problem and at least two opposite delta functions are
needed to satisfy the $n=0$ boundary condition, reflecting that the total
amount of charges in the rectangle is zero, i.e. $\left.  \frac{df_{0}\left(
y\right)  }{dy}\right\vert _{y=b}=\int_{0}^{a}d\xi\int_{0}^{b}d\eta
\varphi\left(  \xi,\eta\right)  =0$, we proceed with the suggestion that the
constant $c=0$ to arrive at the \textquotedblleft quasi\textquotedblright%
-Green's function
\begin{equation}
K\left(  x,y,\xi,\eta\right)  =\frac{1}{a}g_{0}\left(  y,\eta\right)
+\frac{2}{a}\sum_{n=1}^{\infty}g_{n}\left(  y,\eta\right)  \cos\left(
\frac{\pi n}{a}\xi\right)  \cos\left(  \frac{\pi n}{a}x\right)
\label{Green_K_single_sum}%
\end{equation}
where%
\begin{equation}
g_{n}\left(  y,\eta\right)  =\frac{1}{\frac{\pi n}{a}\sinh\left(  \frac{\pi
n}{a}b\right)  }\left\{
\begin{array}
[c]{cc}%
\cosh\frac{\pi n}{a}\left(  b-y\right)  \cosh\frac{\pi n}{a}\eta & \text{for
}0\leq\eta\leq y\\
\cosh\frac{\pi n}{a}\left(  b-\eta\right)  \cosh\frac{\pi n}{a}y & \text{for
}y<\eta\leq b
\end{array}
\right.  \label{Fourier_coeff_gn}%
\end{equation}
The \textquotedblleft quasi\textquotedblright-Green's function thus obeys the
Poisson equation
\[
\Delta K\left(  r;r_{s},r_{d}\right)  =-\left(  \delta\left(  r-r_{s}\right)
-\delta\left(  r-r_{d}\right)  \right)
\]
with a delta function at the source node $s$ and an opposite delta function at
the drain node $d$. The reversal of $y$ and $\eta$ in (\ref{Fourier_coeff_gn})
shows that $g_{n}\left(  y,\eta\right)  $ for $n>0$ in the regime $0\leq
\eta\leq y$ equals $g_{n}\left(  \eta,y\right)  $ in the interval $y<\eta\leq
b$.

Finally, due to symmetry, replacing $x$ by $y$, $\xi$ by $\eta$ and $a$ by $b$
(and vice versa) results into the dual single quasi-Green's function%
\begin{equation}
K\left(  x,y,\xi,\eta\right)  =\frac{1}{b}g_{0}\left(  x,\xi\right)  +\frac
{2}{b}\sum_{n=1}^{\infty}h_{n}\left(  x,\xi\right)  \cos\left(  \frac{\pi
n}{b}\eta\right)  \cos\left(  \frac{\pi n}{b}y\right)
\label{Green_K_single_sum_dual}%
\end{equation}
where%
\begin{equation}
h_{n}\left(  x,\xi\right)  =\frac{1}{\frac{\pi n}{b}\sinh\left(  \frac{\pi
n}{b}a\right)  }\left\{
\begin{array}
[c]{cc}%
\cosh\frac{\pi n}{b}\left(  a-x\right)  \cosh\frac{\pi n}{b}\xi & \text{for
}0\leq\xi\leq x\\
\cosh\frac{\pi n}{b}x\cosh\frac{\pi n}{b}\left(  a-\xi\right)  & \text{for
}x<\xi\leq a
\end{array}
\right.  \label{Fourier_coeff_hn}%
\end{equation}
The correctness of (\ref{Green_K_single_sum_dual}) may be verified by starting
from the proposal $F\left(  x,y\right)  =\sum_{n=0}^{\infty}r_{n}\left(
x\right)  \cos\left(  \frac{\pi n}{b}y\right)  $ in the Poisson equation
(\ref{Poisson_2D}) and mimicking the steps followed that led to the
\textquotedblleft quasi\textquotedblright-Green's function
(\ref{Green_K_single_sum}).

Numerical computations indicate that the voltage (\ref{Voltage_2D_Green}),
computed by the single sum Green's function (\ref{Green_K_single_sum}), only
differs from the voltage (\ref{Voltage_2D}), which equals the voltage
(\ref{Voltage_2D_Green}) computed by the double sum Green's function
(\ref{Green'sFunction_K_2D}), by a constant that only depends upon the
injection and ejection points $\left(  s,d\right)  $. Thus, for a given pair
of injection points $\left(  s,d\right)  $ the voltage is only defined with
respect to a reference potential. In conclusion, we may proceed with the much
faster converging single Green's function (\ref{Green_K_single_sum}).

\subsection{Convergence of the quasi-Green's function
(\ref{Green_K_single_sum}) and (\ref{Green_K_single_sum_dual})}

\label{sec_discontinuities}The Fourier coefficients $g_{n}\left(
y,\eta\right)  $ in the quasi-Green's function (\ref{Green_K_single_sum}) are
peculiar around the point $y=\eta$, which corresponds to the $y$-coordinate of
a point charge or delta function. The derivative%
\[
\frac{dg_{n}\left(  y,\eta\right)  }{dy}=\frac{1}{\sinh\left(  \frac{\pi n}%
{a}b\right)  }\left\{
\begin{array}
[c]{cc}%
-\sinh\frac{\pi n}{a}\left(  b-y\right)  \cosh\frac{\pi n}{a}\eta & \text{for
}0\leq\eta\leq y\\
\sinh\frac{\pi n}{a}y\cosh\frac{\pi n}{a}\left(  b-\eta\right)  & \text{for
}y<\eta\leq b
\end{array}
\right.
\]
indicates that the left-derivative at $y=\eta$%
\[
\lim_{y\rightarrow\eta-\varepsilon}\frac{dg_{n}\left(  y,\eta\right)  }%
{dy}=-\frac{\sinh\frac{\pi n}{a}\left(  b-\eta\right)  \cosh\frac{\pi n}%
{a}\eta}{\sinh\left(  \frac{\pi n}{a}b\right)  }%
\]
and the right-derivative at $y=\eta$%
\[
\lim_{y\rightarrow\eta+\varepsilon}\frac{dg_{n}\left(  y,\eta\right)  }%
{dy}=\frac{\cosh\frac{\pi n}{a}\left(  b-\eta\right)  \sinh\frac{\pi n}{a}%
\eta}{\sinh\left(  \frac{\pi n}{a}b\right)  }%
\]
are different, implying that the function $g_{n}\left(  y,\eta\right)  $ is
not differentiable at $y=\eta$. Although not differentiable, the function
$g_{n}\left(  y,\eta\right)  $ in (\ref{Fourier_coeff_gn}) remains continuous
at $y=\eta$,
\[
\lim_{y\rightarrow\eta}g_{n}\left(  y,\eta\right)  =g_{n}\left(  \eta
,\eta\right)  =\frac{\cosh\frac{\pi n}{a}\left(  b-\eta\right)  \cosh\frac{\pi
n}{a}\eta}{\frac{\pi n}{a}\sinh\left(  \frac{\pi n}{a}b\right)  }%
\]

Perhaps, the more important impact is the computation of the Fourier series at
$y=\eta$. The function $g_{n}\left(  y,\eta\right)  $ in
(\ref{Fourier_coeff_gn}), rewritten as%
\[
g_{n}\left(  y,\eta\right)  =\frac{1}{\frac{\pi n}{a}\left(  e^{\frac{\pi
n}{a}b}-e^{-\frac{\pi n}{a}b}\right)  }\left\{
\begin{array}
[c]{cc}%
\left(  e^{\frac{\pi n}{a}\left(  b-y\right)  }+e^{-\frac{\pi n}{a}\left(
b-y\right)  }\right)  \left(  e^{\frac{\pi n}{a}\eta}+e^{-\frac{\pi n}{a}\eta
}\right)  & \text{for }0\leq\eta<y\\
\left(  e^{\frac{\pi n}{a}\left(  b-\eta\right)  }+e^{-\frac{\pi n}{a}\left(
b-\eta\right)  }\right)  \left(  e^{\frac{\pi n}{a}y}+e^{-\frac{\pi n}{a}%
y}\right)  & \text{for }y<\eta\leq b
\end{array}
\right.
\]
simplifies, for large $n$, to
\[
g_{n}\left(  y,\eta\right)  \simeq\frac{1}{\frac{\pi n}{a}}\left(
\begin{array}
[c]{cc}%
e^{\frac{\pi n}{a}\left(  \eta-y\right)  } & \text{for }0\leq\eta<y\\
e^{\frac{\pi n}{a}\left(  y-\eta\right)  } & \text{for }y<\eta\leq b
\end{array}
\right)  =\frac{a}{\pi n}e^{-\frac{\pi n}{a}\left\vert \eta-y\right\vert }%
\]
Hence, the Fourier series (\ref{Green_K_single_sum}), for $0\leq x\leq a$ and
$0\leq y\leq b$, is approximated by
\[
K\left(  x,y,\xi,\eta\right)  \simeq\frac{1}{a}g_{0}\left(  y,\eta\right)
+\frac{2}{a}\sum_{n=1}^{n_{o}-1}g_{n}\left(  y,\eta\right)  \cos\left(
\frac{\pi n}{a}\xi\right)  \cos\left(  \frac{\pi n}{a}x\right)  +\frac{2}{\pi
}\sum_{n=n_{o}}^{\infty}\frac{e^{-\frac{\pi n}{a}\left\vert \eta-y\right\vert
}}{n}\cos\left(  \frac{\pi n}{a}\xi\right)  \cos\left(  \frac{\pi n}%
{a}x\right)
\]
The last sum%
\begin{align*}
\left\vert \sum_{n=n_{o}}^{\infty}\frac{e^{-\frac{\pi n}{a}\left\vert
\eta-y\right\vert }}{n}\cos\left(  \frac{\pi n}{a}\xi\right)  \cos\left(
\frac{\pi n}{a}x\right)  \right\vert  &  \leq\sum_{n=n_{o}}^{\infty}%
\frac{\left(  e^{-\frac{\pi}{a}\left\vert \eta-y\right\vert }\right)  ^{n}}%
{n}<\sum_{n=1}^{\infty}\frac{\left(  e^{-\frac{\pi}{a}\left\vert
\eta-y\right\vert }\right)  ^{n}}{n}\\
&  =-\log\left(  1-e^{-\frac{\pi}{a}\left\vert \eta-y\right\vert }\right)
\end{align*}
indicates that, due to the decreasing exponential factor $e^{-\frac{\pi n}%
{a}\left\vert \eta-y\right\vert }$, the convergence is fast, provided that
$\left\vert \eta-y\right\vert $ is sufficiently large. However, at $y=\eta$,
the convergence of the Fourier series is very slow; the Fourier series even
diverges for $x=\xi$ as shown below. The Fourier series
(\ref{Green_K_single_sum}) at $y=\eta$ becomes%
\[
K\left(  x,\eta,\xi,\eta\right)  \simeq\frac{2}{a}\sum_{n=1}^{n_{o}-1}%
g_{n}\left(  \eta,\eta\right)  \cos\left(  \frac{\pi n}{a}\xi\right)
\cos\left(  \frac{\pi n}{a}x\right)  +\frac{2}{\pi}\sum_{n=n_{o}}^{\infty
}\frac{\cos\left(  \frac{\pi n}{a}\xi\right)  \cos\left(  \frac{\pi n}%
{a}x\right)  }{n}%
\]
After replacing $\cos\left(  \frac{\pi n}{a}\xi\right)  \cos\left(  \frac{\pi
n}{a}x\right)  =\frac{1}{2}\cos\left(  \frac{\pi n}{a}\left(  \xi+x\right)
\right)  +\frac{1}{2}\cos\left(  \frac{\pi n}{a}\left(  \xi-x\right)  \right)
$, the last series can be approximated as%
\begin{align*}
\sum_{n=n_{o}}^{\infty}\frac{\cos\left(  \frac{\pi n}{a}\xi\right)
\cos\left(  \frac{\pi n}{a}x\right)  }{n}  &  =\frac{1}{2}\sum_{n=n_{o}%
}^{\infty}\frac{\cos\left(  \frac{\pi n}{a}\left(  \xi+x\right)  \right)  }%
{n}+\frac{1}{2}\sum_{n=n_{o}}^{\infty}\frac{\cos\left(  \frac{\pi n}{a}\left(
\xi-x\right)  \right)  }{n}\\
&  \approx-\frac{1}{2}\log\left(  4\left\vert \sin\frac{\pi\left(
\xi+x\right)  }{2a}\right\vert \left\vert \sin\frac{\pi\left(  \xi-x\right)
}{2a}\right\vert \right)
\end{align*}
where the last line has approximated $n_{0}=1$ and invoked
(\ref{cosine_Fourier_harmonic_series}) in Appendix
\ref{sec_Fourier_series_S5(A)}. The approximation clearly demonstrates that a
point charge or delta function at $\left(  x,y\right)  =\left(  \xi
,\eta\right)  $, then $K\left(  \xi,\eta,\xi,\eta\right)  \rightarrow\infty$.

In stead of approximating the sum as above, we compute the single Fourier
series \emph{exactly} in Appendix \ref{sec_q_series}. The surprise was the
discovery of the $q$-series in (\ref{quasi_Green_in_Tq_elliptic}) for the
quasi-Green's function (\ref{Green_K_single_sum}), that are intimately related
to Gaussian polynomials and the Jacobi theta functions, which are the building
system for elliptic functions as explained Appendix \ref{sec_Gaussianpoly}. As
deduced in Appendix \ref{sec_q_analog_Green}, the single sum quasi-Green's
function (\ref{Green_K_single_sum_dual}) in the interval $0\leq\eta\leq y$ can
be written, with $q=e^{-2\frac{\pi b}{a}}$ as%
\begin{equation}
K\left(  x,y,\xi,\eta\right)  =\frac{1}{a}g_{0}\left(  y,\eta\right)
+T_{q}\left(  x,y;\xi,\eta\right)  +T_{q}\left(  x,y;\xi,-\eta\right)
+T_{q}\left(  -x,y;\xi,\eta\right)  +T_{q}\left(  -x,y;\xi,-\eta\right)
\label{quasi_Green_in_Tq_elliptic}%
\end{equation}
where the $q$-function $T_{q}\left(  x,y;\xi,\eta\right)  $ is%
\begin{align}
T_{q}\left(  x,y;\xi,\eta\right)   &  =-\frac{1}{4\pi}\log\left(
1-2\cos\left(  \frac{\pi}{a}\left(  \xi+x\right)  \right)  e^{-\frac{\pi}%
{a}\left(  y-\eta\right)  }+e^{-2\frac{\pi}{a}\left(  y-\eta\right)  }\right)
\label{def_Tq_function}\\
&  \hspace{0.5cm}-\frac{1}{4\pi}\log\prod_{k=1}^{\infty}\left(  1-2\cos\left(
\frac{\pi}{a}\left(  \xi+x\right)  \right)  e^{\frac{\pi}{a}\left(
y-\eta\right)  }q^{k}+e^{2\frac{\pi}{a}\left(  y-\eta\right)  }q^{2k}\right)
\nonumber\\
&  \hspace{1cm}\times\left(  1-2\cos\left(  \frac{\pi}{a}\left(  \xi+x\right)
\right)  e^{-\frac{\pi}{a}\left(  y-\eta\right)  }q^{k}+e^{-2\frac{\pi}%
{a}\left(  y-\eta\right)  }q^{2k}\right) \nonumber
\end{align}
Due to the fast converging sum (\ref{binomial-q_series_product_in_w(z,A)}) of
the basic $q$-product in Section \ref{sec_function_w(z,A)} and Section
\ref{sec_Gaussianpoly}, the infinite product in $T_{q}\left(  x,y;\xi
,\eta\right)  $ converges very quickly and only a few factors in product
(\ref{def_Tq_function}) are needed.

\subsection{Current density in 2D}

We proceed to compute the current density in the rectangle, for $0\leq x\leq
a$ and $0\leq y\leq b$,
\[
\overrightarrow{j}\left(  x,y\right)  =\sigma\operatorname{grad}%
V=-\sigma\left(  \frac{\partial V}{\partial x},\frac{\partial V}{\partial
y}\right)
\]
which is a two dimensional vector with magnitude%
\[
\left\vert \overrightarrow{j}\left(  x,y\right)  \right\vert =\left\vert
j\left(  x,y\right)  \right\vert =\sigma\sqrt{\left(  \frac{\partial
V}{\partial x}\right)  ^{2}+\left(  \frac{\partial V}{\partial y}\right)
^{2}}%
\]

The partial derivative with respect to $x$ of the single \textquotedblleft
quasi\textquotedblright- Green's function (\ref{Green_K_single_sum}) is%
\[
\frac{\partial K\left(  x,y,\xi,\eta\right)  }{\partial x}=\frac{2\pi}{a^{2}%
}\sum_{n=1}^{\infty}n\ g_{n}\left(  y,\eta\right)  \cos\left(  \frac{\pi n}%
{a}\xi\right)  \sin\left(  \frac{\pi n}{a}x\right)
\]
However, $\frac{\partial K\left(  x,y,\xi,\eta\right)  }{\partial y}$ of the
single \textquotedblleft quasi\textquotedblright- Green's function
(\ref{Green_K_single_sum}) is problematic at $y=\eta$, where the derivative of
$g_{n}\left(  y,\eta\right)  $ does not exist. Fortunately, we have the dual
\textquotedblleft quasi\textquotedblright- Green's function
(\ref{Green_K_single_sum_dual}) to our disposal, from which we find that%
\[
\frac{\partial K\left(  x,y,\xi,\eta\right)  }{\partial y}=\frac{2\pi}{b^{2}%
}\sum_{n=1}^{\infty}n\ h_{n}\left(  x,\xi\right)  \cos\left(  \frac{\pi n}%
{b}\eta\right)  \sin\left(  \frac{\pi n}{b}y\right)
\]
The corresponding partial derivatives for the voltage are%
\[
\left\{
\begin{array}
[c]{c}%
\frac{\partial V\left(  x,y,x_{in},y_{in},x_{out},y_{out}\right)  }{\partial
x}=\frac{2\pi}{a^{2}}\sum_{n=1}^{\infty}n\ \left(  g_{n}\left(  y,y_{in}%
\right)  \cos\left(  \frac{\pi n}{a}x_{in}\right)  -g_{n}\left(
y,y_{out}\right)  \cos\left(  \frac{\pi n}{a}x_{out}\right)  \right)
\sin\left(  \frac{\pi n}{a}x\right) \\
\frac{\partial V\left(  x,y,x_{in},y_{in},x_{out},y_{out}\right)  }{\partial
y}=\frac{2\pi}{b^{2}}\sum_{n=1}^{\infty}n\ \left(  h_{n}\left(  x,x_{in}%
\right)  \cos\left(  \frac{\pi n}{b}y_{in}\right)  -h_{n}\left(
x,x_{out}\right)  \cos\left(  \frac{\pi n}{b}y_{out}\right)  \right)
\cos\left(  \frac{\pi n}{b}y\right)
\end{array}
\right.
\]
which specify the magnitude of current density%
\begin{equation}
\left\vert j\left(  x,y;\left(  x_{in},y_{in}\right)  ,\left(  x_{out}%
,y_{out}\right)  \right)  \right\vert =\sigma\sqrt{\left(  \frac{\partial
V\left(  x,y,x_{in},y_{in},x_{out},y_{out}\right)  }{\partial x}\right)
^{2}+\left(  \frac{\partial V\left(  x,y,x_{in},y_{in},x_{out},y_{out}\right)
}{\partial y}\right)  ^{2}} \label{Single_sum_current_density}%
\end{equation}

In contrast to the partial derivatives of the double Fourier series
(\ref{Voltage_2D}), the above single series converges considerably faster and
can be used for numerical computations to compute the magnitude $\left\vert
j\left(  x,y\right)  \right\vert $ of the current density. Perhaps, more
importantly, our mathematical derivation is consistent with numerical
computations of the magnitude $\left\vert j\left(  x,y\right)  \right\vert $
of the current density, demonstrating that the single Fourier series
(\ref{Green_K_single_sum}) returns precisely the same results as the double
Fourier series (\ref{Voltage_2D}).

\begin{figure}
[h]
\begin{center}
\includegraphics[
trim=0.000000in 3.576300in 0.000000in 4.955129in,
height=5.8408cm,
width=15.1281cm
]%
{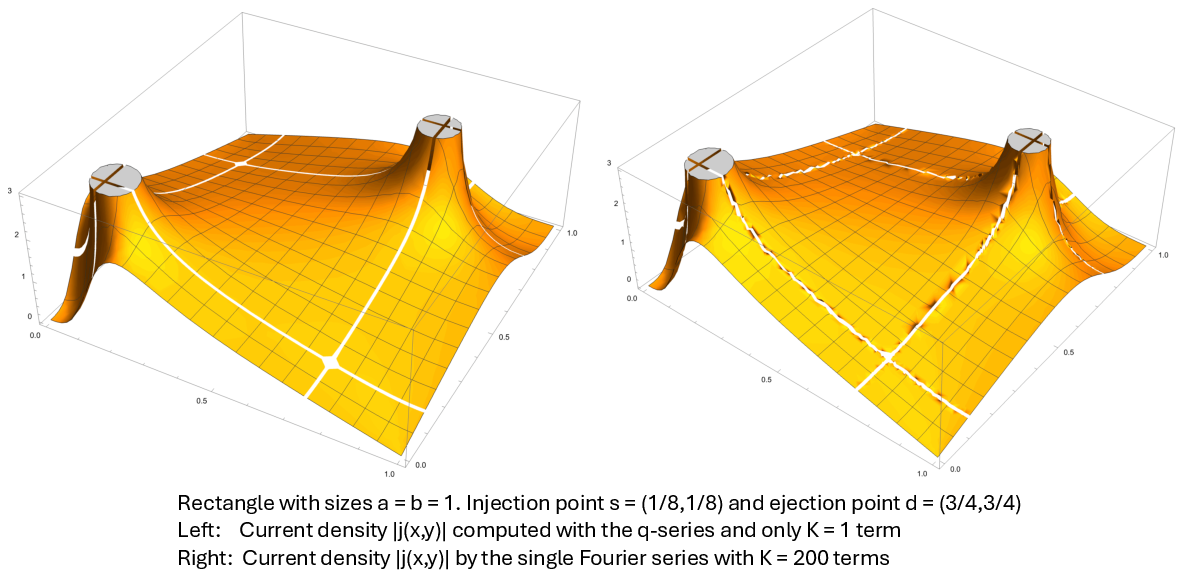}%
\caption{Again the same example computation that compares the single Fourier
series (right) computation with the $q$-series (left). Only the $k=0$ term in
the series (\ref{series_derivative_w_maal_z}) gives the same \textquotedblleft
visual\textquotedblright\ (i.e. at least two digits) accuracy as evaluations
with more terms (tested up to $k\leq K=10$)!}%
\label{Fig_currentdensity_fourier_qserie}%
\end{center}
\end{figure}

Finally, the $q$-analogue computation (\ref{partial_K_w'_function}) of
(\ref{Single_sum_current_density}) leads to the best accuracy with least
computational effort as illustrated in Fig.
\ref{Fig_currentdensity_fourier_qserie}.%
%TCIMACRO{\FRAME{fhFU}{15.1281cm}{5.8408cm}{0pt}{\Qcb{Again the same example
%computation that compares the single Fourier series (right) computation with
%the $q$-series (left). Only the $k=0$ term in the series
%(\ref{series_derivative_w_maal_z}) gives the same \textquotedblleft
%visual\textquotedblright\ (i.e. at least two digits) accuracy as evaluations
%with more terms (tested up to $k\leq K=10$)!}}%
%{\Qlb{Fig_currentdensity_fourier_qserie}}{currentdensity_fourier_qseries.ps}%
%width 15.1281cm;  height 5.8408cm;  depth 0pt;  original-width 8.2538in;
%original-height 11.6949in;  cropleft "0";  croptop "0.5763";  cropright "1";
%cropbottom "0.3058";
%filename 'currentdensity_Fourier_qSeries.ps';file-properties "XNPEU";}} }%
%BeginExpansion
%EndExpansion
Fig. \ref{Fig_currentdensity_fourier_qserie} shows that the $q$-series
computation at the white lines, where the derivative of $g_{n}\left(
y,\eta\right)  $ does not exist, leads to smooth lines. Precisely at these
lines, the Fourier series converges badly (as shown in Section
\ref{sec_discontinuities}).

\section{Summary}

Our main result is the exact and efficient computation of the potential in
(\ref{Voltage_2D_Green}), via the Green's function $K\left(  x,y,\xi
,\eta\right)  $ expressed in (\ref{quasi_Green_in_Tq_elliptic}) in terms of
the $q$-function $T_{q}\left(  x,y;\xi,\eta\right)  $ in
(\ref{def_Tq_function}) and of the magnitude $\left\vert j(x,y)\right\vert $
of the current density in (\ref{Single_sum_current_density}), where the
derivatives $\frac{\partial K\left(  x,y,\xi,\eta\right)  }{\partial x}$ and
$\frac{\partial K\left(  x,y,\xi,\eta\right)  }{\partial y}$ are computed by
(\ref{partial_K_w'_function}) together with the \textquotedblleft%
$q$-series\textquotedblright\ (\ref{series_derivative_w_maal_z}), derived in
Appendix \ref{sec_q_analog_Green_afgeleide}.

The astonishingly fast convergence of \textquotedblleft$q$%
-series\textquotedblright, in general, reduces the computation to only a few
terms in the \textquotedblleft$q$-series\textquotedblright%
\ (\ref{series_derivative_w_maal_z}). Even only one term, i.e. the $k=0$ term,
is used to produce Fig. \ref{Fig_currentdensity_fourier_qserie} and visually
no difference with more terms in the \textquotedblleft$q$%
-series\textquotedblright\ (\ref{series_derivative_w_maal_z}) is observable.
The $k=0$ term approximation provides a simple, closed analytic formula, which
is practical for quick estimates (of a quantities that are otherwise
computationally very expensive).

\begin{figure}
[h]
\begin{center}
\includegraphics[
trim=0.000000in 4.951621in 0.000000in 3.854639in,
height=5.34cm,
width=15.1281cm
]%
{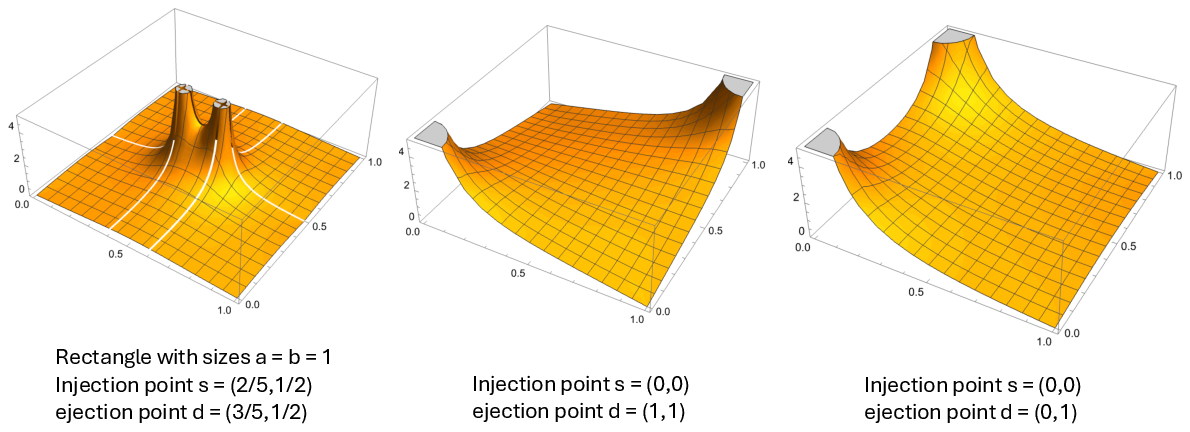}%
\caption{The magnitude $\left\vert j(x,y)\right\vert $ in
(\ref{Single_sum_current_density}), computed with (\ref{partial_K_w'_function}%
) and one term $k=0$ in the \textquotedblleft$q$-series\textquotedblright%
\ (\ref{series_derivative_w_maal_z}), for different source-destination pairs
$\left(  s,d\right)  $.}%
\label{Fig_currentdensity_fourier_qseries_verschillende_source_des}%
\end{center}
\end{figure}

Finally, we also include a drawing in Fig.
\ref{Fig_currentdensity_fourier_qseries_verschillende_source_des} of the
magnitude $\left\vert j(x,y)\right\vert $ in (\ref{Single_sum_current_density}%
), computed with (\ref{partial_K_w'_function}) and one term $k=0$ in the
\textquotedblleft$q$-series\textquotedblright%
\ (\ref{series_derivative_w_maal_z}), for different source-destination pairs
$\left(  s,d\right)  $.%
%TCIMACRO{\FRAME{fhFU}{15.1281cm}{5.34cm}{0pt}{\Qcb{The magnitude $\left\vert
%j(x,y)\right\vert $ in (\ref{Single_sum_current_density}), computed with
%(\ref{partial_K_w'_function}) and one term $k=0$ in the \textquotedblleft%
%$q$-series\textquotedblright\ (\ref{series_derivative_w_maal_z}), for
%different source-destination pairs $\left(  s,d\right)  $.}}%
%{\Qlb{Fig_currentdensity_fourier_qseries_verschillende_source_des}%
%}{currentdensity_fourier_qseries_verschillende_source_dest.ps}%
%width 15.1281cm;  height 5.34cm;  depth 0pt;  original-width 8.2538in;
%original-height 11.6949in;  cropleft "0";  croptop "0.6704";  cropright "1";
%cropbottom "0.4234";
%filename 'currentdensity_Fourier_qSeries_verschillende_source_dest.ps';file-properties "XNPEU";}%
%} }%
%BeginExpansion
%EndExpansion

\medskip\textbf{Acknowledgement}

I am very grateful to Wim van Horsen for his expert views and suggestions on
the mathematical derivations. I am supported by the European Research Council
(ERC) under the European Union's Horizon 2020 research and innovation
programme (grant agreement No 101019718).

{\footnotesize
\bibliographystyle{unsrt}
\bibliography{cac,MATH,misc,net,pvm,QTH,tel}
}

\appendix{}

\section{Elements from Fourier theory}

\subsection{Orthogonality properties}

\label{sec_orthogonality_relations}The well-known orthogonality property in
the theory of Fourier series is
\[
\int_{0}^{\pi}\cos m\theta\cos k\theta d\theta=0\text{\hspace{0.5cm}for }k\neq
m
\]
and%
\[
\int_{0}^{\pi}\cos^{2}k\theta d\theta=\pi1_{\left\{  k=0\right\}  }+\frac{\pi
}{2}1_{\left\{  k\neq0\right\}  }%
\]
as well as%
\[
\left\{
\begin{array}
[c]{ll}%
\int_{0}^{\pi}\sin m\theta\sin k\theta d\theta=0 & \text{for }k\neq m\\
\int_{-\pi}^{\pi}\sin^{2}k\theta d\theta=\pi &
\end{array}
\right.
\]
and $\int_{-\pi}^{\pi}\sin m\theta\cos k\theta d\theta=0$.

Indeed, with $\cos m\theta\cos k\theta=\frac{1}{2}\cos\left(  m-k\right)
\theta+\frac{1}{2}\cos\left(  m+k\right)  \theta$, it holds that%
\begin{align*}
\int_{0}^{\pi}\cos m\theta\cos k\theta d\theta &  =\frac{1}{2}\int_{0}^{\pi
}\cos\left(  m-k\right)  \theta d\theta+\frac{1}{2}\int_{0}^{\pi}\cos\left(
m+k\right)  \theta d\theta\\
&  =\frac{1}{2}\frac{\sin\left(  m-k\right)  \pi}{\left(  m-k\right)  }%
+\frac{1}{2}\frac{\sin\left(  m+k\right)  \pi}{\left(  m+k\right)  }%
\end{align*}
which vanishes if $k\neq m$, but $\lim_{m\rightarrow k}$ $\frac{1}{2}%
\frac{\sin\left(  m-k\right)  \pi}{\left(  m-k\right)  }=\frac{\pi}{2}$,
except if $m=k=0$. The other integrals similarly follow from $\sin m\theta\sin
k\theta=\frac{1}{2}\cos\left(  m-k\right)  \theta-\frac{1}{2}\cos\left(
m+k\right)  \theta$.

\subsection{The Fourier series $S\left(  x\right)  =\sum_{m=1}^{\infty}%
\frac{\sin xm}{m}$}

\label{sec_Fourier_series_S5(A)}The sum $S\left(  x\right)  =\sum
_{m=1}^{\infty}\frac{\sin xm}{m}=\operatorname{Im}\sum_{m=1}^{\infty}%
\frac{\left(  e^{ix}\right)  ^{m}}{m}$ is periodic in the real $x$ with
periodic $2\pi$. Clearly, $S\left(  0\right)  =S\left(  \pi\right)  =0$. Since
the Taylor series $\sum_{m=1}^{\infty}\frac{z^{m}}{m}=-\log\left(  1-z\right)
$ converges for $\left\vert z\right\vert \leq1$, except for the point $z=1$,
it holds, for $0<x<2\pi$, that%
\[
\sum_{m=1}^{\infty}\frac{e^{ixm}}{m}=-\log\left(  1-e^{ix}\right)
\]
from which $S\left(  x\right)  =\operatorname{Im}\sum_{m=1}^{\infty}%
\frac{\left(  e^{ix}\right)  ^{m}}{m}=-\operatorname{Im}\log\left(
1-e^{ix}\right)  $.

Now,
\begin{align*}
1-e^{ix}  &  =1-\cos x-i\sin x=\sqrt{\left(  1-\cos x\right)  ^{2}+\sin^{2}%
x}e^{-i\arctan\frac{\sin x}{1-\cos x}}\\
&  =\sqrt{2-2\cos x}e^{-i\arctan\frac{2\sin\frac{x}{2}\cos\frac{x}{2}}%
{2-2\cos^{2}\frac{x}{2}}}=\sqrt{2\left(  1-\left(  1-2\sin^{2}\frac{x}%
{2}\right)  \right)  }e^{-i\arctan\frac{\cos\frac{x}{2}}{\sin\frac{x}{2}}}\\
&  =2\left\vert \sin\frac{x}{2}\right\vert e^{-i\arctan\frac{1}{\tan\frac
{x}{2}}}%
\end{align*}
With%
\[
\arctan x+\arctan\frac{1}{x}=\lim_{z_{2}\rightarrow\frac{1}{x}}\arctan\left(
\frac{x+\frac{1}{x}}{1-xz_{2}}\right)  =\left\{
\begin{array}
[c]{cc}%
\frac{\pi}{2} & \text{for }x\geq0\\
-\frac{\pi}{2} & \text{for }x<0
\end{array}
\right.
\]
and the function sign$\left(  x\right)  =1$ if $x\geq0$, else sign$\left(
x\right)  =-1$, we find that%
\[
\log\left(  1-e^{ix}\right)  =\log\left(  2\left\vert \sin\frac{x}%
{2}\right\vert e^{-i\left(  \frac{\pi}{2}\text{sign}\left(  x\right)
-\frac{A}{2}\right)  }\right)  =\log\left(  2\left\vert \sin\frac{x}%
{2}\right\vert \right)  -i\left(  \frac{\pi}{2}\text{sign}\left(  x\right)
-\frac{x}{2}\right)
\]
Separating real and imaginary part in $\sum_{m=1}^{\infty}\frac{e^{ixm}}%
{m}=-\log\left(  1-e^{ix}\right)  $ shows, for $0\leq x\leq2\pi$, that
\[
S\left(  x\right)  =\sum_{m=1}^{\infty}\frac{\sin xm}{m}=\frac{\pi}%
{2}\text{sign}\left(  x\right)  -\frac{x}{2}%
\]
while, for $0<x<2\pi$,%
\begin{equation}
\sum_{m=1}^{\infty}\frac{\cos xm}{m}=-\log\left(  2\left\vert \sin\frac{x}%
{2}\right\vert \right)  \label{cosine_Fourier_harmonic_series}%
\end{equation}

Both the sine and cosine Fourier series converge very slowly, due to their
close proximity to the diverging harmonic series for $x=0$ in $\sum
_{m=1}^{\infty}\frac{\cos xm}{m}=\lim_{N\rightarrow\infty}\sum_{m=1}^{N}%
\frac{1}{m}$, where $H_{N}=\sum_{m=1}^{N}\frac{1}{m}$ is the harmonic series,
discussed in \cite[Part II]{PVM_Mittag-Leffler_Gamma}. This type of very
slowly converging series occurs in the current density derived from the double
Fourier series (\ref{Green'sFunction_K_2D}).

\section{Green's function in 2D}

\label{sec_Green_functions}The idea of Green's function (or influence function
\cite[\$14, p. 351]{Courant_HilbertI}) stems from the linearity of the Poisson
equation $\Delta V=-\varphi\left(  r\right)  $, where $r$ denotes the vector
in $D$ dimensions. If the density of charges $\varphi\left(  r\right)  $ is
\textquotedblleft piecewise continuous\textquotedblright, then $\varphi\left(
r\right)  $ can be decomposed as a continuous sum over point charges. Suppose
that $K\left(  r,r_{p}\right)  $ is a solution of the Poisson equation
\begin{equation}
\Delta K\left(  r,r_{p}\right)  =-\delta\left(  r-r_{p}\right)
\label{Green_diff_equation_Poisson}%
\end{equation}
with same boundary conditions are in the original problem, then%
\begin{equation}
V\left(  r\right)  =\int_{vol}\varphi\left(  r_{p}\right)  K\left(
r,r_{p}\right)  dr_{p} \label{generalV_solution_Green_function}%
\end{equation}
where $vol$ is the region in $D$ dimension containing the charges, formally
(but proved in \cite[p. 353-354]{Courant_HilbertI}) obeys the Poisson equation
$\Delta V=-\varphi\left(  r\right)  $, because%
\[
\Delta V\left(  r\right)  =\int_{vol}\varphi\left(  r_{p}\right)  \Delta
K\left(  r,r_{p}\right)  dr_{p}=-\int_{vol}\varphi\left(  r_{p}\right)
\delta\left(  r-r_{p}\right)  dr_{p}=-\varphi\left(  r\right)
\]
The function $K\left(  r,r_{p}\right)  $ is called the influence or Green's
function. The theory of Green's function is rich, see e.g. \cite[Chapter
7]{Zauderer}, \cite[Chapter 7]{Morse}, \cite{Duff_Naylor}. The symmetry
relation $K\left(  r,r_{p}\right)  =K\left(  r_{p},r\right)  $ is key.

Another building block in the theory of Green's functions are eigenfunctions
of the differential equation $\Delta V\left(  r\right)  +\alpha V\left(
r\right)  =-\varphi\left(  r\right)  $. Let the eigenfunction $u_{m}\left(
r\right)  $ obey $\Delta u_{m}\left(  r\right)  =\lambda_{m}w^{\prime}\left(
r\right)  u_{m}\left(  r\right)  $ as well as the boundary conditions, where
$w\left(  r\right)  $ is the associated, positive weight function \cite[art
350]{PVM_graphspectra_second_edition} with derivative $w^{\prime}\left(
r\right)  $ that defines the scalar product%
\[
\left(  f,g\right)  =\int_{vol}f\left(  r\right)  g\left(  r\right)  dw\left(
r\right)
\]
and orthogonality among eigenfunctions $\left(  u_{m},u_{k}\right)
=\delta_{km}$. Courant \& Hilbert \cite[p. 360]{Courant_HilbertI} show that
the eigenfunctions $\left\{  u_{m}\left(  r\right)  \right\}  _{m\geq1}$ form
a complete orthogonal set. Hence, the Green's function can be written as a
linear combination of eigenfunctions,%
\[
K\left(  r,r_{p}\right)  =\sum_{m=1}^{\infty}c_{m}u_{m}\left(  r\right)
\]
where $c_{m}=\left(  K\left(  r,r_{p}\right)  ,u_{m}\left(  r\right)  \right)
=\int_{vol}K\left(  r,r_{p}\right)  u_{m}\left(  r\right)  dw\left(  r\right)
$. Multiplying $\Delta K\left(  r,r_{p}\right)  =-\delta\left(  r-r_{p}%
\right)  $ by $u_{k}\left(  r\right)  $ and integrating over the region%
\[
\int_{vol}u_{k}\left(  r\right)  \Delta K\left(  r,r_{p}\right)
dr=-\int_{vol}u_{k}\left(  r\right)  \delta\left(  r-r_{p}\right)
dr=-u_{k}\left(  r_{p}\right)
\]
yields, with $K\left(  r,r_{p}\right)  =\sum_{m=1}^{\infty}c_{m}u_{m}\left(
r\right)  $, $\Delta u_{m}\left(  r\right)  =\lambda_{m}\frac{dw\left(
r\right)  }{dr}u_{m}\left(  r\right)  $ and $\left(  u_{m},u_{k}\right)
=\delta_{km}$,
\[
-u_{k}\left(  r_{p}\right)  =\int_{vol}u_{k}\left(  r\right)  \sum
_{m=1}^{\infty}c_{m}\Delta u_{m}\left(  r\right)  d\left(  r\right)
=\sum_{m=1}^{\infty}c_{m}\lambda_{m}\int_{vol}u_{k}\left(  r\right)
u_{m}\left(  r\right)  dw\left(  r\right)  =c_{k}\lambda_{k}%
\]
leading to $c_{k}=-\frac{u_{k}\left(  r_{p}\right)  }{\lambda_{k}}$. Hence, we
arrive at the eigenfunction expansion of the Green's function%
\[
K\left(  r,r_{p}\right)  =-\sum_{m=1}^{\infty}\frac{u_{m}\left(  r_{p}\right)
u_{m}\left(  r\right)  }{\lambda_{m}}%
\]
which demonstrates symmetry $K\left(  r,r_{p}\right)  =K\left(  r_{p}%
,r\right)  $ between observer at the point $r$ and the point charge at point
$r_{p}$.

\subsection{Examples}

Morse and Feshbach \cite{Morse} shows that, for large sizes $a$ and $b$ the
rectangle, the Green's function should behave as $K\left(  x,y;x_{in}%
,y_{in}\right)  \approx-2\log\left(  \sqrt{\left(  x-x_{in}\right)
^{2}+\left(  y-y_{in}\right)  ^{2}}\right)  $.

Zauderer \cite[Example 5, p. 389-392]{Zauderer} presents the solution of%
\[
\frac{\partial^{2}K}{\partial x^{2}}+\frac{\partial^{2}K}{\partial y^{2}%
}=-\delta\left(  x\,_{in},y_{in}\right)  =-\delta\left(  x-\xi\right)
\delta\left(  y-\eta\right)
\]
with Dirichlet boundary conditions $V=0$ at the boundary of the rectangle
$\left(  0,0\right)  ,\left(  0,b\right)  ,\left(  a,b\right)  $ and $\left(
a,0\right)  $. There are two representations,%

\[
K\left(  x,y;\xi,\eta\right)  =\frac{4}{ab}\sum_{n=1}^{\infty}\sum
_{m=1}^{\infty}\frac{\sin\left(  \frac{\pi n}{a}x\right)  \sin\left(
\frac{\pi n}{a}\xi\right)  \sin\left(  \frac{\pi m}{b}y\right)  \sin\left(
\frac{\pi m}{b}\eta\right)  }{\left(  \frac{\pi n}{a}\right)  ^{2}+\left(
\frac{\pi m}{b}\right)  ^{2}}%
\]
and the faster converging series \cite[p. 800]{Morse}%
\[
K\left(  x,y;\xi,\eta\right)  =\frac{2}{\pi}\sum_{k=1}^{\infty}\frac
{\sin\left(  \frac{\pi k}{a}x\right)  \sin\left(  \frac{\pi k}{a}\xi\right)
}{k\sinh\left(  \frac{\pi k}{a}b\right)  }q_{k}\left(  y\right)
\]
with%
\[
q_{k}\left(  y\right)  =\left\{
\begin{array}
[c]{cc}%
\sinh\left(  \frac{\pi k}{a}y\right)  \sinh\left(  \frac{\pi k}{a}\left(
b-\eta\right)  \right)  & \eta>y\\
\sinh\left(  \frac{\pi k}{a}\eta\right)  \sinh\left(  \frac{\pi k}{a}\left(
b-y\right)  \right)  & \eta<y
\end{array}
\right.
\]

The solution \cite[p. 399]{Zauderer}, \cite[p. 281]{Duff_Naylor} of%
\[
\frac{\partial^{2}K}{\partial x^{2}}+\frac{\partial^{2}K}{\partial y^{2}%
}=-\delta\left(  x\,_{in},y_{in}\right)  =-\delta\left(  x-\xi\right)
\delta\left(  y-\eta\right)
\]
with Neumann boundary conditions $\frac{\partial V}{\partial n}=0$ at the
boundary of the rectangle $\left(  0,0\right)  ,\left(  0,b\right)  ,\left(
a,b\right)  $ and $\left(  a,0\right)  $ is%
\begin{align}
K\left(  x,y;\xi,\eta\right)   &  =\frac{2}{ab}\sum_{n=1}^{\infty}\frac
{\cos\left(  \frac{\pi n}{a}x\right)  \cos\left(  \frac{\pi n}{a}\xi\right)
}{\left(  \frac{\pi n}{a}\right)  ^{2}}+\frac{2}{ab}\sum_{m=1}^{\infty}%
\frac{\cos\left(  \frac{\pi m}{b}y\right)  \cos\left(  \frac{\pi m}{b}%
\eta\right)  }{\left(  \frac{\pi m}{b}\right)  ^{2}}\nonumber\\
&  +\frac{4}{ab}\sum_{n=1}^{\infty}\sum_{m=1}^{\infty}\frac{\cos\left(
\frac{\pi n}{a}x\right)  \cos\left(  \frac{\pi n}{a}\xi\right)  \cos\left(
\frac{\pi m}{b}y\right)  \cos\left(  \frac{\pi m}{b}\eta\right)  }{\left(
\frac{\pi n}{a}\right)  ^{2}+\left(  \frac{\pi m}{b}\right)  ^{2}}
\label{Green'sFunction_K_2D}%
\end{align}
which is associated to (\ref{Voltage_beam}) in 3D with current injection as
$K\left(  x,y;x_{in},y_{in}\right)  -K\left(  x,y;x_{out},y_{out}\right)  $.
The shorter notation as a double Fourier series of (\ref{Green'sFunction_K_2D}%
) in \cite[p. 281]{Duff_Naylor} is%
\begin{equation}
K\left(  x,y;\xi,\eta\right)  =\frac{4}{ab}\sum_{n=0}^{\infty}\sum
_{m=0}^{\infty}\gamma_{nm}\frac{\cos\left(  \frac{\pi n}{a}x\right)
\cos\left(  \frac{\pi n}{a}\xi\right)  \cos\left(  \frac{\pi m}{b}y\right)
\cos\left(  \frac{\pi m}{b}\eta\right)  }{\left(  \frac{\pi n}{a}\right)
^{2}+\left(  \frac{\pi m}{b}\right)  ^{2}} \label{Green'sFunction_K_2D_short}%
\end{equation}
with $\gamma_{00}=0$, $\gamma_{0m}=\gamma_{n0}=\frac{1}{2}$ and $\gamma
_{nm}=1$ for $n>0$ and $m>0$. Unfortunately, it seems that a single and faster
converging series as for the boundary conditions $V=0$ above is not available,
which has motivated our computation in Section \ref{sec_Green_2D_single_sum}.

\subsection{Polar coordinates}

The Laplacian operator in polar coordinates, $x=r\cos\varphi$ and
$y=r\sin\varphi$, is\footnote{A tedious computation of the chain rule,
starting with%
\begin{align*}
\frac{\partial}{\partial r}  &  =\frac{\partial}{\partial x}\frac{\partial
x}{\partial r}+\frac{\partial}{\partial y}\frac{\partial y}{\partial r}%
=\frac{\partial}{\partial x}\cos\varphi+\frac{\partial}{\partial y}\sin
\varphi\\
\frac{\partial}{\partial\varphi}  &  =\frac{\partial}{\partial x}%
\frac{\partial x}{\partial\varphi}+\frac{\partial}{\partial y}\frac{\partial
y}{\partial\varphi}=-r\frac{\partial}{\partial x}\sin\varphi+r\frac{\partial
}{\partial y}\cos\varphi
\end{align*}
finally leads to (\ref{Laplace_polar})..}%
\begin{equation}
\frac{\partial^{2}}{\partial x^{2}}+\frac{\partial^{2}}{\partial y^{2}}%
=\frac{\partial^{2}}{\partial r^{2}}+\frac{1}{r}\frac{\partial}{\partial
r}+\frac{1}{r^{2}}\frac{\partial^{2}}{\partial\varphi^{2}}
\label{Laplace_polar}%
\end{equation}
The solution of%
\[
\frac{\partial^{2}V\left(  r,\varphi\right)  }{\partial r^{2}}+\frac{1}%
{r}\frac{\partial V\left(  r,\varphi\right)  }{\partial r}+\frac{1}{r^{2}%
}\frac{\partial^{2}V\left(  r,\varphi\right)  }{\partial\varphi^{2}}=0
\]
with separation of coordinates $V\left(  r,\varphi\right)  =R\left(  r\right)
F\left(  \varphi\right)  $ becomes%
\[
R^{\prime\prime}\left(  r\right)  F\left(  \varphi\right)  +\frac{1}%
{r}R^{\prime}\left(  r\right)  F\left(  \varphi\right)  +\frac{1}{r^{2}%
}R\left(  r\right)  F^{\prime\prime}\left(  \varphi\right)  =0
\]
Dividing by $V\left(  r,\varphi\right)  =R\left(  r\right)  F\left(
\varphi\right)  >0$ and multiplying by $r^{2}$ yields%
\[
r^{2}\frac{R^{\prime\prime}\left(  r\right)  }{R\left(  r\right)  }%
+r\frac{R^{\prime}\left(  r\right)  }{R\left(  r\right)  }+\frac
{F^{\prime\prime}\left(  \varphi\right)  }{F\left(  \varphi\right)  }=0
\]
The separate differentiations%
\[
\left\{
\begin{array}
[c]{c}%
\frac{d}{dr}\left(  r^{2}\frac{R^{\prime\prime}\left(  r\right)  }{R\left(
r\right)  }+r\frac{R^{\prime}\left(  r\right)  }{R\left(  r\right)  }\right)
=0\\
\frac{d}{d\varphi}\frac{F^{\prime\prime}\left(  \varphi\right)  }{F\left(
\varphi\right)  }=0
\end{array}
\right.
\]
lead, after integration, to%
\[
\left\{
\begin{array}
[c]{c}%
r^{2}\frac{R^{\prime\prime}\left(  r\right)  }{R\left(  r\right)  }%
+r\frac{R^{\prime}\left(  r\right)  }{R\left(  r\right)  }=\lambda_{r}\\
\frac{F^{\prime\prime}\left(  \varphi\right)  }{F\left(  \varphi\right)
}=\lambda_{\varphi}%
\end{array}
\right.
\]
where the integration constants satisfy $\lambda_{r}+\lambda_{\varphi}=0$. The
differential equation $r^{2}R^{\prime\prime}\left(  r\right)  +rR^{\prime
}\left(  r\right)  =\lambda_{r}R\left(  r\right)  $ transforms with $\frac
{d}{dr}\left(  rR^{\prime}\left(  r\right)  \right)  =rR^{\prime\prime}\left(
r\right)  +R^{\prime}\left(  r\right)  $ to $\frac{d}{dr}\left(  rR^{\prime
}\left(  r\right)  \right)  =\frac{\lambda_{r}}{r}R\left(  r\right)  $.
Further progress requires the boundary conditions.

The governing equation (\ref{Green_diff_equation_Poisson}) of Green's function
becomes%
\[
\frac{\partial^{2}K\left(  r,\varphi\right)  }{\partial r^{2}}+\frac{1}%
{r}\frac{\partial K\left(  r,\varphi\right)  }{\partial r}+\frac{1}{r^{2}%
}\frac{\partial^{2}K\left(  r,\varphi\right)  }{\partial\varphi^{2}}=-\frac
{1}{r}\delta\left(  r-r_{s}\right)  \delta\left(  \varphi-\varphi_{s}\right)
\]
because $1=\int_{-\infty}^{\infty}dx\int_{-\infty}^{\infty}dy\delta\left(
x-x_{s}\right)  \delta\left(  y-y_{s}\right)  =\int_{0}^{\infty}rdr\int
_{0}^{2\pi}\delta\varphi\left(  \frac{1}{r}\delta\left(  r-r_{s}\right)
\right)  \delta\left(  \varphi-\varphi_{s}\right)  $. If $K\left(
r,\varphi\right)  $ operates in the entire 2D plane, then Duff and Naylor
\cite[Sec. 7.2]{Duff_Naylor} argue that, by symmetry, the Green's function
$K\left(  r,\varphi\right)  $ will only depend upon $r$, the distance of
source and observer. If the source is at the origin, they pose%
\begin{equation}
\frac{d^{2}K\left(  r\right)  }{dr^{2}}+\frac{1}{r}\frac{dK\left(  r\right)
}{dr}=-\frac{1}{2\pi r}\delta\left(  r\right)  \label{dvgl_polar_punt_bron}%
\end{equation}
which reduces to%
\[
-\frac{1}{2\pi}\delta\left(  r\right)  =r\frac{d^{2}K\left(  r\right)
}{dr^{2}}+\frac{dK\left(  r\right)  }{dr}=\frac{d}{dr}\left(  r\frac{dK\left(
r\right)  }{dr}\right)
\]
Integrating both sides from $0$ to $r$ yields, with $\lim_{r\rightarrow
0}r\frac{dK\left(  r\right)  }{dr}=c$, to $r\frac{dK\left(  r\right)  }%
{dr}-c=-\frac{1}{2\pi}$ or
\[
\frac{dK\left(  r\right)  }{dr}=\frac{c-\frac{1}{2\pi}}{r}%
\]
Another integration leads, with $\widetilde{c}=c-\frac{1}{2\pi}$, to%
\begin{equation}
K\left(  r\right)  =\widetilde{c}\log r+K\left(  1\right)
\label{Green_function_entire_plane}%
\end{equation}
Duff and Naylor \cite[Sec. 7.2]{Duff_Naylor} give the particular solution
$K\left(  r\right)  =-\frac{1}{2\pi}\log r$ and remark that $\lim
_{r\rightarrow\infty}K\left(  r\right)  =-\infty$ in 2D, rather than 0 as in 3D.

\subsection{Conformal mapping}

Duff and Naylor \cite[p. 271]{Duff_Naylor} invoke conformal mapping
$w=f\left(  z\right)  $ that maps a region $D$ into the unit circle
$\left\vert w\right\vert <1$, with a given source point $s$ being mapped into
the origin $w=0$ of the circle such that $f\left(  z_{s}\right)  =0$ and
$f^{\prime}\left(  z_{s}\right)  \neq0$. The Green's function of the region
$D$ is%
\begin{equation}
K\left(  x,y\right)  =\widetilde{c}\log\left\vert f\left(  x+iy\right)
\right\vert =\widetilde{c}\operatorname{Re}\left(  \log f\left(  z\right)
\right)  \label{Green_2D_Dirichlet_boundary}%
\end{equation}
Indeed, the real (or/and imaginary) part of an analytic function is, by the
Cauchy-Riemann equations, an harmonic functions satisfying the Laplacian. At
the point $z_{s}=x_{s}+iy_{s}$, the Green's function $K\left(  x_{s}%
,y_{s}\right)  $ has a logarithmic singularity and satisfies the Poisson
equation (\ref{Green_diff_equation_Poisson}) as demonstrated by
(\ref{Green_function_entire_plane}). Also, a point $\left(  x,y\right)  $ on
the boundary $\partial D$ of the region $D$ is mapped to a point $\left\vert
w\right\vert =1$ on the unit circle, for which $K\left(  x,y\right)  =0$,
thus, satisfying the Dirichlet boundary conditions. As an example, Duff and
Naylor \cite[p. 271]{Duff_Naylor} consider as the region $D$, a unit circle
and the M\"{o}bius map $w=\frac{z-z_{s}}{z_{s}^{\ast}z-1}$, where the complex
conjugate $z_{s}^{\ast}=x_{s}-iy_{s}$, maps the disk $\left\vert z\right\vert
\leq1$ into the unit disk $\left\vert w\right\vert \leq1$ where $z_{s}$ is
mapped to the origin. Hence, the Green's function
(\ref{Green_2D_Dirichlet_boundary}) becomes%
\[
K\left(  x,y;x_{s},y_{s}\right)  =\widetilde{c}\log\left\vert \frac{z-z_{s}%
}{z_{s}^{\ast}z-1}\right\vert
\]
In the same vein, Courant and Hilbert \cite[p. 384-386]{Courant_HilbertI}
illustrate a beautiful elliptic solution (in terms of Weierstrass's $\sigma$-
function) of the two dimension partial differential equation
(\ref{Green_diff_equation_Poisson}) with Dirichlet boundary conditions $V=0$
on the boundary of the rectangle, which transforms via the Cauchy-Riemann
equations to analytic functions in the complex plane.

Our $q$-series solution (\ref{quasi_Green_in_Tq_elliptic}) together with
(\ref{def_Tq_function}) is a similar example for the Neumann boundary
conditions and quasi-Green's functions (that require two opposite delta functions).

\section{Theta-function evaluation of the Green's function
(\ref{Green'sFunction_K_2D})}

\label{sec_Theta_functions_in_Green_function}After using $\frac{1}{s}=\int
_{0}^{\infty}e^{-st}dt$ for $\operatorname{Re}\left(  s\right)  >0$, the
Green's function (\ref{Green'sFunction_K_2D}) is transformed into%
\begin{align*}
K\left(  x,y;\xi,\eta\right)   &  =\frac{2}{ab}\int_{0}^{\infty}\sum
_{n=1}^{\infty}\cos\left(  \frac{\pi n}{a}x\right)  \cos\left(  \frac{\pi
n}{a}\xi\right)  e^{-\left(  \frac{\pi n}{a}\right)  ^{2}t}dt+\frac{2}{ab}%
\int_{0}^{\infty}\sum_{m=1}^{\infty}\cos\left(  \frac{\pi m}{b}y\right)
\cos\left(  \frac{\pi m}{b}\eta\right)  e^{-\left(  \frac{\pi m}{b}\right)
^{2}t}dt\\
&  +\frac{4}{ab}\int_{0}^{\infty}\sum_{n=1}^{\infty}\cos\left(  \frac{\pi
n}{a}x\right)  \cos\left(  \frac{\pi n}{a}\xi\right)  e^{-\left(  \frac{\pi
n}{a}\right)  ^{2}t}\sum_{m=1}^{\infty}\cos\left(  \frac{\pi m}{b}y\right)
\cos\left(  \frac{\pi m}{b}\eta\right)  e^{-\left(  \frac{\pi m}{b}\right)
^{2}t}dt
\end{align*}
We write%
\[
\sum_{n=1}^{\infty}\cos\left(  \frac{\pi n}{a}x\right)  \cos\left(  \frac{\pi
n}{a}\xi\right)  e^{-\left(  \frac{\pi n}{a}\right)  ^{2}t}=\operatorname{Re}%
\left(  \sum_{n=1}^{\infty}\left(  e^{i\frac{\pi}{a}\left(  x+\xi\right)
}\right)  ^{n}e^{-\left(  \frac{\pi n}{a}\right)  ^{2}t}+\sum_{n=1}^{\infty
}\left(  e^{i\frac{\pi}{a}\left(  x-\xi\right)  }\right)  ^{n}e^{-\left(
\frac{\pi n}{a}\right)  ^{2}t}\right)
\]
The appearing series $\sum_{n=1}^{\infty}z^{n}e^{-\left(  \frac{\pi n}%
{a}\right)  ^{2}t}$ can be reformulated in terms of the theta $\theta
_{3}\left(  z|\tau\right)  =\sum_{m=-\infty}^{\infty}q^{m^{2}}e^{2m\pi iz}$
function \cite[art. 160]{Tannery_MolkII} with $q=e^{\pi i\tau}$ and
$\operatorname{Im}\left(  \tau\right)  >0$, also written as $\theta_{3}\left(
z|\tau\right)  =1+2\sum_{m=1}^{\infty}q^{m^{2}}\cos\left(  2m\pi z\right)  $,
as%
\begin{align*}
\sum_{n=1}^{\infty}\cos\left(  \frac{\pi n}{a}x\right)  \cos\left(  \frac{\pi
n}{a}\xi\right)  e^{-\left(  \frac{\pi n}{a}\right)  ^{2}t}  &  =\sum
_{n=1}^{\infty}\cos\left(  \frac{\pi n}{a}\left(  x+\xi\right)  \right)
e^{-\left(  \frac{\pi n}{a}\right)  ^{2}t}+\sum_{n=1}^{\infty}\cos\left(
\frac{\pi n}{a}\left(  x-\xi\right)  \right)  e^{-\left(  \frac{\pi n}%
{a}\right)  ^{2}t}\\
&  =\frac{1}{2}\left(  \theta_{3}\left(  \left.  \frac{\left(  x+\xi\right)
}{2a}\right\vert i\frac{\pi}{a^{2}}t\right)  +\theta_{3}\left(  \left.
\frac{\left(  x-\xi\right)  }{2a}\right\vert i\frac{\pi}{a^{2}}t\right)
\right)  -1
\end{align*}
Thus, $q=e^{-\left(  \frac{\pi}{a}\right)  ^{2}t}$. After simplification,
where the constant term $-1$ disappeared, we find that%
\begin{align}
K\left(  x,y;\xi,\eta\right)   &  =-\frac{1}{ab}\int_{0}^{\infty}\left(
\theta_{3}\left(  \left.  \frac{\left(  x+\xi\right)  }{2a}\right\vert
i\frac{\pi}{a^{2}}t\right)  +\theta_{3}\left(  \left.  \frac{\left(
x-\xi\right)  }{2a}\right\vert i\frac{\pi}{a^{2}}t\right)  \right)
dt\nonumber\\
&  \hspace{0.5cm}-\frac{1}{ab}\int_{0}^{\infty}\left(  \theta_{3}\left(
\left.  \frac{\left(  y+\eta\right)  }{2b}\right\vert i\frac{\pi}{b^{2}%
}t\right)  +\theta_{3}\left(  \left.  \frac{\left(  y-\eta\right)  }%
{2b}\right\vert i\frac{\pi}{b^{2}}t\right)  \right)  dt\nonumber\\
&  \hspace{0.5cm}+\frac{1}{ab}\int_{0}^{\infty}\left(  \theta_{3}\left(
\left.  \frac{\left(  x+\xi\right)  }{2a}\right\vert i\frac{\pi}{a^{2}%
}t\right)  +\theta_{3}\left(  \left.  \frac{\left(  x-\xi\right)  }%
{2a}\right\vert i\frac{\pi}{a^{2}}t\right)  \right) \nonumber\\
&  \hspace{1cm}\times\left(  \theta_{3}\left(  \left.  \frac{\left(
y+\eta\right)  }{2b}\right\vert i\frac{\pi}{b^{2}}t\right)  +\theta_{3}\left(
\left.  \frac{\left(  y-\eta\right)  }{2b}\right\vert i\frac{\pi}{b^{2}%
}t\right)  \right)  dt \label{Green_2D_theta_functions}%
\end{align}

This theta function approach (\ref{Green_2D_theta_functions}) may be regarded
as the more fundamental, because a rectangle can be transformed into a
half-plane by an elliptic function (see e.g. \cite{Sansone}). Any elliptic
function can be written in terms of the four theta functions. The theory of
elliptic functions is a pearl of complex function theory, established in the
19-th century, starting with Gauss's work on the arithmetic-geometric mean
\cite{PVM_AGM_Gauss}, followed by Abel, Jacobi, Weierstrass and others.

Although exact, (\ref{Green_2D_theta_functions}) requires the integral over
all positive $t$. Numerically, we found that contour plots computed from
(\ref{Green_2D_theta_functions}) with Mathematica were computationally (i.e.
in time) demanding. In addition, the partial derivatives with respect to $x$
and $y$ are not straightforwardly written in terms of the four theta-function.
Therefore, we have not further exploited (\ref{Green_2D_theta_functions}).

\section{Evaluation of the single sum quasi-Green's function
(\ref{Green_K_single_sum_dual})}

\label{sec_q_series}

\subsection{Voltage: Single sum quasi-Green's function
(\ref{Green_K_single_sum_dual})}

\label{sec_q_analog_Green}

We confine ourselves to the single sum quasi-Green's function
(\ref{Green_K_single_sum_dual}) in the interval $0\leq\eta\leq y$, where
(\ref{Fourier_coeff_gn}) indicates that
\[
K\left(  x,y,\xi,\eta\right)  =\frac{1}{a}g_{0}\left(  y,\eta\right)
+\frac{2}{\pi}\sum_{n=1}^{\infty}\frac{\cosh\frac{\pi n}{a}\left(  b-y\right)
\cosh\frac{\pi n}{a}\eta}{n\sinh\left(  \frac{\pi n}{a}b\right)  }\cos\left(
\frac{\pi n}{a}\xi\right)  \cos\left(  \frac{\pi n}{a}x\right)
\]
After invoking $2\cos a\cos b=\cos\left(  a+b\right)  +\cos\left(  a-b\right)
$, the $n$-sum is%
\begin{align*}
S  &  =\frac{1}{\pi}\sum_{n=1}^{\infty}\frac{\cosh\frac{\pi n}{a}\left(
b-y\right)  \cosh\frac{\pi n}{a}\eta}{n\sinh\left(  \frac{\pi n}{a}b\right)
}\cos\left(  \frac{\pi n}{a}\left(  \xi+x\right)  \right) \\
&  \hspace{0.5cm}+\frac{1}{\pi}\sum_{n=1}^{\infty}\frac{\cosh\frac{\pi n}%
{a}\left(  b-y\right)  \cosh\frac{\pi n}{a}\eta}{n\sinh\left(  \frac{\pi n}%
{a}b\right)  }\cos\left(  \frac{\pi n}{a}\left(  \xi-x\right)  \right)
\end{align*}
Similarly, with $2\cosh a\cosh b=\cosh\left(  a+b\right)  +\cosh\left(
a-b\right)  $, we have that%
\begin{align*}
S  &  =\frac{1}{2\pi}\sum_{n=1}^{\infty}\frac{\cosh\frac{\pi n}{a}\left(
b-y+\eta\right)  }{n\sinh\left(  \frac{\pi n}{a}b\right)  }\cos\left(
\frac{\pi n}{a}\left(  \xi+x\right)  \right)  +\\
&  \hspace{0.5cm}+\frac{1}{2\pi}\sum_{n=1}^{\infty}\frac{\cosh\frac{\pi n}%
{a}\left(  b-y-\eta\right)  }{n\sinh\left(  \frac{\pi n}{a}b\right)  }%
\cos\left(  \frac{\pi n}{a}\left(  \xi+x\right)  \right) \\
&  \hspace{0.5cm}+\frac{1}{2\pi}\sum_{n=1}^{\infty}\frac{\cosh\frac{\pi n}%
{a}\left(  b-y+\eta\right)  }{n\sinh\left(  \frac{\pi n}{a}b\right)  }%
\cos\left(  \frac{\pi n}{a}\left(  \xi-x\right)  \right) \\
&  \hspace{0.5cm}+\frac{1}{2\pi}\sum_{n=1}^{\infty}\frac{\cosh\frac{\pi n}%
{a}\left(  b-y-\eta\right)  }{n\sinh\left(  \frac{\pi n}{a}b\right)  }%
\cos\left(  \frac{\pi n}{a}\left(  \xi-x\right)  \right)
\end{align*}
Finally, with $\cosh a=\cos ia$, we find the reduction%
\begin{align*}
S  &  =\frac{1}{4\pi}\sum_{n=1}^{\infty}\frac{\cos\frac{\pi n}{a}\left(
\left(  \xi+x\right)  +i\left(  b-y+\eta\right)  \right)  }{n\sinh\left(
\frac{\pi n}{a}b\right)  }+\frac{1}{4\pi}\sum_{n=1}^{\infty}\frac{\cos
\frac{\pi n}{a}\left(  \left(  \xi+x\right)  -i\left(  b-y+\eta\right)
\right)  }{n\sinh\left(  \frac{\pi n}{a}b\right)  }\\
&  \hspace{0.5cm}+\frac{1}{4\pi}\sum_{n=1}^{\infty}\frac{\cos\frac{\pi n}%
{a}\left(  \left(  \xi+x\right)  +i\left(  b-y-\eta\right)  \right)  }%
{n\sinh\left(  \frac{\pi n}{a}b\right)  }+\frac{1}{4\pi}\sum_{n=1}^{\infty
}\frac{\cos\frac{\pi n}{a}\left(  \left(  \xi+x\right)  -i\left(
b-y-\eta\right)  \right)  }{n\sinh\left(  \frac{\pi n}{a}b\right)  }\\
&  \hspace{0.5cm}+\frac{1}{4\pi}\sum_{n=1}^{\infty}\frac{\cos\frac{\pi n}%
{a}\left(  \left(  \xi-x\right)  +i\left(  b-y+\eta\right)  \right)  }%
{n\sinh\left(  \frac{\pi n}{a}b\right)  }+\frac{1}{4\pi}\sum_{n=1}^{\infty
}\frac{\cos\frac{\pi n}{a}\left(  \left(  \xi-x\right)  -i\left(
b-y+\eta\right)  \right)  }{n\sinh\left(  \frac{\pi n}{a}b\right)  }\\
&  \hspace{0.5cm}+\frac{1}{4\pi}\sum_{n=1}^{\infty}\frac{\cos\frac{\pi n}%
{a}\left(  \left(  \xi-x\right)  +i\left(  b-y-\eta\right)  \right)  }%
{n\sinh\left(  \frac{\pi n}{a}b\right)  }+\frac{1}{4\pi}\sum_{n=1}^{\infty
}\frac{\cos\frac{\pi n}{a}\left(  \left(  \xi-x\right)  -i\left(
b-y-\eta\right)  \right)  }{n\sinh\left(  \frac{\pi n}{a}b\right)  }%
\end{align*}
With $\cos x=\frac{e^{ix}+e^{-ix}}{2}=\operatorname{Re}\left(  e^{ix}\right)
$ and $A=\frac{\pi b}{a}$, each of the sums is of the type%
\begin{equation}
w\left(  z,A\right)  =\sum_{n=1}^{\infty}\frac{z^{n}}{n\sinh\left(  An\right)
}=2\sum_{n=1}^{\infty}\frac{z^{n}}{n\left(  e^{An}-e^{-An}\right)  }
\label{w(z,A)_series_def}%
\end{equation}
with $w\left(  0,A\right)  =0$. Thus, the single sum quasi-Green's function
(\ref{Green_K_single_sum_dual}) in the interval $0\leq\eta\leq y$ is
represented as%
\begin{align*}
K\left(  x,y,\xi,\eta\right)   &  =\frac{1}{a}g_{0}\left(  y,\eta\right)
+\frac{\operatorname{Re}\left\{  w\left(  e^{i\frac{\pi}{a}\left(
\xi+x\right)  }e^{-\frac{\pi}{a}\left(  b-y+\eta\right)  },\frac{\pi b}%
{a}\right)  +w\left(  e^{i\frac{\pi}{a}\left(  \xi+x\right)  }e^{\frac{\pi}%
{a}\left(  b-y+\eta\right)  },\frac{\pi b}{a}\right)  \right\}  }{4\pi}\\
&  \hspace{0.5cm}+\frac{\operatorname{Re}\left\{  w\left(  e^{i\frac{\pi}%
{a}\left(  \xi+x\right)  }e^{-\frac{\pi}{a}\left(  b-y-\eta\right)  }%
,\frac{\pi b}{a}\right)  +w\left(  e^{i\frac{\pi}{a}\left(  \xi+x\right)
}e^{\frac{\pi}{a}\left(  b-y-\eta\right)  },\frac{\pi b}{a}\right)  \right\}
}{4\pi}\\
&  \hspace{0.5cm}+\frac{\operatorname{Re}\left\{  w\left(  e^{i\frac{\pi}%
{a}\left(  \xi-x\right)  }e^{-\frac{\pi}{a}\left(  b-y+\eta\right)  }%
,\frac{\pi b}{a}\right)  +w\left(  e^{i\frac{\pi}{a}\left(  \xi-x\right)
}e^{\frac{\pi}{a}\left(  b-y+\eta\right)  },\frac{\pi b}{a}\right)  \right\}
}{4\pi}\\
&  \hspace{0.5cm}+\frac{\operatorname{Re}\left\{  w\left(  e^{i\frac{\pi}%
{a}\left(  \xi-x\right)  }e^{-\frac{\pi}{a}\left(  b-y-\eta\right)  }%
,\frac{\pi b}{a}\right)  +w\left(  e^{i\frac{\pi}{a}\left(  \xi-x\right)
}e^{\frac{\pi}{a}\left(  b-y-\eta\right)  },\frac{\pi b}{a}\right)  \right\}
}{4\pi}%
\end{align*}
where the function $w\left(  z,A\right)  $ is studied in Appendix
\ref{sec_function_w(z,A)}.

In the sequel, we simplify the real part and present, perhaps, the simplest
possible form in (\ref{def_Tq_function}) below. With $w\left(  z,A\right)
=-2\log\prod_{k=0}^{\infty}\left(  1-ze^{-A}q^{k}\right)  $ and $q=e^{-2A}$ in
(\ref{w_log_product}), taking into account that $0\leq\eta\leq y$,%
\begin{align*}
\operatorname{Re}w\left(  e^{i\frac{\pi}{a}\left(  \xi+x\right)  }%
e^{-\frac{\pi}{a}\left(  b-y+\eta\right)  },\frac{\pi b}{a}\right)   &
=-2\operatorname{Re}\left(  \log\prod_{k=0}^{\infty}\left(  1-e^{i\frac{\pi
}{a}\left(  \xi+x\right)  }e^{-\frac{\pi}{a}\left(  b-y+\eta\right)
}e^{-\frac{\pi b}{a}}\left(  e^{-\frac{2\pi b}{a}}\right)  ^{k}\right)
\right) \\
&  =-2\operatorname{Re}\left(  \log\prod_{k=0}^{\infty}\left(  1-e^{i\frac
{\pi}{a}\left(  \xi+x\right)  }e^{\frac{\pi}{a}\left(  y-\eta\right)  }\left(
e^{-\frac{2\pi b}{a}}\right)  ^{k+1}\right)  \right) \\
&  =-2\log\prod_{k=1}^{\infty}\left\vert 1-e^{i\frac{\pi}{a}\left(
\xi+x\right)  }e^{\frac{\pi}{a}\left(  y-\eta\right)  }\left(  e^{-\frac{2\pi
b}{a}}\right)  ^{k}\right\vert
\end{align*}
and%
\begin{align*}
\left\vert 1-e^{i\frac{\pi}{a}\left(  \xi+x\right)  }e^{\frac{\pi}{a}\left(
y-\eta\right)  }\left(  e^{-\frac{2\pi b}{a}}\right)  ^{k}\right\vert  &
=\left\vert 1-e^{i\frac{\pi}{a}\left(  \xi+x\right)  }e^{\frac{\pi}{a}\left(
y-\eta\right)  }q^{k}\right\vert \\
&  =\left\vert 1-\cos\left(  \frac{\pi}{a}\left(  \xi+x\right)  \right)
e^{\frac{\pi}{a}\left(  y-\eta\right)  }q^{k}-i\sin\left(  \frac{\pi}%
{a}\left(  \xi+x\right)  \right)  e^{\frac{\pi}{a}\left(  y-\eta\right)
}q^{k}\right\vert \\
&  =\sqrt{\left(  1-\cos\left(  \frac{\pi}{a}\left(  \xi+x\right)  \right)
e^{\frac{\pi}{a}\left(  y-\eta\right)  }q^{k}\right)  ^{2}+\left(  \sin\left(
\frac{\pi}{a}\left(  \xi+x\right)  \right)  e^{\frac{\pi}{a}\left(
y-\eta\right)  }q^{k}\right)  ^{2}}\\
&  =\sqrt{1-2\cos\left(  \frac{\pi}{a}\left(  \xi+x\right)  \right)
e^{\frac{\pi}{a}\left(  y-\eta\right)  }q^{k}+e^{2\frac{\pi}{a}\left(
y-\eta\right)  }q^{2k}}%
\end{align*}
Thus,%
\[
\operatorname{Re}w\left(  e^{i\frac{\pi}{a}\left(  \xi+x\right)  }%
e^{-\frac{\pi}{a}\left(  b-y+\eta\right)  },\frac{\pi b}{a}\right)
=-\log\prod_{k=1}^{\infty}\left(  1-2\cos\left(  \frac{\pi}{a}\left(
\xi+x\right)  \right)  e^{\frac{\pi}{a}\left(  y-\eta\right)  }q^{k}%
+e^{2\frac{\pi}{a}\left(  y-\eta\right)  }q^{2k}\right)
\]
On the other hand,%
\begin{align*}
\operatorname{Re}w\left(  e^{i\frac{\pi}{a}\left(  \xi+x\right)  }e^{\frac
{\pi}{a}\left(  b-y+\eta\right)  },\frac{\pi b}{a}\right)   &
=-2\operatorname{Re}\left(  \log\prod_{k=0}^{\infty}\left(  1-e^{i\frac{\pi
}{a}\left(  \xi+x\right)  }e^{\frac{\pi}{a}\left(  b-y+\eta\right)  }%
e^{-\frac{\pi b}{a}}\left(  e^{-\frac{2\pi b}{a}}\right)  ^{k}\right)  \right)
\\
&  =-2\operatorname{Re}\left(  \log\prod_{k=0}^{\infty}\left(  1-e^{i\frac
{\pi}{a}\left(  \xi+x\right)  }e^{-\frac{\pi}{a}\left(  y-\eta\right)
}\left(  e^{-\frac{2\pi b}{a}}\right)  ^{k}\right)  \right) \\
&  =-2\log\prod_{k=0}^{\infty}\left\vert 1-e^{i\frac{\pi}{a}\left(
\xi+x\right)  }e^{-\frac{\pi}{a}\left(  y-\eta\right)  }\left(  e^{-\frac{2\pi
b}{a}}\right)  ^{k}\right\vert
\end{align*}
is%
\[
\operatorname{Re}w\left(  e^{i\frac{\pi}{a}\left(  \xi+x\right)  }e^{\frac
{\pi}{a}\left(  b-y+\eta\right)  },\frac{\pi b}{a}\right)  =-\log\prod
_{k=0}^{\infty}\left(  1-2\cos\left(  \frac{\pi}{a}\left(  \xi+x\right)
\right)  e^{-\frac{\pi}{a}\left(  y-\eta\right)  }q^{k}+e^{-2\frac{\pi}%
{a}\left(  y-\eta\right)  }q^{2k}\right)
\]
Hence,%
\begin{align*}
T_{q}\left(  x,y;\xi,\eta\right)   &  =\frac{\operatorname{Re}\left\{
w\left(  e^{i\frac{\pi}{a}\left(  \xi+x\right)  }e^{-\frac{\pi}{a}\left(
b-y+\eta\right)  },\frac{\pi b}{a}\right)  +w\left(  e^{i\frac{\pi}{a}\left(
\xi+x\right)  }e^{\frac{\pi}{a}\left(  b-y+\eta\right)  },\frac{\pi b}%
{a}\right)  \right\}  }{4\pi}\\
&  =-\frac{1}{4\pi}\log\prod_{k=1}^{\infty}\left(  1-2\cos\left(  \frac{\pi
}{a}\left(  \xi+x\right)  \right)  e^{\frac{\pi}{a}\left(  y-\eta\right)
}q^{k}+e^{2\frac{\pi}{a}\left(  y-\eta\right)  }q^{2k}\right) \\
&  \hspace{0.5cm}-\frac{1}{4\pi}\log\prod_{k=0}^{\infty}\left(  1-2\cos\left(
\frac{\pi}{a}\left(  \xi+x\right)  \right)  e^{-\frac{\pi}{a}\left(
y-\eta\right)  }q^{k}+e^{-2\frac{\pi}{a}\left(  y-\eta\right)  }q^{2k}\right)
\end{align*}
which leads to (\ref{def_Tq_function}). The argument of the log in the first
term of (\ref{def_Tq_function}) is%
\begin{align*}
1-2\cos\left(  \frac{\pi}{a}\left(  \xi\pm x\right)  \right)  e^{-\frac{\pi
}{a}\left(  y-\eta\right)  }+e^{-2\frac{\pi}{a}\left(  y-\eta\right)  }  &
=1-2e^{-\frac{\pi}{a}\left(  y-\eta\right)  }\left(  1-2\sin^{2}\left(
\frac{\pi}{2a}\left(  \xi\pm x\right)  \right)  \right)  +e^{-2\frac{\pi}%
{a}\left(  y-\eta\right)  }\\
&  =1-2e^{-\frac{\pi}{a}\left(  y-\eta\right)  }+4e^{-\frac{\pi}{a}\left(
y-\eta\right)  }\sin^{2}\left(  \frac{\pi}{2a}\left(  \xi\pm x\right)
\right)  +e^{-2\frac{\pi}{a}\left(  y-\eta\right)  }\\
&  =\left(  1-e^{-\frac{\pi}{a}\left(  y-\eta\right)  }\right)  ^{2}+\left(
2e^{-\frac{\pi}{2a}\left(  y-\eta\right)  }\sin\left(  \frac{\pi}{2a}\left(
\xi\pm x\right)  \right)  \right)  ^{2}%
\end{align*}
which more clearly indicates that $e^{-\frac{\pi}{a}\left(  y-\eta\right)
}\leq1$ in the interval $0\leq\eta\leq y$ and only equal to 1 if $y=\eta$.

\subsection{Current density}

\label{sec_q_analog_Green_afgeleide}In the interval $0\leq\eta\leq y$, the
partial derivative of the single sum quasi-Green's function
(\ref{Green_K_single_sum_dual}) is
\[
\frac{\partial K\left(  x,y,\xi,\eta\right)  }{\partial x}=\frac{2}{a}%
\sum_{n=1}^{\infty}\frac{\cosh\frac{\pi n}{a}\left(  b-y\right)  \cosh
\frac{\pi n}{a}\eta}{\sinh\left(  \frac{\pi n}{a}b\right)  }\cos\left(
\frac{\pi n}{a}\xi\right)  \sin\left(  \frac{\pi n}{a}x\right)
\]
We repeat the method of Section \ref{sec_q_analog_Green}. After invoking
$2\sin a\cos b=\sin\left(  a+b\right)  +\sin\left(  a-b\right)  $, the $n$-sum
is%
\begin{align*}
S^{\prime}  &  =\frac{1}{a}\sum_{n=1}^{\infty}\frac{\cosh\frac{\pi n}%
{a}\left(  b-y\right)  \cosh\frac{\pi n}{a}\eta}{\sinh\left(  \frac{\pi n}%
{a}b\right)  }\sin\left(  \frac{\pi n}{a}\left(  \xi+x\right)  \right) \\
&  \hspace{0.5cm}+\frac{1}{a}\sum_{n=1}^{\infty}\frac{\cosh\frac{\pi n}%
{a}\left(  b-y\right)  \cosh\frac{\pi n}{a}\eta}{\sinh\left(  \frac{\pi n}%
{a}b\right)  }\sin\left(  \frac{\pi n}{a}\left(  \xi-x\right)  \right)
\end{align*}
Similarly, with $2\cosh a\cosh b=\cosh\left(  a+b\right)  +\cosh\left(
a-b\right)  $, we have that%
\begin{align*}
S^{\prime}  &  =\frac{1}{2a}\sum_{n=1}^{\infty}\frac{\cosh\frac{\pi n}%
{a}\left(  b-y+\eta\right)  }{\sinh\left(  \frac{\pi n}{a}b\right)  }%
\sin\left(  \frac{\pi n}{a}\left(  \xi+x\right)  \right)  +\\
&  \hspace{0.5cm}+\frac{1}{2a}\sum_{n=1}^{\infty}\frac{\cosh\frac{\pi n}%
{a}\left(  b-y-\eta\right)  }{n\sinh\left(  \frac{\pi n}{a}b\right)  }%
\sin\left(  \frac{\pi n}{a}\left(  \xi+x\right)  \right) \\
&  \hspace{0.5cm}+\frac{1}{2a}\sum_{n=1}^{\infty}\frac{\cosh\frac{\pi n}%
{a}\left(  b-y+\eta\right)  }{\sinh\left(  \frac{\pi n}{a}b\right)  }%
\sin\left(  \frac{\pi n}{a}\left(  \xi-x\right)  \right) \\
&  \hspace{0.5cm}+\frac{1}{2a}\sum_{n=1}^{\infty}\frac{\cosh\frac{\pi n}%
{a}\left(  b-y-\eta\right)  }{\sinh\left(  \frac{\pi n}{a}b\right)  }%
\sin\left(  \frac{\pi n}{a}\left(  \xi-x\right)  \right)
\end{align*}
Finally, with $\cosh a=\cos ia$ and again $2\sin a\cos b=\sin\left(
a+b\right)  +\sin\left(  a-b\right)  $, we find the reduction%
\begin{align*}
S  &  =\frac{1}{4a}\sum_{n=1}^{\infty}\frac{\sin\frac{\pi n}{a}\left(  \left(
\xi+x\right)  +i\left(  b-y+\eta\right)  \right)  }{\sinh\left(  \frac{\pi
n}{a}b\right)  }+\frac{1}{4a}\sum_{n=1}^{\infty}\frac{\sin\frac{\pi n}%
{a}\left(  \left(  \xi+x\right)  -i\left(  b-y+\eta\right)  \right)  }%
{\sinh\left(  \frac{\pi n}{a}b\right)  }\\
&  \hspace{0.5cm}+\frac{1}{4a}\sum_{n=1}^{\infty}\frac{\sin\frac{\pi n}%
{a}\left(  \left(  \xi+x\right)  +i\left(  b-y-\eta\right)  \right)  }%
{\sinh\left(  \frac{\pi n}{a}b\right)  }+\frac{1}{4a}\sum_{n=1}^{\infty}%
\frac{\sin\frac{\pi n}{a}\left(  \left(  \xi+x\right)  -i\left(
b-y-\eta\right)  \right)  }{\sinh\left(  \frac{\pi n}{a}b\right)  }\\
&  \hspace{0.5cm}+\frac{1}{4a}\sum_{n=1}^{\infty}\frac{\sin\frac{\pi n}%
{a}\left(  \left(  \xi-x\right)  +i\left(  b-y+\eta\right)  \right)  }%
{\sinh\left(  \frac{\pi n}{a}b\right)  }+\frac{1}{4a}\sum_{n=1}^{\infty}%
\frac{\sin\frac{\pi n}{a}\left(  \left(  \xi-x\right)  -i\left(
b-y+\eta\right)  \right)  }{\sinh\left(  \frac{\pi n}{a}b\right)  }\\
&  \hspace{0.5cm}+\frac{1}{4a}\sum_{n=1}^{\infty}\frac{\sin\frac{\pi n}%
{a}\left(  \left(  \xi-x\right)  +i\left(  b-y-\eta\right)  \right)  }%
{\sinh\left(  \frac{\pi n}{a}b\right)  }+\frac{1}{4a}\sum_{n=1}^{\infty}%
\frac{\sin\frac{\pi n}{a}\left(  \left(  \xi-x\right)  -i\left(
b-y-\eta\right)  \right)  }{\sinh\left(  \frac{\pi n}{a}b\right)  }%
\end{align*}
With $\sin x=\frac{e^{ix}-e^{-ix}}{2i}=\operatorname{Im}\left(  e^{ix}\right)
$, $A=\frac{\pi b}{a}$ and the definition (\ref{w(z,A)_series_def}) of the
function $w\left(  z,A\right)  =\sum_{n=1}^{\infty}\frac{z^{n}}{n\sinh\left(
An\right)  }$, each of the sums is of the type%
\[
z\frac{\partial w\left(  z,A\right)  }{\partial z}=\sum_{n=1}^{\infty}%
\frac{z^{n}}{\sinh\left(  An\right)  }%
\]
which also equals, after derivation of (\ref{w(z,A)_series_logs}),%
\begin{equation}
W\left(  z,A\right)  =z\frac{\partial w\left(  z,A\right)  }{\partial
z}=zw^{\prime}\left(  z,A\right)  =2\sum_{k=0}^{\infty}\frac{ze^{-A}q^{k}%
}{1-ze^{-A}q^{k}}=2\sum_{k=0}^{\infty}\frac{z}{e^{A\left(  1+2k\right)  }-z}
\label{series_derivative_w_maal_z}%
\end{equation}
Thus, the partial $x$-derivative of single sum quasi-Green's function
(\ref{Green_K_single_sum_dual}) in the interval $0\leq\eta\leq y$ is%
\begin{align}
\frac{\partial K\left(  x,y,\xi,\eta\right)  }{\partial x}  &  =\frac
{\operatorname{Im}\left\{  W\left(  e^{i\frac{\pi}{a}\left(  \xi+x\right)
}e^{-\frac{\pi}{a}\left(  b-y+\eta\right)  },\frac{\pi b}{a}\right)  +W\left(
e^{i\frac{\pi}{a}\left(  \xi+x\right)  }e^{\frac{\pi}{a}\left(  b-y+\eta
\right)  },\frac{\pi b}{a}\right)  \right\}  }{4a}\nonumber\\
&  \hspace{0.5cm}+\frac{\operatorname{Im}\left\{  W\left(  e^{i\frac{\pi}%
{a}\left(  \xi+x\right)  }e^{-\frac{\pi}{a}\left(  b-y-\eta\right)  }%
,\frac{\pi b}{a}\right)  +W\left(  e^{i\frac{\pi}{a}\left(  \xi+x\right)
}e^{\frac{\pi}{a}\left(  b-y-\eta\right)  },\frac{\pi b}{a}\right)  \right\}
}{4a}\nonumber\\
&  \hspace{0.5cm}+\frac{\operatorname{Im}\left\{  W\left(  e^{i\frac{\pi}%
{a}\left(  \xi-x\right)  }e^{-\frac{\pi}{a}\left(  b-y+\eta\right)  }%
,\frac{\pi b}{a}\right)  +W\left(  e^{i\frac{\pi}{a}\left(  \xi-x\right)
}e^{\frac{\pi}{a}\left(  b-y+\eta\right)  },\frac{\pi b}{a}\right)  \right\}
}{4a}\nonumber\\
&  \hspace{0.5cm}+\frac{\operatorname{Im}\left\{  W\left(  e^{i\frac{\pi}%
{a}\left(  \xi-x\right)  }e^{-\frac{\pi}{a}\left(  b-y-\eta\right)  }%
,\frac{\pi b}{a}\right)  +W\left(  e^{i\frac{\pi}{a}\left(  \xi-x\right)
}e^{\frac{\pi}{a}\left(  b-y-\eta\right)  },\frac{\pi b}{a}\right)  \right\}
}{4a} \label{partial_K_w'_function}%
\end{align}
For complex $z=X+iY$, it holds that%
\[
w^{\prime}\left(  X+iY,A\right)  =2\sum_{k=0}^{\infty}\frac{1}{e^{A\left(
1+2k\right)  }-X-iY}=2\sum_{k=0}^{\infty}\frac{e^{A\left(  1+2k\right)
}-X+iY}{\left(  e^{A\left(  1+2k\right)  }-X\right)  ^{2}+Y^{2}}%
\]
and%
\begin{align*}
W\left(  X+iY,A\right)   &  =\left(  X+iY\right)  w^{\prime}\left(
X+iY,A\right)  =2\sum_{k=0}^{\infty}\frac{\left(  \left(  e^{A\left(
1+2k\right)  }-X\right)  +iY\right)  \left(  X+iY\right)  }{\left(
e^{A\left(  1+2k\right)  }-X\right)  ^{2}+Y^{2}}\\
&  =2\sum_{k=0}^{\infty}\frac{e^{A\left(  1+2k\right)  }-\left(  X^{2}%
+Y^{2}\right)  +iYe^{A\left(  1+2k\right)  }}{\left(  e^{A\left(  1+2k\right)
}-X\right)  ^{2}+Y^{2}}%
\end{align*}
from which it follows that
\[
\operatorname{Im}\left(  W\left(  X+iY,A\right)  \right)  =2\sum_{k=0}%
^{\infty}\frac{Ye^{A\left(  1+2k\right)  }}{\left(  e^{A\left(  1+2k\right)
}-X\right)  ^{2}+Y^{2}}%
\]
allowing the computation of $\frac{\partial K\left(  x,y,\xi,\eta\right)
}{\partial x}$ in (\ref{partial_K_w'_function}).

Since the rectangle can be considered as the fundamental parallelogram $P$,
whose repetition covers the entire complex plane, the real functions
$h_{\operatorname{Re}}\left(  x,y\right)  =\frac{\partial K\left(
x,y,x_{in},y_{in}\right)  }{\partial x}-\frac{\partial K\left(  x,y,x_{out}%
,y_{out}\right)  }{\partial x}$ and $h_{\operatorname{Im}}\left(  x,y\right)
=\frac{\partial K\left(  x,y,x_{in},y_{in}\right)  }{\partial y}%
-\frac{\partial K\left(  x,y,x_{out},y_{out}\right)  }{\partial y}$ are the
real and imaginary part of a complex function $h\left(  z\right)  $ in the
complex number $z=x+iy,$%
\[
h\left(  z\right)  =h\left(  x+iy\right)  =h_{\operatorname{Re}}\left(
x,y\right)  +ih_{\operatorname{Im}}\left(  x,y\right)
\]
that is zero on the boundaries of the rectangle. Analytic functions achieve
their extrema at the boundary. The complex function $h\left(  z\right)  $ is
thus double-periodic in the complex complex plane and is an elliptic function
\cite[Chapter I \& II]{Chandrasekharan}. Indeed, $h\left(  z\right)  $ is not
an entire function, for it would be bounded in the rectangle and due to
periodicity, $h\left(  z\right)  $ would be bounded over the entire complex
plane in which case Liouville's theorem then indicates that $h\left(
z\right)  $ is a constant. Second, $h\left(  z\right)  $ must have an even
number of poles and zeros in the rectangle, because $\int_{P}h\left(
z\right)  dz=0$. The derivative of an elliptic function is also elliptic, but
its integral is not necessarily an elliptic function \cite{Chandrasekharan}.

\section{The effective resistance in the continuum RGG}

\label{sec_effective_resistance_2D_3D}The corresponding resistance of the 3D
beam, due to current injection at point $s=\left(  x_{in},y_{in},0\right)  $
and ejection at point $d=\left(  x_{out},y_{out},0\right)  $ follows from the
law of Ohm with the voltage (\ref{Voltage_beam}) as%
\begin{align}
R_{sd}  &  =\frac{V\left(  x_{in},y_{in},0\right)  -V\left(  x_{out}%
,y_{out},0\right)  }{I}\label{resistance_3D}\\
&  =\frac{4\varrho}{ab}\sum_{k=0}^{\infty}\sum_{l=0}^{\infty}\gamma_{kl}%
\frac{\left(  \cos\frac{k\pi}{a}x_{in}\cos\frac{l\pi}{b}y_{in}-\cos\frac{k\pi
}{a}x_{out}\cos\frac{l\pi}{b}y_{out}\right)  ^{2}}{\sqrt{\left(  \frac{k\pi
}{a}\right)  ^{2}+\left(  \frac{l\pi}{b}\right)  ^{2}}}\coth\left(
c\sqrt{\left(  \frac{k\pi}{a}\right)  ^{2}+\left(  \frac{l\pi}{b}\right)
^{2}}\right) \nonumber
\end{align}
where $\gamma_{00}=0,\gamma_{k0}=\gamma_{0l}=\frac{1}{2}$ and $\gamma_{kl}=1$
for $k>0$ and $l>0$; a notation that is also used in
(\ref{Green'sFunction_K_2D_short}). The last relation shows that the
resistance is always non-negative, i.e. $R_{sd}\geq0$.

The resistance $R_{sd}$ is the proportionality constant that connects the
voltage difference between the injection point $s$ and the ejection point $d$
and the current $I$ in the physical 3D world, due to Ohm's law. Choosing a
unit current of $1$ Ampere with its associated voltage difference simplifies
the above definition of the \textquotedblleft effective\textquotedblright%
\ resistance between point $s$ and $d$ to $\omega_{sd}=\left.  R_{sd}%
\right\vert _{I=1\text{ Ampere}}$. In particular in the limit towards 2D,
which is not physical\footnote{Although electrons are tiny, they live in 3D
and do not exist in 2D$.$}, but a mathematical generalization, we are led to
define the \textquotedblleft effective\textquotedblright\ resistance between
node $s$ and node $d$ in the continuum RGG as%
\begin{equation}
\omega_{sd}=V\left(  x_{in},y_{in}\right)  -V\left(  x_{out},y_{out}\right)
\label{effective_resistance_pq}%
\end{equation}
The definition (\ref{effective_resistance_pq}) circumvents difficulties due
the renormalization $\frac{\varrho I}{c}=1$ in Section \ref{sec_potential_2D}.
Thus, similarly as in Section \ref{sec_potential_2D}, we find, either from the
3D resistance (\ref{resistance_3D}) with the renormalization $\frac{\varrho
I}{c}=1$ or from (\ref{effective_resistance_pq}) and (\ref{Voltage_2D}), the
double Fourier series for \textquotedblleft effective\textquotedblright%
\ resistance between node $s$ and node $d$ in the continuum RGG%
\begin{align}
\omega_{sd}  &  =\frac{2}{ab}\lim_{K\rightarrow\infty}\sum_{k=1}^{K}%
\frac{\left(  \cos\frac{k\pi}{a}x_{in}\;-\cos\frac{k\pi}{a}x_{out}\;\right)
^{2}}{\left(  \frac{k\pi}{a}\right)  ^{2}}+\frac{2}{ab}\sum_{l=1}^{K}%
\frac{\left(  \cos\frac{l\pi}{b}y_{in}-\cos\frac{l\pi}{b}y_{out}\right)  ^{2}%
}{\left(  \frac{l\pi}{b}\right)  ^{2}}\nonumber\\
&  \hspace{0.5cm}+\frac{4}{ab}\sum_{k=1}^{K}\sum_{l=1}^{K}\frac{\left(
\cos\frac{k\pi}{a}x_{in}\;\cos\frac{l\pi}{b}y_{in}-\cos\frac{k\pi}{a}%
x_{out}\;\cos\frac{l\pi}{b}y_{out}\right)  ^{2}}{\left(  \frac{k\pi}%
{a}\right)  ^{2}+\left(  \frac{l\pi}{b}\right)  ^{2}}
\label{effective_resistance_RGG_continuun}%
\end{align}

\subsection{The effective\ resistance $\omega_{sd}$ in
(\ref{effective_resistance_RGG_continuun}) does not exist}

\label{sec_effective_res_Eisenstein}The single sums in
(\ref{effective_resistance_RGG_continuun}) converge (and can be written in
terms of the polylogarithm function Li$_{k}\left(  z\right)  =\sum
_{n=1}^{\infty}\frac{z^{n}}{n^{k}}$), but we will demonstrate here that the
double sum $S$ in (\ref{effective_resistance_RGG_continuun}) diverges.

Indeed, all the terms are non-negative. Moreover, only a finite number of
terms are zero, i.e.%
\[
t_{kl}=\cos\frac{k\pi}{a}x_{in}\;\cos\frac{l\pi}{b}y_{in}-\cos\frac{k\pi}%
{a}x_{out}\;\cos\frac{l\pi}{b}y_{out}=0
\]
which happens, for $m_{j}\in\mathbb{Z}$ for $j=1,2,3,4$, if
\begin{align*}
&  \left\{  \left\{  k\frac{x_{in}}{a}=\frac{1}{2}+m_{1}\right\}  \cap\left\{
k\frac{x_{out}}{a}=\frac{1}{2}+m_{2}\right\}  \right\} \\
&  \cup\left\{  \left\{  k\frac{x_{in}}{a}=\frac{1}{2}+m_{1}\right\}
\cap\left\{  l\frac{y_{out}}{b}=\frac{1}{2}+m_{3}\right\}  \right\} \\
&  \cup\left\{  \left\{  l\frac{y_{in}}{b}=\frac{1}{2}+m_{4}\right\}
\cap\left\{  k\frac{x_{out}}{a}=\frac{1}{2}+m_{2}\right\}  \right\} \\
&  \cup\left\{  \left\{  l\frac{y_{in}}{b}=\frac{1}{2}+m_{4}\right\}
\cap\left\{  l\frac{y_{out}}{b}=\frac{1}{2}+m_{3}\right\}  \right\}
\end{align*}
Now, $k\frac{2x_{in}}{a}=1+m_{1}$ requires that $\frac{2x_{in}}{a}=\frac
{n_{1}}{n_{2}}\in\mathbb{Q}$, where both $n_{1}$ and $n_{2}$ are non-negative
integers, because $0\leq x_{in}\leq a$. If $\frac{x_{in}}{a},\frac{x_{out}}%
{a},\frac{y_{in}}{b}$ and $\frac{y_{out}}{b}$ are not rational numbers (that
form the overwhelming majority of cases), then there does not exist any pair
$\left(  k,l\right)  $ of positive integers for which any term $t_{kl}$ is
zero. Hence, there exist a positive number%
\[
\beta=\frac{4}{ab}\min_{\left(  k,l\right)  }\left(  \cos\frac{k\pi}{a}%
x_{in}\;\cos\frac{l\pi}{b}y_{in}-\cos\frac{k\pi}{a}x_{out}\;\cos\frac{l\pi}%
{b}y_{out}\right)  ^{2}%
\]
so that double sum $S$ in (\ref{effective_resistance_RGG_continuun}) is lower
bounded by%
\begin{equation}
S>\beta\sum_{k=1}^{\infty}\sum_{l=1}^{\infty}\frac{1}{\left(  \frac{k\pi}%
{a}\right)  ^{2}+\left(  \frac{l\pi}{b}\right)  ^{2}}=\beta H\left(
a,b\right)  \label{lowerBound_S}%
\end{equation}

The complex number $\frac{k\pi}{a}+i\frac{l\pi}{b}$ has the square norm
$\left\vert \frac{k\pi}{a}+i\frac{l\pi}{b}\right\vert ^{2}=\left(  \frac{k\pi
}{a}\right)  ^{2}+\left(  \frac{l\pi}{b}\right)  ^{2}$. We now review some
basics from elliptic function theory and follow Rademacher \cite[Chapter
8]{Rademacher}. Let $\omega_{1}$ and $\omega_{2}$ be two complex numbers,
different from zero, such that their quotient $\tau=\frac{\omega_{2}}%
{\omega_{1}}$ is not real. The totality $\Omega$ of all complex numbers
$m_{1}\omega_{1}+m_{2}\omega_{2}$, with integers $m_{1}$ and $m_{2}$, forms a
point-lattice in the complex plane. The Eisenstein series%
\begin{equation}
G_{p}\left(  \omega_{1},\omega_{2}\right)  =\sum_{\left(  m_{1},m_{2}\right)
\neq\left(  0,0\right)  }\frac{1}{\left(  m_{1}\omega_{1}+m_{2}\omega
_{2}\right)  ^{p}} \label{Eisenstein_series_r}%
\end{equation}
sums over all pairs $\left(  m_{1},m_{2}\right)  $ of integers, except for the
pair $m_{1}=m_{2}=0$. Rademacher \cite[Lemma, p. 117]{Rademacher} proofs
Eisenstein's convergence lemma that, for $p>2$, the series%
\[
\sum_{\left(  m_{1},m_{2}\right)  \neq\left(  0,0\right)  }\frac{1}{\left\vert
m_{1}\omega_{1}+m_{2}\omega_{2}\right\vert ^{p}}%
\]
is convergent. We return to our series $H\left(  a,b\right)  $, which
satisfies%
\[
\left\vert \sum_{k=1}^{\infty}\sum_{l=1}^{\infty}\frac{1}{\left(  \frac{k\pi
}{a}+i\frac{l\pi}{b}\right)  ^{2}}\right\vert \leq\sum_{k=1}^{\infty}%
\sum_{l=1}^{\infty}\frac{1}{\left(  \frac{k\pi}{a}\right)  ^{2}+\left(
\frac{l\pi}{b}\right)  ^{2}}%
\]
The left-hand side is a special case of an Eisenstein series
(\ref{Eisenstein_series_r}) with $\omega_{1}=\frac{\pi}{a}$ and $\omega
_{2}=i\frac{\pi}{b}$, which does not converge absolutely, because $p=2$,
whereas convergence demands a strictly larger value of $p$ than 2. A
convergent series with exponent $p=2$ has led
Weierstrass\footnote{Weierstrass's $\wp\left(  u\right)  $ function can be
regarded as a generalization to a two-dimensions complex lattice of Gauss's
digamma function $\psi\left(  u\right)  =\frac{d}{du}\log\Gamma\left(
u\right)  $ in \cite[Chapter 12]{PVM_Mittag-Leffler_Gamma}, because%
\[
\psi\left(  u+1\right)  =-\gamma-\sum_{n=1}^{\infty}\left(  \frac{1}%
{u+n}-\frac{1}{n}\right)
\]
} to propose his famous, basic elliptic function
\begin{equation}
\wp\left(  u\right)  =\frac{1}{u^{2}}+\sum_{\left(  m_{1},m_{2}\right)
\neq\left(  0,0\right)  }\left(  \frac{1}{\left(  u-\left(  m_{1}\omega
_{1}+m_{2}\omega_{2}\right)  \right)  ^{2}}-\frac{1}{\left(  m_{1}\omega
_{1}+m_{2}\omega_{2}\right)  ^{2}}\right)  \label{Pe_Weierstrass}%
\end{equation}

In summary, we arrive at the conclusion that the double sum $H\left(
a,b\right)  \rightarrow\infty$ and that the effective\ resistance $\omega
_{sd}$ in (\ref{effective_resistance_RGG_continuun}) does not exist. An
alternative proof of the non-existence of the effective\ resistance
$\omega_{sd}$, based on our $q$-series, is given in Appendix
\ref{sec_non_existence_effective_res_q_analysis}.

Due to the fact that $p=2$ in double sum $H\left(  a,b\right)  $ in
(\ref{lowerBound_S}), whereas Eisenstein's convergence lemma requires
$p=2+\varepsilon$, with $\varepsilon>0$ but arbitrarily small, numerical
evaluations of the double Fourier series $\omega_{sd}$ in
(\ref{effective_resistance_RGG_continuun}) diverge extremely slowly. For
example, for a source node $s=\left(  1/8,3/8\right)  $ and a drain node
$d=\left(  3/4,3/4\right)  $, the effective resistance $\omega_{sd}\left(
K\right)  $ increase with the number terms $K$ as $\omega_{sd}\left(
1000\right)  =2.79642$, $\omega_{sd}\left(  2000\right)  =3.0169$,
$\omega_{sd}\left(  3000\right)  =3.14591\approx\pi$, $\omega_{sd}\left(
4000\right)  =3.23745$, $\omega_{sd}\left(  5000\right)  =3.30846$ and
$\omega_{sd}\left(  10^{4}\right)  =3.52907$.

\subsection{Alternative demonstration of non-existence of the
effective\ resistance $\omega_{sd}$ in
(\ref{effective_resistance_RGG_continuun})}

\label{sec_non_existence_effective_res_q_analysis}The \textquotedblleft
effective\textquotedblright\ resistance (\ref{effective_resistance_pq}) in 2D
between node $s=\left(  x_{in},y_{in}\right)  $ and node $d=\left(
x_{out},y_{out}\right)  $ becomes with the Green's function $V\left(
x,y\right)  =K\left(  x,y;x_{in},y_{in}\right)  -K\left(  x,y;x_{out}%
,y_{out}\right)  $ in (\ref{Voltage_2D_Green})%
\begin{equation}
\omega_{sd}=K\left(  x_{in},y_{in};x_{in},y_{in}\right)  -K\left(
x_{in},y_{in};x_{out},y_{out}\right)  -K\left(  x_{out},y_{out};x_{in}%
,y_{in}\right)  +K\left(  x_{out},y_{out};x_{out},y_{out}\right)
\label{effective_resistance_Green}%
\end{equation}
Invoking the $q$-extension (\ref{quasi_Green_in_Tq_elliptic})%
\[
K\left(  x,y,\xi,\eta\right)  =\frac{1}{a}g_{0}\left(  y,\eta\right)
+T_{q}\left(  x,y;\xi,\eta\right)  +T_{q}\left(  x,y;\xi,-\eta\right)
+T_{q}\left(  -x,y;\xi,\eta\right)  +T_{q}\left(  -x,y;\xi,-\eta\right)
\]
then yields%
\[
K\left(  \xi,\eta,\xi,\eta\right)  =T_{q}\left(  \xi,\eta;\xi,\eta\right)
+T_{q}\left(  -\xi,\eta;\xi,\eta\right)  +T_{q}\left(  \xi,\eta;\xi
,-\eta\right)  +T_{q}\left(  -\xi,\eta;\xi,-\eta\right)
\]

The definition (\ref{def_Tq_function}) of $T_{q}\left(  x,y;\xi,\eta\right)  $%
\begin{align*}
T_{q}\left(  x,y;\xi,\eta\right)   &  =-\frac{1}{4\pi}\log\left(
1-2\cos\left(  \frac{\pi}{a}\left(  \xi+x\right)  \right)  e^{-\frac{\pi}%
{a}\left(  y-\eta\right)  }+e^{-2\frac{\pi}{a}\left(  y-\eta\right)  }\right)
\\
&  \hspace{0.5cm}-\frac{1}{4\pi}\log\prod_{k=1}^{\infty}\left(  1-2\cos\left(
\frac{\pi}{a}\left(  \xi+x\right)  \right)  e^{\frac{\pi}{a}\left(
y-\eta\right)  }q^{k}+e^{2\frac{\pi}{a}\left(  y-\eta\right)  }q^{2k}\right)
\\
&  \hspace{1cm}\times\left(  1-2\cos\left(  \frac{\pi}{a}\left(  \xi+x\right)
\right)  e^{-\frac{\pi}{a}\left(  y-\eta\right)  }q^{k}+e^{-2\frac{\pi}%
{a}\left(  y-\eta\right)  }q^{2k}\right)
\end{align*}
simplifies for $y=\eta$ to
\[
T_{q}\left(  x,\eta;\xi,\eta\right)  =-\frac{1}{2\pi}\log\left(  2\sin\left(
\frac{\pi}{2a}\left(  \xi+x\right)  \right)  \right)  -\frac{1}{2\pi}\log
\prod_{k=1}^{\infty}\left(  1-2\cos\left(  \frac{\pi}{a}\left(  \xi+x\right)
\right)  q^{k}+q^{2k}\right)
\]
Further, we find that%
\[
T_{q}\left(  \xi,\eta;\xi,\eta\right)  =-\frac{1}{2\pi}\log\left(  2\left\vert
\sin\left(  \frac{\pi}{a}\xi\right)  \right\vert \right)  -\frac{1}{2\pi}%
\log\prod_{k=1}^{\infty}\left(  1+q^{2k}-2q^{k}\cos\left(  \frac{2\pi}{a}%
\xi\right)  \right)
\]
and that%
\[
\lim_{x\rightarrow-\xi}T_{q}\left(  x,\eta;\xi,\eta\right)  =-\frac{1}{2\pi
}\lim_{x\rightarrow-\xi}\log\left(  2\sin\left(  \frac{\pi}{2a}\left(
\xi+x\right)  \right)  \right)  -\frac{1}{\pi}\log\prod_{k=1}^{\infty}\left(
1-q^{k}\right)
\]
Similarly, the definition (\ref{def_Tq_function}) of $T_{q}\left(
x,y;\xi,\eta\right)  $ indicates that%
\begin{align*}
T_{q}\left(  x,\eta;\xi,-\eta\right)   &  =-\frac{1}{4\pi}\log\left(
1-2\cos\left(  \frac{\pi}{a}\left(  \xi+x\right)  \right)  e^{-\frac{2\pi}%
{a}\eta}+e^{-\frac{4\pi}{a}\eta}\right) \\
&  \hspace{0.5cm}-\frac{1}{4\pi}\log\prod_{k=1}^{\infty}\left(  1-2\cos\left(
\frac{\pi}{a}\left(  \xi+x\right)  \right)  e^{\frac{2\pi}{a}\eta}%
q^{k}+e^{\frac{4\pi}{a}\eta}q^{2k}\right) \\
&  \hspace{0.5cm}\times\left(  1-2\cos\left(  \frac{\pi}{a}\left(
\xi+x\right)  \right)  e^{-\frac{2\pi}{a}\eta}q^{k}+e^{-\frac{4\pi}{a}\eta
}q^{2k}\right)
\end{align*}
from which%
\begin{align*}
T_{q}\left(  \xi,\eta;\xi,-\eta\right)   &  =-\frac{1}{4\pi}\log\left(
1-2\cos\left(  \frac{2\pi}{a}\xi\right)  e^{-\frac{2\pi}{a}\eta}%
+e^{-\frac{4\pi}{a}\eta}\right) \\
&  \hspace{0.5cm}-\frac{1}{4\pi}\log\prod_{k=1}^{\infty}\left(  1-2\cos\left(
\frac{2\pi}{a}\xi\right)  e^{\frac{2\pi}{a}\eta}q^{k}+e^{\frac{4\pi}{a}\eta
}q^{2k}\right)  \left(  1-2\cos\left(  \frac{2\pi}{a}\xi\right)
e^{-\frac{2\pi}{a}\eta}q^{k}+e^{-\frac{4\pi}{a}\eta}q^{2k}\right)
\end{align*}
and%
\begin{align*}
T_{q}\left(  -\xi,\eta;\xi,-\eta\right)   &  =-\frac{1}{2\pi}\log\left(
1-e^{-\frac{2\pi}{a}\eta}\right) \\
&  \hspace{0.5cm}-\frac{1}{4\pi}\log\prod_{k=1}^{\infty}\left(  1-2e^{\frac
{2\pi}{a}\eta}q^{k}+e^{\frac{4\pi}{a}\eta}q^{2k}\right)  \left(
1-2e^{-\frac{2\pi}{a}\eta}q^{k}+e^{-\frac{4\pi}{a}\eta}q^{2k}\right)
\end{align*}
Hence, we observe that only in the term $\lim_{x\rightarrow-\xi}T_{q}\left(
x,\eta;\xi,\eta\right)  $, a logarithmic divergence occurs of the type
$-\frac{1}{2\pi}\lim_{x\rightarrow-\xi}\log\left(  2\sin\left(  \frac{\pi}%
{2a}\left(  \xi+x\right)  \right)  \right)  $. Such a term occurs both in
$K\left(  x_{in},y_{in};x_{in},y_{in}\right)  $ and $K\left(  x_{out}%
,y_{out};x_{out},y_{out}\right)  $ with the same sign. Since the remaining
contributions $K\left(  x_{in},y_{in};x_{out},y_{out}\right)  +K\left(
x_{out},y_{out};x_{in},y_{in}\right)  $ do not contain a divergent term, we
again find that the effective resistance $\omega_{sd}$ in
(\ref{effective_resistance_Green}) does not exist.

\subsection{The resistance $R_{sd}$ in (\ref{resistance_3D}) in 3D does not
exist either}

Since $\coth\left(  c\sqrt{\left(  \frac{k\pi}{a}\right)  ^{2}+\left(
\frac{l\pi}{b}\right)  ^{2}}\right)  $ rapidly tends to 1 for sufficiently
large $k$ and $l$ in the 3D resistance (\ref{resistance_3D}) the corresponding
double sum tends to%
\[
S_{3D}\simeq\frac{4\varrho}{ab}\sum_{k=K_{k}}^{\infty}\sum_{l=K_{k}}^{\infty
}\frac{\left(  \cos\frac{k\pi}{a}x_{in}\cos\frac{l\pi}{b}y_{in}-\cos\frac
{k\pi}{a}x_{out}\cos\frac{l\pi}{b}y_{out}\right)  ^{2}}{\sqrt{\left(
\frac{k\pi}{a}\right)  ^{2}+\left(  \frac{l\pi}{b}\right)  ^{2}}}%
\]
By a similar argument as Section \ref{sec_effective_res_Eisenstein} and
Eisenstein's convergence lemma, we conclude that also the resistance $R_{sd}$
in (\ref{resistance_3D}) in 3D does not exist!

\subsection{Consequence for the effective resistance $\omega_{sd}$ in RGG}

Although our mathematical derivations are correct, the computation of a simple
physical concept as the resistance reveals an error. The error is a
\textquotedblleft modeling fault\textquotedblright\ of reality: the assumption
of current in- and ejection in a point relies on delta functions that simplify
the mathematical equations. However, forcing a current to pass through a
single point with zero dimensions causes the resistance between injection
point $s$ and ejection point $d$ to grow infinitely large, for any
configuration of the pair $\left(  s,d\right)  $.

We replace the point injection by a small finite area. The simplest form in
our Euclidian coordinate frame is a rectangular in- and ejection area over
which the current is homogeneously spread with constant the current density
$j=\frac{I}{h_{xin}h_{yin}}$, which results in the function $g\left(
x,y\right)  $ in Section \ref{sec_current_point_injection_3D},
\[
g\left(  x,y\right)  =\left\{
\begin{array}
[c]{c}%
\frac{\varrho I}{h_{xin}h_{yin}}\text{ for }x_{in}-\frac{h_{xin}}{2}\leq x\leq
x_{in}+\frac{h_{xin}}{2}\text{ and }y_{in}-\frac{h_{yin}}{2}\leq y\leq
y_{in}+\frac{h_{yin}}{2}\\
\frac{\varrho I}{h_{xout}h_{yout}}\text{ for }x_{out}-\frac{h_{xout}}{2}\leq
x\leq x_{out}+\frac{h_{xout}}{2}\text{ and }y_{out}-\frac{h_{yout}}{2}\leq
x\leq y_{out}+\frac{h_{yout}}{2}%
\end{array}
\right.
\]
The correspondingly modified Fourier coefficients follow from
(\ref{Fourier_coeff_psi}) and require the computations of the double integral%
\begin{align*}
I_{2}  &  =\int_{0}^{a}dx\;\int_{0}^{b}dy\;g\left(  x,y\right)  \cos\frac{k\pi
x}{a}\cos\frac{l\pi y}{b}=\frac{\varrho I}{h_{xin}h_{yin}}\int_{x_{in}%
-\frac{h_{xin}}{2}}^{x_{in}+\frac{h_{xin}}{2}}\cos\frac{k\pi x}{a}%
dx\;\int_{y_{in}-\frac{h_{yin}}{2}}^{y_{in}+\frac{h_{yin}}{2}}\cos\frac{l\pi
y}{b}dy\\
&  =\frac{\varrho I}{h_{xin}h_{yin}}\frac{\sin\frac{k\pi}{a}\left(
x_{in}+\frac{h_{xin}}{2}\right)  -\sin\frac{k\pi}{a}\left(  x_{in}%
-\frac{h_{xin}}{2}\right)  }{\frac{k\pi}{a}}\frac{\sin\frac{l\pi}{b}\left(
y_{in}+\frac{h_{xin}}{2}\right)  -\sin\frac{l\pi}{b}\left(  x_{in}%
-\frac{h_{xin}}{2}\right)  }{\frac{l\pi}{b}}\\
&  =\varrho I\cos\left(  \frac{k\pi}{a}x_{in}\right)  \cos\left(  \frac{l\pi
}{b}y_{in}\right)  \times\frac{\sin\left(  \frac{k\pi}{a}\frac{h_{xin}}%
{2}\right)  }{\frac{k\pi}{a}\frac{h_{xin}}{2}}\frac{\sin\left(  \frac{l\pi}%
{b}\frac{h_{yin}}{2}\right)  }{\frac{l\pi}{b}\frac{h_{yin}}{2}}%
\end{align*}
Comparing the point injection in Section \ref{sec_current_point_injection_3D}
with a non-zero rectangular area injection reveals that only a product of
sinc$\left(  X\right)  =\frac{\sin X}{X}$ is added. Provided that the
rectangular injection area completely lies within the $xy$-face of the beam,
i.e. $\left\{  0\leq x_{in}-\frac{h_{xin}}{2}\right\}  \cap\left\{
x_{in}+\frac{h_{xin}}{2}\leq a\right\}  $ and $\left\{  0\leq y_{in}%
-\frac{h_{yin}}{2}\right\}  \cap\left\{  y_{in}+\frac{h_{yin}}{2}\leq
b\right\}  $ and similarly for the rectangular ejection area, the modified 3D
voltage is, corresponding to (\ref{Voltage_beam}), is%
\begin{align}
V\left(  x,y,z\right)  +V_{\text{ref}}  &  =\frac{2\varrho I}{ab}\sum
_{k=1}^{\infty}\frac{\left(  \cos\frac{k\pi}{a}x_{in}\frac{\sin\left(
\frac{k\pi}{a}\frac{h_{xin}}{2}\right)  }{\frac{k\pi}{a}\frac{h_{xin}}{2}%
}\;-\cos\frac{k\pi}{a}x_{out}\;\frac{\sin\left(  \frac{k\pi}{a}\frac{h_{xout}%
}{2}\right)  }{\frac{k\pi}{a}\frac{h_{xout}}{2}}\right)  }{\frac{k\pi}{a}%
\sinh\left(  c\frac{k\pi}{a}\right)  }\cos\frac{k\pi}{a}x\cosh\left(
\frac{k\pi}{a}\left(  z-c\right)  \right) \nonumber\\
&  \hspace{0.5cm}+\frac{2\varrho I}{ab}\sum_{l=1}^{\infty}\frac{\left(
\cos\frac{l\pi}{b}y_{in}\frac{\sin\left(  \frac{l\pi}{b}\frac{h_{yin}}%
{2}\right)  }{\frac{l\pi}{b}\frac{h_{yin}}{2}}-\cos\frac{l\pi}{b}y_{out}%
\frac{\sin\left(  \frac{l\pi}{b}\frac{h_{yout}}{2}\right)  }{\frac{l\pi}%
{b}\frac{h_{yout}}{2}}\right)  }{\frac{l\pi}{b}\sinh\left(  c\frac{l\pi}%
{b}\right)  }\cos\frac{l\pi}{b}y\cosh\left(  \frac{l\pi}{b}\left(  z-c\right)
\right) \nonumber\\
&  \hspace{0.5cm}+\frac{4\varrho I}{ab}\sum_{k=1}^{\infty}\sum_{l=1}^{\infty
}\frac{\left(
\begin{array}
[c]{c}%
\cos\frac{k\pi}{a}x_{in}\;\cos\frac{l\pi}{b}y_{in}\frac{\sin\left(  \frac
{k\pi}{a}\frac{h_{xin}}{2}\right)  }{\frac{k\pi}{a}\frac{h_{xin}}{2}}%
\frac{\sin\left(  \frac{l\pi}{b}\frac{h_{yin}}{2}\right)  }{\frac{l\pi}%
{b}\frac{h_{yin}}{2}}\\
-\cos\frac{k\pi}{a}x_{out}\;\cos\frac{l\pi}{b}y_{out}\frac{\sin\left(
\frac{k\pi}{a}\frac{h_{xout}}{2}\right)  }{\frac{k\pi}{a}\frac{h_{xout}}{2}%
}\frac{\sin\left(  \frac{l\pi}{b}\frac{h_{yout}}{2}\right)  }{\frac{l\pi}%
{b}\frac{h_{yout}}{2}}%
\end{array}
\right)  }{\sinh\left(  c\sqrt{\left(  \frac{k\pi}{a}\right)  ^{2}+\left(
\frac{l\pi}{b}\right)  ^{2}}\right)  \sqrt{\left(  \frac{k\pi}{a}\right)
^{2}+\left(  \frac{l\pi}{b}\right)  ^{2}}}\nonumber\\
&  \hspace{0.5cm}\hspace{1cm}\times\cos\frac{k\pi}{a}x\cos\frac{l\pi}{b}%
y\cosh\left(  \sqrt{\left(  \frac{k\pi}{a}\right)  ^{2}+\left(  \frac{l\pi}%
{b}\right)  ^{2}}\left(  z-c\right)  \right)
\label{Voltage_beam_rectangular_contacts}%
\end{align}
Shrinking the rectangle's dimensions $h_{xin}$ and $h_{yin}$ with center point
$s=\left(  x_{in},y_{in}\right)  $ to zero, re-establishes the point injection
at $\left(  x_{in},y_{in}\right)  $ and the voltage
(\ref{Voltage_beam_rectangular_contacts}) reduces to (\ref{Voltage_beam}). The
modified resistance $R_{sd}=\frac{V\left(  x_{in},y_{in},0\right)  -V\left(
x_{out},y_{out},0\right)  }{I}$, with rectangular current in- and ejection
area centered around the source $s$ and drain $d$, respectively, and lying
entirely within the beam's faces, is%
\begin{equation}
R_{sd}^{\prime}=\frac{4\varrho}{ab}\sum_{k=0}^{\infty}\sum_{l=0}^{\infty
}\gamma_{kl}\frac{\left(
\begin{array}
[c]{c}%
\cos\frac{k\pi}{a}x_{in}\;\cos\frac{l\pi}{b}y_{in}\frac{\sin\left(  \frac
{k\pi}{a}\frac{h_{xin}}{2}\right)  }{\frac{k\pi}{a}\frac{h_{xin}}{2}}%
\frac{\sin\left(  \frac{l\pi}{b}\frac{h_{yin}}{2}\right)  }{\frac{l\pi}%
{b}\frac{h_{yin}}{2}}\\
-\cos\frac{k\pi}{a}x_{out}\;\cos\frac{l\pi}{b}y_{out}\frac{\sin\left(
\frac{k\pi}{a}\frac{h_{xout}}{2}\right)  }{\frac{k\pi}{a}\frac{h_{xout}}{2}%
}\frac{\sin\left(  \frac{l\pi}{b}\frac{h_{yout}}{2}\right)  }{\frac{l\pi}%
{b}\frac{h_{yout}}{2}}%
\end{array}
\right)  ^{2}}{\sqrt{\left(  \frac{k\pi}{a}\right)  ^{2}+\left(  \frac{l\pi
}{b}\right)  ^{2}}}\coth\left(  c\sqrt{\left(  \frac{k\pi}{a}\right)
^{2}+\left(  \frac{l\pi}{b}\right)  ^{2}}\right)
\label{resistance_3D_rectangular_area}%
\end{equation}
For sufficiently large $k$ and $l$, the double sum in
(\ref{resistance_3D_rectangular_area}) tends to%
\[
S_{3D}^{\prime}\simeq\frac{4\varrho}{ab}\sum_{k=K_{k}}^{\infty}\sum_{l=K_{k}%
}^{\infty}\frac{\left(  \cos\frac{k\pi}{a}x_{in}\cos\frac{l\pi}{b}y_{in}%
\frac{\sin\left(  \frac{k\pi}{a}\frac{h_{xin}}{2}\right)  }{\frac{h_{xin}}{2}%
}\frac{\sin\left(  \frac{l\pi}{b}\frac{h_{yin}}{2}\right)  }{\frac{h_{yin}}%
{2}}-\cos\frac{k\pi}{a}x_{out}\cos\frac{l\pi}{b}y_{out}\frac{\sin\left(
\frac{k\pi}{a}\frac{h_{xout}}{2}\right)  }{\frac{h_{xout}}{2}}\frac
{\sin\left(  \frac{l\pi}{b}\frac{h_{yout}}{2}\right)  }{\frac{h_{yout}}{2}%
}\right)  ^{2}}{\frac{k\pi}{a}\frac{l\pi}{b}\sqrt{\left(  \frac{k\pi}%
{a}\right)  ^{2}+\left(  \frac{l\pi}{b}\right)  ^{2}}}%
\]
which always converges for all involved \emph{positive} parameters, because
both sums are convergent series of the type $\sum_{n=K_{n}}^{\infty}\frac
{1}{n^{1+\varepsilon}}$ with $\varepsilon>0$. Taking the limit $c\rightarrow0$
combined with the renormalization $\frac{\varrho I}{c}=1$ transforms the 3D
voltage (\ref{Voltage_beam_rectangular_contacts}) with rectangular current in-
and ejection as well as the corresponding resistance
(\ref{resistance_3D_rectangular_area}) to the 2D analogon (whose Fourier
equations are also verified to convergence always for positive rectangle sizes
$h_{xin}>0$, $h_{yin}>0$, $h_{xout}>0$ and $h_{yout}>0$.)

In spite of the mathematical elegant modeling of point in- and ejection,
caution is always required to deduce conclusions beyond the (sometimes
implicitly made) modeling assumptions.

\section{The function $w\left(  z,A\right)  $ in (\ref{w(z,A)_series_def}) and
$q$-analysis}

\label{sec_function_w(z,A)}We deduce properties of the function $w\left(
z,A\right)  $, defined by the series
\[
w\left(  z,A\right)  =\sum_{n=1}^{\infty}\frac{z^{n}}{n\sinh\left(  An\right)
}%
\]
in (\ref{w(z,A)_series_def}). Provided that $A\neq0$, the series $\sum
_{n=1}^{\infty}\frac{z^{n}}{n\sinh\left(  An\right)  }$ in
(\ref{w(z,A)_series_def}) converges for $\left\vert ze^{-A}\right\vert <1$.
Indeed, the absolute value is%
\[
\left\vert w\left(  z,A\right)  \right\vert \leq2\sum_{n=1}^{\infty}%
\frac{\left\vert z\right\vert ^{n}}{n\left\vert e^{An}-e^{-An}\right\vert
}=2\sum_{n=1}^{\infty}\frac{\left\vert ze^{-A}\right\vert ^{n}}{n\left\vert
1-e^{-2An}\right\vert }%
\]
Since $\left\vert 1-e^{-2A}\right\vert <\left\vert 1-e^{-2An}\right\vert
\leq1$ for all integer $n>0$, we conclude, for $\left\vert ze^{-A}\right\vert
<1$, that%
\[
\left\vert w\left(  z,A\right)  \right\vert \leq\frac{-2\log\left(
1-\left\vert ze^{-A}\right\vert \right)  }{\left\vert 1-e^{-2A}\right\vert }%
\]

We rewrite (\ref{w(z,A)_series_def}) as%
\[
w\left(  z,A\right)  =2\sum_{n=1}^{\infty}\frac{\left(  ze^{-A}\right)  ^{n}%
}{n\left(  1-e^{-2An}\right)  }=2\sum_{n=1}^{\infty}\sum_{k=0}^{\infty}\left(
e^{-2An}\right)  ^{k}\frac{\left(  ze^{-A}\right)  ^{n}}{n}%
\]
We can reverse the sums due to absolute convergence,%
\[
w\left(  z,A\right)  =2\sum_{k=0}^{\infty}\sum_{n=1}^{\infty}\frac{\left(
ze^{-A\left(  1+2k\right)  }\right)  ^{n}}{n}=-2\sum_{k=0}^{\infty}\log\left(
1-ze^{-A\left(  1+2k\right)  }\right)  =-2\log\prod_{k=0}^{\infty}\left(
1-ze^{-A\left(  1+2k\right)  }\right)
\]
Hence, for $\left\vert ze^{-A}\right\vert <1$, we find that%
\begin{align}
w\left(  z,A\right)   &  =\sum_{n=1}^{\infty}\frac{z^{n}}{n\sinh\left(
An\right)  }=-2\sum_{k=0}^{\infty}\log\left(  1-ze^{-A\left(  1+2k\right)
}\right) \label{w(z,A)_series_logs}\\
&  =-2\log\left(  1-ze^{-A}\right)  -2\log\left(  1-ze^{-3A}\right)
-2\log\left(  1-ze^{-5A}\right)  +\cdots\nonumber
\end{align}

The definition (\ref{w(z,A)_series_def}) of $w\left(  z,A\right)  =2\sum
_{n=1}^{\infty}\frac{z^{n}}{n\left(  e^{An}-e^{-An}\right)  }$, for
$\left\vert z\right\vert <1$, indicates that%
\[
w\left(  ze^{A},A\right)  -w\left(  ze^{-A},A\right)  =2\sum_{n=1}^{\infty
}\frac{z^{n}\left(  e^{An}-e^{-An}\right)  }{n\left(  e^{An}-e^{-An}\right)
}=2\sum_{n=1}^{\infty}\frac{z^{n}}{n}%
\]
from which we obtain%
\[
w\left(  ze^{A},A\right)  -w\left(  ze^{-A},A\right)  =-2\log\left(
1-z\right)
\]
Replace $x=ze^{A}$, then $z=e^{-A}x$ and $w\left(  x,A\right)  =-2\log\left(
1-e^{-A}x\right)  +w\left(  xe^{-2A},A\right)  $. Writing again $z$ for $x$,
results in the recursion%
\begin{equation}
w\left(  z,A\right)  =-2\log\left(  1-e^{-A}z\right)  +w\left(  ze^{-2A}%
,A\right)  \label{w(z,A)_recursion}%
\end{equation}
After $p$ iterations of the recursion (\ref{w(z,A)_recursion}), we find that%
\begin{equation}
w\left(  z,A\right)  =-2\sum_{k=0}^{p}\log\left(  1-ze^{-A\left(  1+2k\right)
}\right)  +w\left(  ze^{-2pA},A\right)  \label{w(z,A)_p_times_iterated}%
\end{equation}
Since $e^{-2pA}$ rapidly tends to zero if $p$ increases and $w\left(
0,A\right)  =0$, the $p$-th iteration (\ref{w(z,A)_p_times_iterated}) leads
again to (\ref{w(z,A)_series_logs}) in the limit $p\rightarrow\infty$. In
other words, the recursion (\ref{w(z,A)_recursion}) allows us to write the sum
in $w\left(  z,A\right)  $ up to any desired accuracy.

Interestingly, the product in $w\left(  z,A\right)  =$ $-2\log\left(
\prod_{k=0}^{\infty}\left(  1-ze^{-A\left(  1+2k\right)  }\right)  \right)  $
is written, with $q=e^{-2A}$, as
\[
\prod_{k=0}^{\infty}\left(  1-ze^{-A\left(  1+2k\right)  }\right)
=\prod_{k=0}^{\infty}\left(  1-ze^{-A}q^{k}\right)  =\left(  ze^{-A};q\right)
_{\infty}%
\]
where%
\[
\left(  a;q\right)  _{n}=\prod_{k=0}^{n-1}\left(  1-aq^{k}\right)
\]
is known as the $q$-Pochhammer symbol, with $\left(  a;q\right)  _{0}=1$.
Then, with $q=e^{-2A}$, we have
\begin{equation}
w\left(  z,A\right)  =-2\log\prod_{k=0}^{\infty}\left(  1-ze^{-A}q^{k}\right)
=-2\log\left(  z\sqrt{q};q\right)  _{\infty} \label{w_log_product}%
\end{equation}
Invoking $\prod_{m=0}^{\infty}(1+q^{m}\,z)=\sum_{m=0}^{\infty}\frac
{q^{m(m-1)/2}}{\prod_{j=1}^{m}(1-q^{j})}\,z^{m}$ in (\ref{prop_Gausspol_infty}%
) shows, with $q=e^{-2A}$, that%
\begin{equation}
\prod_{k=0}^{\infty}\left(  1-ze^{-A}q^{k}\right)  =\sum_{m=0}^{\infty}%
\frac{q^{\frac{m^{2}}{2}}}{\prod_{j=1}^{m}(1-q^{j})}\,\left(  -z\right)  ^{m}
\label{binomial-q_series_product_in_w(z,A)}%
\end{equation}
Hence, the product in (\ref{w_log_product}) is%
\begin{align*}
\prod_{k=0}^{\infty}\left(  1-ze^{-A}q^{k}\right)   &  =\sum_{m=0}^{\infty
}\frac{e^{-Am^{2}}}{\prod_{j=1}^{m}(1-e^{-2Aj})}\,\left(  -z\right)  ^{m}\\
&  =1-\frac{e^{-A}}{1-e^{-2A}}\,z+\frac{e^{-4A}}{(1-e^{-2A})(1-e^{-4A}%
)}\,z^{2}+\sum_{m=3}^{\infty}\frac{e^{-Am^{2}}}{\prod_{j=1}^{m}(1-e^{-2Aj}%
)}\,\left(  -z\right)  ^{m}%
\end{align*}
In our setting $q=e^{-2A}=e^{-2\pi\frac{b}{a}}$ is, for a square where $a=b$,
equal to $q=e^{-2\pi}=\frac{1}{535.491}\approx2.10^{-3}$ and $e^{-A}=e^{-\pi
}\approx4.10^{-2}$, which illustrates that the series
(\ref{binomial-q_series_product_in_w(z,A)}) converges very fast, implying that
only a few terms are sufficient for numerical computations. Thus, for
$q=e^{-2A}$, the product form (\ref{w_log_product}) of $w\left(  z,A\right)
$
\begin{align*}
w\left(  z,A\right)   &  =-2\log\prod_{k=0}^{\infty}\left(  1-ze^{-A}%
q^{k}\right) \\
&  \simeq-2\log\left(  1-\frac{e^{-A}}{1-e^{-2A}}\,z+\sum_{m=2}^{K}%
\frac{e^{-Am^{2}}}{\prod_{j=1}^{m}(1-e^{-2Aj})}\,\left(  -z\right)
^{m}\right)
\end{align*}
already yields a very good accuracy for a small integer $K\approx5$, provided
$A$ is not too small.

The function (\ref{w(z,A)_series_def})%
\[
w\left(  z,A\right)  =2\sum_{n=1}^{\infty}\frac{z^{n}}{n\left(  e^{An}%
-e^{-An}\right)  }=2\sum_{n=1}^{\infty}\frac{\left(  ze^{A}\right)  ^{n}%
}{n\left(  e^{2An}-1\right)  }%
\]
for small $A$ is computed with the Bernoulli expansion \cite[Appendix]%
{PVM_Mittag-Leffler_Gamma}, valid for $\left\vert t\right\vert <2\pi$,
\begin{equation}
\frac{1\,}{e^{t}-1}=\frac{1}{t}-\frac{1}{2}+\sum_{m=1}^{\infty}B_{2m}%
\,\frac{t^{2m-1}}{(2m)!} \label{gf_Bernoullinumbers}%
\end{equation}
where $B_{n}$ is the $n$-th Bernoulli number, as%
\begin{align*}
w\left(  z,A\right)   &  =2\sum_{n=1}^{\infty}\frac{\left(  ze^{A}\right)
^{n}}{n\left(  e^{2An}-1\right)  }\\
&  =2\lim_{N\rightarrow\infty}\sum_{n=1}^{N}\frac{\left(  ze^{A}\right)  ^{n}%
}{n}\left(  \frac{1}{2An}-\frac{1}{2}+\sum_{m=1}^{\infty}B_{2m}\,\frac{\left(
2An\right)  ^{2m-1}}{(2m)!}\right)
\end{align*}
Assuming that $A<\frac{\pi}{N}$, reversal of $n$- and $m$-sum shows that%
\begin{align*}
w\left(  z,A\right)   &  =\lim_{N\rightarrow\infty}\left(  \frac{1}{A}%
\sum_{n=1}^{N}\frac{\left(  ze^{A}\right)  ^{n}}{n^{2}}-\sum_{n=1}^{N}%
\frac{\left(  ze^{A}\right)  ^{n}}{n}+2\sum_{m=1}^{\infty}B_{2m}%
\,\frac{\left(  2A\right)  ^{2m-1}\sum_{n=1}^{N}\frac{\left(  ze^{A}\right)
^{n}}{n^{2\left(  1-m\right)  }}}{(2m)!}\right) \\
&  =\lim_{N\rightarrow\infty}\left(  \frac{1}{A}\sum_{n=1}^{N}\frac{\left(
ze^{A}\right)  ^{n}}{n^{2}}-\sum_{n=1}^{N}\frac{\left(  ze^{A}\right)  ^{n}%
}{n}+O\left(  A\right)  \right)
\end{align*}
Since the two remaining sums exist for $N\rightarrow\infty$ and%
\[
\sum_{n=1}^{\infty}\frac{x^{n}}{n^{2}}=\int_{0}^{x}\sum_{n=1}^{\infty}%
\frac{u^{n-1}}{n}du=-\int_{0}^{x}\log\left(  1-u\right)  \frac{du}{u}=\int
_{0}^{-\log\left(  1-x\right)  }\frac{t}{e^{t}-1}dt
\]
we arrive, for small $A\rightarrow0$ and $\left\vert ze^{A}\right\vert <1$, at%
\[
w\left(  z,A\right)  =\frac{1}{A}\int_{0}^{-\log\left(  1-z\right)  }\frac
{t}{e^{t}-1}dt+\log\left(  1-z\right)  +O\left(  A\right)
\]
where the first term clearly dominates.

\section{Gaussian polynomials}

\label{sec_Gaussianpoly}

\subsection{Definition}

\label{sec_def_Gaussianpoly}The Gaussian polynomial for $k\geq0$ and $l>0$ is
defined \cite[pp. 250]{Rademacher} as
\begin{align}
\left[
\begin{array}
[c]{c}%
k\\
l
\end{array}
\right]  (q)  &  =\frac{(1-q^{k})(1-q^{k-1})\cdots(1-q^{k-l+1})}%
{(1-q)(1-q^{2})\cdots(1-q^{l})}=\frac{\prod_{j=k-l+1}^{k}(1-q^{j})}%
{\prod_{j=1}^{l}(1-q^{j})}=\prod_{j=1}^{l}\frac{(1-q^{k-j+1})}{(1-q^{j}%
)}\label{def_Gausspolynoom}\\
&  =\frac{\prod_{j=1}^{k}(1-q^{j})}{\prod_{j=1}^{l}(1-q^{j})\;\prod
_{j=1}^{k-l}(1-q^{j})} \label{Gauss_binomialcoeff}%
\end{align}
with $\left[
\begin{array}
[c]{c}%
0\\
l
\end{array}
\right]  (q)=\delta_{0l}$, $\left[
\begin{array}
[c]{c}%
k\\
0
\end{array}
\right]  (q)=1$ and $\left[
\begin{array}
[c]{c}%
k\\
-l
\end{array}
\right]  (q)=0$. Obviously, it holds that $\left[
\begin{array}
[c]{c}%
k\\
l
\end{array}
\right]  (0)=1$. We list some Gaussian polynomials,
\begin{align*}
\left[
\begin{array}
[c]{c}%
3\\
2
\end{array}
\right]  (q)  &  =1+q+q^{2}\\
\left[
\begin{array}
[c]{c}%
4\\
2
\end{array}
\right]  (q)  &  =1+q+2\,q^{2}+q^{3}+q^{4}\\
\left[
\begin{array}
[c]{c}%
5\\
2
\end{array}
\right]  (q)  &  =1+q+2\,q^{2}+2\,q^{3}+2\,q^{4}+q^{5}+q^{6}%
\end{align*}
The Gaussian polynomials (\ref{def_Gausspolynoom}) are particularly
interesting due to their intimate relation to the binomial coefficients,
exemplified by the form of (\ref{Gauss_binomialcoeff}) and by
\begin{equation}
\left[
\begin{array}
[c]{c}%
k\\
l
\end{array}
\right]  (1)={\binom{k}{l}} \label{binomial_Gausspoly}%
\end{equation}
Indeed,%
\begin{align*}
\lim_{q\rightarrow1}\left[
\begin{array}
[c]{c}%
k\\
l
\end{array}
\right]  (q)  &  =\lim_{q\rightarrow1}\prod_{j=1}^{l}\frac{(1-q^{k-j+1}%
)}{(1-q^{j})}=\prod_{j=1}^{l}\lim_{q\rightarrow1}\frac{(1-q^{k-j+1})}%
{(1-q^{j})}=\prod_{j=1}^{l}\frac{(k-j+1)}{j}\\
&  =\frac{1}{l!}\frac{\prod_{j=1}^{k}(k-j+1)}{\prod_{j=l+1}^{k}(k-j+1)}%
=\frac{1}{l!}\frac{\prod_{r=1}^{k}r}{\prod_{r=k-l}^{1}r}=\frac{k!}{l!\left(
k-l\right)  !}={\binom{k}{l}}%
\end{align*}
The same argument shows that%
\[
\lim_{q\rightarrow1}\prod_{j=0}^{n-1}\frac{(1-q^{x+j})}{(1-q)}=\prod
_{j=0}^{n-1}(x+j)=\frac{\Gamma\left(  x+n\right)  }{\Gamma\left(  x\right)
}=\left(  x\right)  _{n}%
\]
suggesting, with the notation $\left(  a;q\right)  _{n}=\prod_{k=0}%
^{n-1}\left(  1-aq^{k}\right)  $, that $\frac{\left(  q^{x};q\right)  _{n}%
}{\left(  1-q\right)  ^{n}}$ is the $q$-extension of the Pochhammer symbol
$\left(  x\right)  _{n}$.

The other important limit, for $\left\vert q\right\vert <1$, is%
\begin{equation}
\lim_{k\rightarrow\infty}\left[
\begin{array}
[c]{c}%
k\\
l
\end{array}
\right]  (q)=\frac{1}{\prod_{j=1}^{l}(1-q^{j})}\hspace{1cm}(|q|<1)
\label{lim_Gausspoly}%
\end{equation}
which directly follows, for finite $l$, from $\left[
\begin{array}
[c]{c}%
k\\
l
\end{array}
\right]  (q)=\frac{\prod_{j=1}^{k}(1-q^{j})}{\prod_{j=1}^{l}(1-q^{j}%
)\;\prod_{j=1}^{k-l}(1-q^{j})}$.

We will demonstrate that many of the well-known relations involving binomial
coefficients readily follow as special cases (i.e. evaluations at $q=1$) of
the more general results based on Gaussian polynomials. Combinatorial results
based on Gaussian polynomials are discussed by Goulden and Jackson
\cite{Goulden_Jackson}.

\subsection{Generating function for the Gaussian Polynomials}

\label{sec_genfunc_gauss} The fundamental cornerstone in the theory of
Gaussian polynomials is
\begin{equation}
Q_{k}(z,x)=\prod_{m=0}^{k-1}(x+q^{m}\,z)=\sum_{m=0}^{k}\left[
\begin{array}
[c]{c}%
k\\
m
\end{array}
\right]  (q)\;q^{m(m-1)/2}\,z^{m}\,x^{k-m} \label{prop_Gausspol}%
\end{equation}
which bears a striking resemblance to Newton's binomium
\cite{Goulden_Jackson,Rademacher}; see also \cite[Chapter 2, Sec.
3.3]{Andrews_Partitions}. We define $Q_{k}(-x,x)=\delta_{0k}$ in
correspondence to the first factor for $m=0$ in the product. Relation
(\ref{prop_Gausspol}) is derived via induction from the recursion
$Q_{k}(z,x)=(x+q^{k-1}\,z)\,Q_{k-1}(z,x)$ for $k>0$ and $Q_{0}(x,z)=1$.

When $k$ tends to infinity, (\ref{prop_Gausspol}) leads for $|q|<1$ with
(\ref{lim_Gausspoly}) to
\begin{equation}
\prod_{m=0}^{\infty}(1+q^{m}\,z)=\sum_{m=0}^{\infty}\frac{q^{m(m-1)/2}}%
{\prod_{j=1}^{m}(1-q^{j})}\,z^{m} \label{prop_Gausspol_infty}%
\end{equation}
This relation immediately leads to two formulas due to Euler \cite[sec.
19.5]{Hardy_Wright}. There is wealth of nice properties of the product
$\left(  a;q\right)  _{n}=\prod_{k=0}^{n-1}\left(  1-aq^{k}\right)  $. The
finite product can be treated by the infinite variant as
\begin{align*}
\left(  a;q\right)  _{n}  &  =\prod_{k=0}^{n-1}\left(  1-aq^{k}\right)
=\frac{\prod_{k=0}^{n-1}\left(  1-aq^{k}\right)  \prod_{k=n}^{\infty}\left(
1-aq^{k}\right)  }{\prod_{k=n}^{\infty}\left(  1-aq^{k}\right)  }\\
&  =\frac{\prod_{k=0}^{\infty}\left(  1-aq^{k}\right)  }{\prod_{k=0}^{\infty
}\left(  1-aq^{k+n}\right)  }=\frac{\left(  a;q\right)  _{\infty}}{\left(
aq^{n};q\right)  _{\infty}}%
\end{align*}
We add the important relation of Gaussian polynomials
\begin{equation}
\prod_{j=1}^{n}(1+z\,q^{2j-1})(1+z^{-1}\,q^{2j-1})=\sum_{k=-n}^{n}\left[
\begin{array}
[c]{c}%
2n\\
n+k
\end{array}
\right]  (q^{2})\;q^{k^{2}}\,z^{k} \label{finite_jacobi}%
\end{equation}
which is a finite version of famous Jacobi's triple product\footnote{There
seems to exist a large variety of such identities. Another is Watson's
quintuple product
\begin{equation}
\sum_{k=-\infty}^{\infty}(z^{3k}-z^{-3k-1})\,q^{k(3k+2)/2}=\prod_{m=1}%
^{\infty}\left(  1-q^{m}\right)  \,\left(  1-z\,q^{m}\right)  \,\left(
1-z^{-1}\,q^{m-1}\right)  \,\left(  1-z^{2}\,q^{2m-1}\right)  \,\left(
1-z^{-2}\,q^{2m-1}\right)  \label{Watson_quintuple}%
\end{equation}
} \cite[sec. 19.8]{Hardy_Wright},
\begin{equation}
\theta_{3}(z,q)=\sum_{k=-\infty}^{\infty}z^{k}\,q^{k^{2}}=\prod_{m=1}^{\infty
}\left(  1-q^{2m}\right)  \,\left(  1+z\,q^{2m-1}\right)  \,\left(
1+z^{-1}\,q^{2m-1}\right)  \label{theta3}%
\end{equation}
that illustrates the intimate connection of the product $\left(  a;q\right)
_{n}=\prod_{k=0}^{n-1}\left(  1-aq^{k}\right)  $ and elliptic functions!

We must refer to the large literature (see e.g.
\cite{Goulden_Jackson,Rademacher,Andrews_q_series,Andrews_Partitions}) on the
product $\left(  a;q\right)  _{n}=\prod_{k=0}^{n-1}\left(  1-aq^{k}\right)  $.

\subsection{Examples of $q$-extension}

\subsubsection{The $q$-Gamma function}

The $q$-Gamma function is defined by%
\begin{equation}
\Gamma_{q}\left(  x\right)  =\left(  1-q\right)  ^{1-x}\frac{\left(
q;q\right)  _{\infty}}{\left(  q^{x};q\right)  _{\infty}}=\left(  1-q\right)
^{1-x}\frac{\prod_{k=0}^{\infty}\left(  1-q^{k+1}\right)  }{\prod
_{k=0}^{\infty}\left(  1-q^{k+x}\right)  }=\left(  1-q\right)  ^{1-x}%
\prod_{k=0}^{\infty}\frac{\left(  1-q^{k+1}\right)  }{\left(  1-q^{k+x}%
\right)  } \label{def_q-Gamma}%
\end{equation}
Explicitly,%
\[
\Gamma_{q}\left(  x\right)  =\left(  1-q\right)  ^{1-x}\lim_{m\rightarrow
\infty}\frac{\left(  1-q\right)  }{\left(  1-q^{x}\right)  }\frac{\left(
1-q^{2}\right)  }{\left(  1-q^{1+x}\right)  }\frac{\left(  1-q^{3}\right)
}{\left(  1-q^{2+x}\right)  }\cdots\frac{\left(  1-q^{m+1}\right)  }{\left(
1-q^{m+x}\right)  }%
\]
If $x=n\in\mathbb{N}$, then%
\begin{align*}
\Gamma_{q}\left(  n\right)   &  =\left(  1-q\right)  ^{1-n}\frac{\prod
_{k=0}^{\infty}\left(  1-q^{k+1}\right)  }{\prod_{k=0}^{\infty}\left(
1-q^{k+n}\right)  }=\left(  1-q\right)  ^{1-n}\frac{\prod_{k=1}^{\infty
}\left(  1-q^{k}\right)  }{\prod_{k=n}^{\infty}\left(  1-q^{k}\right)  }\\
&  =\left(  1-q\right)  ^{1-n}\prod_{k=1}^{n-1}\left(  1-q^{k}\right)
=\prod_{k=1}^{n-1}\frac{1-q^{k}}{1-q}%
\end{align*}
and%
\[
\lim_{q\rightarrow1}\Gamma_{q}\left(  n\right)  =\prod_{k=1}^{n-1}%
\frac{1-q^{k}}{1-q}=\prod_{k=1}^{n-1}\lim_{q\rightarrow1}\frac{1-q^{k}}%
{1-q}=\prod_{k=1}^{n-1}k=\left(  n-1\right)  !=\Gamma\left(  n\right)
\]

The $q$-Gamma function satisfies the functional equation%
\begin{equation}
\Gamma_{q}\left(  x+1\right)  =\frac{1-q^{x}}{1-q}\Gamma_{q}\left(  x\right)
\label{func_eq_q-Gamma}%
\end{equation}
Indeed, the definition (\ref{def_q-Gamma}) indicates that
\begin{align*}
\frac{\Gamma_{q}\left(  x+1\right)  }{\Gamma_{q}\left(  x\right)  }  &
=\frac{1}{1-q}\prod_{k=0}^{\infty}\frac{\left(  1-q^{k+x}\right)  }{\left(
1-q^{k+x+1}\right)  }=\frac{1}{1-q}\frac{\prod_{k=0}^{\infty}\left(
1-q^{k+x}\right)  }{\prod_{k=0}^{\infty}\left(  1-q^{k+x+1}\right)  }\\
&  =\frac{1}{1-q}\frac{\prod_{k=0}^{\infty}\left(  1-q^{k+x}\right)  }%
{\prod_{k=1}^{\infty}\left(  1-q^{k+x}\right)  }=\frac{1}{1-q}\left(
1-q^{k+x}\right)
\end{align*}
which is the functional equation (\ref{func_eq_q-Gamma}). After the limit
$q\rightarrow1$ in the functional equation (\ref{func_eq_q-Gamma}), we arrive,
with $\Gamma_{1}\left(  x\right)  =\Gamma\left(  x\right)  $, at the
functional equation $\Gamma\left(  x+1\right)  =x\Gamma\left(  x\right)  $ of
the Gamma function.

We follow Gosper's trick in \cite[Appendix A]{Andrews_q_series}, who was
inspired by Gauss' deduction of the Euler-Gauss product of the Gamma function
\cite[Part II]{PVM_Mittag-Leffler_Gamma}%
\begin{equation}
\Gamma\left(  z+1\right)  =\prod_{n=1}^{\infty}\left(  1+\frac{1}{n}\right)
^{z}\left(  1+\frac{z}{n}\right)  ^{-1} \label{Gauss_product_Gamma_function}%
\end{equation}
in inserting%
\[
\left(  1-q\right)  ^{x}=\frac{\prod_{k=0}^{\infty}\left(  1-q^{k+1}\right)
^{x}}{\prod_{k=1}^{\infty}\left(  1-q^{k+1}\right)  ^{x}}=\prod_{k=1}^{\infty
}\left(  \frac{1-q^{k}}{1-q^{k+1}}\right)  ^{x}%
\]
in (\ref{def_q-Gamma}) to obtain%
\begin{align*}
\Gamma_{q}\left(  x+1\right)   &  =\left(  1-q\right)  ^{-x}\prod
_{k=0}^{\infty}\frac{\left(  1-q^{k+1}\right)  }{\left(  1-q^{k+x+1}\right)
}=\prod_{k=1}^{\infty}\left(  \frac{1-q^{k+1}}{1-q^{k}}\right)  ^{x}%
\prod_{k=1}^{\infty}\frac{\left(  1-q^{k}\right)  }{\left(  1-q^{k+x}\right)
}\\
&  =\prod_{k=1}^{\infty}\left(  \frac{1-q^{k}+q^{k}-q^{k+1}}{1-q^{k}}\right)
^{x}\left(  \frac{1-q^{k}}{1-q^{k+x}}\right)
\end{align*}
Further, with $\frac{1-q^{k+1}}{1-q^{k}}=\frac{1-q^{k}+q^{k}-q^{k+1}}{1-q^{k}%
}=1+q^{k}\frac{1-q}{1-q^{k}}$ and $\frac{1-q^{k}}{1-q^{k+x}}=\left(
\frac{1-q^{k+x}}{1-q^{k}}\right)  ^{-1}=\left(  \frac{1-q^{k}+q^{k}-q^{k+x}%
}{1-q^{k}}\right)  ^{-1}=\left(  1+q^{k}\frac{1-q^{x}}{1-q^{k}}\right)  ^{-1}%
$, we find the $q$-extension of the Euler-Gauss product,%
\begin{equation}
\Gamma_{q}\left(  x+1\right)  =\prod_{k=1}^{\infty}\left(  1+q^{k}\frac
{1-q}{1-q^{k}}\right)  ^{x}\left(  1+q^{k}\frac{1-q^{x}}{1-q^{k}}\right)
^{-1} \label{Gauss_product_Gamma_function_q_extension}%
\end{equation}
Immediately, the limit $q\rightarrow1$ reduces to the Euler-Gauss product
(\ref{Gauss_product_Gamma_function}).

\subsubsection{The $q$-sinus}

We compute with (\ref{Gauss_product_Gamma_function_q_extension}) the product
\begin{align*}
\Gamma_{q}\left(  1+x\right)  \Gamma_{q}\left(  1-x\right)   &  =\prod
_{k=1}^{\infty}\left(  1+q^{k}\frac{1-q^{x}}{1-q^{k}}\right)  ^{-1}\left(
1+q^{k}\frac{1-q^{-x}}{1-q^{k}}\right)  ^{-1}\\
&  =\prod_{k=1}^{\infty}\frac{1-q^{k}}{1-q^{k+x}}\frac{1-q^{k}}{1-q^{k-x}%
}=\prod_{k=1}^{\infty}\frac{1-2q^{k}+q^{2k}}{1-q^{k}\left(  q^{x}%
+q^{-x}\right)  +q^{2k}}\\
&  =\prod_{k=1}^{\infty}\frac{1-2q^{k}+q^{2k}}{1-2q^{k}\cosh\left(  x\log
q\right)  +q^{2k}}%
\end{align*}
The infinite product of the sinus function \cite[4.3.89]{Abramowitz}%
\begin{equation}
\frac{\sin\left(  \pi z\right)  }{\pi z}=%
%TCIMACRO{\dprod \limits_{k=1}^{\infty}}%
%BeginExpansion
{\displaystyle\prod\limits_{k=1}^{\infty}}
%EndExpansion
\left(  1-\frac{z^{2}}{k^{2}}\right)  \label{infinite_product_sinz}%
\end{equation}
then suggests that%
\[
\lim_{q\rightarrow1}\prod_{k=1}^{\infty}\frac{1-q^{k}}{1-q^{k+x}}\frac
{1-q^{k}}{1-q^{k-x}}=\prod_{k=1}^{\infty}\frac{k^{2}}{\left(  k+x\right)
\left(  k-x\right)  }=\prod_{k=1}^{\infty}\frac{1}{\left(  1-\frac{x^{2}%
}{k^{2}}\right)  }=\frac{\pi x}{\sin\left(  \pi x\right)  }%
\]
Since $\Gamma_{q}\left(  1+x\right)  \Gamma_{q}\left(  1-x\right)
=\frac{1-q^{x}}{1-q}\Gamma_{q}\left(  x\right)  \Gamma_{q}\left(  1-x\right)
$ is, by $q$-definition, equal to $\Gamma_{q}\left(  x\right)  \Gamma
_{q}\left(  1-x\right)  =\frac{\pi}{\sin_{q}\left(  \pi x\right)  }$, so that
we may conclude that%
\[
\Gamma_{q}\left(  1+x\right)  \Gamma_{q}\left(  1-x\right)  =\frac{1-q^{x}%
}{1-q}\Gamma_{q}\left(  x\right)  \Gamma_{q}\left(  1-x\right)  =\frac
{1-q^{x}}{1-q}\frac{\pi}{\sin_{q}\left(  \pi x\right)  }=\prod_{k=1}^{\infty
}\frac{1-q^{k}}{1-q^{k+x}}\frac{1-q^{k}}{1-q^{k-x}}%
\]
Hence, the $q$-extension of the sinus product (\ref{infinite_product_sinz}) is%
\[
\sin_{q}\left(  \pi x\right)  =\pi\frac{1-q^{x}}{1-q}\prod_{k=1}^{\infty}%
\frac{1-q^{k+x}}{1-q^{k}}\frac{1-q^{k-x}}{1-q^{k}}%
\]

\end{document}